\documentclass[10pt,a4paper]{amsart}
\usepackage{color,graphicx}

\title[]{Optimal Regularity for The Signorini Problem and its Free Boundary.}
\author[J. Andersson]{John Andersson}
\address{Mathematics Institute,
University of Warwick
Coventry CV4 7AL, UK}
\email{j.e.andersson@warwick.ac.uk}

\thanks{$2000$ {\it Mathematics Subject Classification.\/} Primary
35R35, Secondary 35B40, 35J60.}
\thanks{\textsl{Thanks:} }
\date{}

\def\dh{{\> d\mathcal{H}^{n-1}}}
\def\dh1{{\> d\mathcal{H}^1}}
\def\dh2{{\> d\mathcal{H}^2}}

\def\R{{\mathbb{R}}}
\def\u{{\textbf{u}}}
\def\v{\textbf{v}}
\def\w{\textbf{w}}
\def\p{\textbf{p}}
\def\f{\textbf{f}}
\def\q{\textbf{q}}

\def\div{\textrm{div}}

\def\Pr{\mathbf{Pr}}
\newtheorem{lem}{Lemma}
\newtheorem{prop}{Proposition}
\newtheorem{cor}{Corollary}
\newtheorem{thm}{Theorem}
\newtheorem{defi}{Definition}

\begin{document}

\begin{abstract}
We will show optimal regularity for minimizers of the Signorini problem for 
the Lame system. In particular if $\u=(u^1,u^2,u^3)\in W^{1,2}(B_1^+:\R^3)$ 
minimizes
$$
J(\u)=\int_{B_1^+}|\nabla \u+\nabla^\bot \u|^2+\lambda\div(\u)^2
$$
in the convex set
$$
K=\big\{ \u=(u^1,u^2,u^3)\in W^{1,2}(B_1^+:\R^3);\; u^3\ge 0 \textrm{ on }\Pi,
$$
$$
\u=f\in C^\infty(\partial B_1) \textrm{ on }(\partial B_1)^+ \big\},
$$
where $\lambda\ge 0$ say.

Then $\u\in C^{1,1/2}(B_{1/2}^+)$. Moreover the free boundary, given by
$$
\Gamma_\u=\partial \{x;\;u^3(x)=0,\; x_3=0\}\cap B_{1},
$$
will be a $C^{1,\alpha}$ graph close to points where $\u$ is not degenerate.

Similar results have been know before for scalar partial differential equations
(see for instance \cite{AC} and \cite{ACS}).
The novelty of this approach is that it does not rely on maximum principle
methods and is therefore applicable to systems of equations. 
\end{abstract}
\maketitle
\tableofcontents
\section{Introduction}

We are interested in minimizers $\u=(u^1,u^2,u^3)$ in the set
$$
K=\big\{ \u=(u^1,u^2,u^3)\in W^{1,2}(B_1^+:\R^3);\; u^3\ge 0 \textrm{ on }\Pi,
$$
$$
\u=f\in C^\infty(\partial B_1) \textrm{ on }(\partial B_1)^+ \big\}.
$$
of the following functional
\begin{equation}\label{signorini}
J(\u)=\int_{B_1^+}|\nabla \u+\nabla^\bot \u|^2+\lambda\div(\u)^2.
\end{equation}
Here $\nabla^\bot \u=(\nabla \u)^T$ is the transpose of the gradient matrix.
We will always assume that $\lambda\ge 0$. 
We could relax the condition that $\u\in W^{1,2}$ to $\nabla \u+\nabla^\bot \u\in L^2$, as is usually done. But due to Korn's inequality both conditions are 
equivalent. 

If we denote $\Pi=\{x;\; x_3=0\}$ and  $\Lambda_\u=\{x\in \Pi;\; u^3(x)=0 \}$ 
then it is easy to see that the minimizers solves the following Euler-Lagrange 
equations
\begin{equation}\label{Lame}
\begin{array}{ll}
L\u\equiv \Delta \u+\frac{2+\lambda}{2} \nabla \textrm{div}(\u)=0 & \textrm{ in } \R^{3}_+
\end{array}
\end{equation}
\begin{equation}\label{Bdry1}
\begin{array}{ll}
u^3=0 & \textrm{ on }\Lambda_\u
\end{array}
\end{equation}
\begin{equation}\label{Bdry2}
\begin{array}{ll}
\frac{\partial u^3}{\partial x_3}+\frac{\lambda}{4} \textrm{div}(\u)=0 & \textrm{ on }
\Pi\setminus\Lambda_\u
\end{array}
\end{equation}
\begin{equation}\label{Bdry3}
\begin{array}{ll}
\frac{\partial u^i}{\partial x_3}+\frac{\partial u^3}{\partial x_i}=0 & \textrm{ on } \Pi \textrm{ for }i=1,2
\end{array}
\end{equation}
\begin{equation}\label{Bdry4}
\begin{array}{ll}
u^3\ge 0 & \textrm{ on } \Pi 
\end{array}
\end{equation}
\begin{equation}\label{Bdry5}
\begin{array}{ll}
\frac{\partial u^3}{\partial x_3}+\frac{\lambda}{4}\div(\u)\le 0 & \textrm{ on } \Pi.
\end{array}
\end{equation}

It is important to note that this problem is highly nonlinear since the
set $\Lambda_\u$ is not apriori known. The major difficulty in analysing
the regularity of this problem consists in understanding not only the
behavior of the solution but also of the unknown set $\Lambda_\u$.

This minimization problem models the deformation of an elastic body,
which we here assume for simplicity to be the half ball $B_1^+$ when it is
subjected to some deformation $f$ of the curved part of the boundary 
$\partial B_1^+$ and is required to stay above a certain obstacle,
here given by $x_3=0$.

This is of course a version of the Signorini problem. This problem was 
first formulated by
Antonio Signorini in 1933 \cite{Signorini}.
In the oringinal formulation of the problem Signorini assumed Neumann data on 
the boundary and he included the influence of gravity. From
a mathematical point of view adding gravity to the functional $J(\u)$ does not 
result in any new difficulties. Our analysis is entirely local so the boundary 
data on $(\partial B_1)^+$ will play no roll in our analysis.
Signorini where interested in the existence
and uniqueness of solutions. This was solved by G. Fichera \cite{Fichera} in
1963. With the advances in the calculus of variations since the sixties 
the existence and uniqueness is today considered to be quite standard.

Here we are interested in the regularity of minimizers and in the regularity
of the free boundary $\partial \Lambda_\u$. 

\noindent {\bf Mathematical Background: } It took almost 50 years from 
Signorini's formulation of this problem to the first regularity results
was proved by D. Kinderlehrer in 1981 \cite{kinderlehrer}. Kinderlehrer
proves that the solution is $C^{1,\beta}$ in the case $n=2$.

Soon after A.A. Arkhipova and N.N. Uraltseva showed $C^{1,\beta}$
regularity for variational inequalities of diagonal systems in $n$ 
dimensions \cite{arkhipovauraltseva}. The assumption that the system is 
diagonal excludes 
the Signorini problem from their result. The first $C^{1,\beta}$
result for the Signorini problem in general dimensions is due to
R. Schumann who proved $C^{1,\beta}$-regularity for some $\beta>0$ in 1989
\cite{schumann}.

There are several other papers relating to free boundary problems
for systems of equations see for instance M. Fuchs \cite{fuchs} for a 
pleasant proof of regularity and free boundary regularity for a system.

However, all previous proofs of optimal regularity and free boundary regularity
results for systems of equations are based on the reduction of the
system to a scalar problem. To the authors knowledge there are no 
papers that manage to tackle the difficulties of systems without such a
reduction. Let us therefore investigate the development of the regularity
theory for the scalar versions of the Signorini problem - where much more
is known.

There has been significant progress
in the understanding of the regularity questions for the scalar Signorini
problem, also called the thin obstacle problem, in the last decade, 
see \cite{AC} and \cite{ACS}. The thin obstacle problem is the 
minimisation of the Dirichlet energy
\begin{equation}\label{DirichletEnergy}
\int_{B_1^+}|\nabla u|^2
\end{equation}
in the set
\begin{equation}\label{MinSetinACS}
K=\{u\in W^{1,2};\; u\ge 0  \textrm{ on }\Pi,\; u=f \textrm{ on }(\partial B_1)^+\}.
\end{equation}
To show existence of minimizers to the thin obstacle problem is again rather 
standard. But the regularity theory is quite subtle, both with respect to
the regularity of the solution \cite{AC} and its free boundary \cite{ACS}.

We have several good reasons to dwell on the technique used in \cite{AC} 
and \cite{ACS}. First of all, those papers have provided the framework 
for this paper even though the techniques we use will be very different
from theirs. Secondly, to know something about the scalar problem
considered by Athanasopoulos, Caffarelli and Salsa will also help us understand
the difficulties of the vectoral case. In particular we will
be able to understand why Signorini's problem have not been solved
by the techniques developed in the last thirty years.

In \cite{AC} the main result if that minimizers of (\ref{MinSetinACS}) in the 
set (\ref{MinSetinACS}) are $C^{1,1/2}$, which is the optimal regularity. The 
proof is based on the Bernstein technique, a monotonicity formula and an 
iteration. The Bernstein technique is basically to apply the maximum
principle to the function
$$
g(x)=\eta(x)\frac{\partial^2 u}{\partial e_i^2}-\lambda |\nabla u|^2
$$
where $\eta$ is a cut off function and $e_i\in \Pi$. Since the 
maximum principle is not true for the Lame system we can not replicate this 
argument for the vectoral Signorini problem. Neither do the structure of
the Lame system allow us to derive the monotonicity formula that is essential
for the optimal regularity proof.

In \cite{ACS} the authors use comparison and boundary comparison principles 
together with some geometrical insight to show that the free boundary
$\partial \{ u>0\}\cap \Pi$ is $C^{1,\alpha}$ in a small neighbourhood 
around free boundary points $x^0$ where $\sup_{B_r(x^0)}|u|\approx r^{1+\beta}$
for any $\beta<1$. The usage of comparison principles makes it impossible
to apply their technique directly to solutions of the vectoral Signorini 
problem.

The theory in \cite{AC} and \cite{ACS} was later generalized in
\cite{CSS} to more general thin obstacle problems interpreted as obstacle 
problems for the fractional Laplacian. But the methods in \cite{CSS} are 
quite similar to the methods used in \cite{AC} and \cite{ACS}. In particular
the methodology in \cite{CSS} relies heavily on comparison and
maximum principles and is therefore not applicable for our problem.

One can say that this article constitutes the author's attempt to
develop a regularity theory for free boundaries that is not dependent on 
maximum principles.

Instead of maximum principle methods we will rely on the blow-up method,
the Liouiville Theorem and linearization together with some simple
geometric observations and a nice way to control blow-up sequences
that we get from \cite{ASU}. 

The article naturally divides itself into two parts that depend
on different procedures. The first part, after some intorductionary and 
standard considerations, constitutes of section
\ref{GlobalSolPart1}-\ref{sectionFlatness} where we show that 
close to free boundary points where $\u$ growths like $r^{1+\alpha}$ for 
some $\alpha<1$ the free boundary $\Gamma_\u$ is actually flat. That 
$\Gamma_\u$ is flat just means that in a small ball $B_r\cap \Pi$ the
free boundary is contained in a strip of with $\sigma(r)r$ for some
modulus of continuity $\sigma$.

The idea of the proof is quite straightforward and uses a result by
M. Benedicks \cite{Benedicks} that states that the set of positive
harmonic functions vanishing on part of $\Pi$ and has zero Neumann data
on the rest of $\Pi$ is one dimensional. Using this result we may 
deduce that the tangential derivatives of the blow-up of a solution $\u$ are
all multiples of each other. That implies in particular that
the blow-up of $\u$ depend only on two directions, say $x_1$ and $x_n$. 
By means of the Liouiville Theorem we can classify such solutions and thus
calculate the asymptotic profile of solutions close to
points where $\sup_{B_r}|\u|\approx r^{1+\alpha}$. The profile in question
is explicitly calculated in Lemma \ref{eigenfunctionsinR2}. In
polar coordinates the asymptotic profile is
\begin{equation}\label{assprofile}
\begin{array}{l}
u(r,\phi)=r^{3/2}\Big(\frac{18+3\lambda}{40}\cos\big((5/2)\phi \big)-\frac{2+\lambda}{8}\cos\big((1/2)\phi \big)\Big)\\
v(r,\phi)=r^{3/2}\Big( \frac{6-\lambda}{40}\sin\big((5/2)\phi \big) 
+\frac{2+\lambda}{8}\sin\big((1/2)\phi \big)  \Big).
\end{array}
\end{equation}
Since the growth of the asymptotic profile is $r^{3/2}$ we can directly 
conclude that $\u\in C^{1,\beta}(B_{1/2}^+)$ for each $\beta<1/2$, see Lemma \ref{almostoptholderreg} and Corollary \ref{addeddetailin6}.

In order to use Benedicks result we will derive that global solutions
with control of the growth at infinity is actually determined by
a harmonic function. Later, in appendix 2,
we will also use this method to indicate
how to make a simple eigenfunction expansion of the linearized problem.
It is quite possible that one could derive all the regularity
theory for the vectoral Signorini problem by this reduction to 
harmonic functions. We will however only use this reduction in our
proof that the free boundary is flat, see the proof of Corollary 
\ref{taubisiszeroonOmega}, and in appendix 2. We believe that the
result in the appendix is well known and that it could be proved by 
other methods such as spectral theory of operator pencils 
\cite{KozlovMazyaRossman}. I could unfortunately not find any good reference 
to such a result. And the machinery of operator pencils 
\cite{KozlovMazyaRossman} is too heavy to 
introduce in this paper to prove a supporting lemma. Therefore, for convenience,
we use the harmonic reduction again in appendix 2.

The first part of the proof is quite trivial from a technical point of view and
we use mostly standard calculus and elementary pde theory. The result
is however very important for the linearisation that follows.
In in section \ref{GlobalSolPart1}-\ref{sectionFlatness} we show 
that the asymptotics of the solution is uniquely determined
at points of lowest regularity. This allows us to make a
linearisation at all such points which will imply everywhere regularity.
This is quite different, and much stronger, than the standard outcome
of a linearisation argument where an $\epsilon-$closeness assumption
is needed and only partial regularity (which may or may not be optimal) 
can be deduced.

The second part of the paper is far more technical and, unfortunately,
much harder to read. There we prove that the solutions are in fact $C^{1,1/2}$
which is optimal as the above asymptotic profile demonstrates. We also 
prove that the free boundary is $C^{1,\alpha}$ close to points where 
$\u$ has the above asymptotic profile (this includes all the points 
$x^0\in \Gamma_\u$ 
where $\sup_{B_r^+}|\u|> r^{1+\beta}$ for  some $\beta<1$ and $r$ small).

The argument is by linearization and flatness improvement. In particular
if the origin is a free boundary point of $\u$ with the asymptotic profile 
$\p$ where $\p=(u,v)$ as in equation (\ref{assprofile}). Then, heuristically at least, the limit (we will use a slightly different limit later)
\begin{equation}\label{meaningless}
\lim_{r\to 0}\frac{\u(rx)-\p(rx)}{\|\u(rx)-\p(rx)\|_{L^2(B_1^+)}}=\v
\end{equation}
will contain information about the regularity of $\u$ and the free boundary.
The problem is that in order to extract any useful information from 
(\ref{meaningless}) we need the convergence to be strong in $L^2$ we would also 
need to show that $\v$ is better than $\u$. This is a very delicate matter
that will be analysed in sections \ref{fundsec}-\ref{decaysec}.

In the final sections we prove the regularity theorem and free boundary
regularity. We also show, Lemma \ref{Claim1Theorem1}, that the analysis
in the previous sections can be made uniform.

Intuitively there is a gap in 
the eigenvalues of the operator for the Lame system on the sphere for the
linearized problem. Where
the next homogeneous solution after $\p$ as in (\ref{assprofile}) is
homogeneous of second order. That implies that $\sup_{B_r^+}|\u-\p|$ will be of
order $r^2$ which implies that the difference between $\u$ and $\p$
decays geometrically in smaller and smaller balls. This is enough to 
deduce the regularity of the solution and the free boundary.

Throughout this paper we will not use the maximum principle or
comparison principles at all.

There is however a deeper connection with the flatness proof of this paper and 
the methodology used in \cite{AC}, \cite{ACS} and \cite{CSS}. As mentioned
before, the flatness proof we use is based on a result by Benedicks 
\cite{Benedicks}. Both Benedicks paper and the monotonicity formula used in
\cite{AC} is based on a paper by Friedland and Hayman \cite{FriedlandHayman}.
Friedland and Hayman's paper also use some structural properties 
of the Laplacian that goes back at least to Alfred von Huber \cite{vonHuber}.
The same methods is used to derive a frequency formula in \cite{CSS}.

Even though there are some technical connections between this paper
and previous work on thin obstacle problems. Most of the material here is 
essentially new. It is the authors hope that the linearisation technique
will prove useful also for other free boundary problems involving
systems of equations.

Our main regularity result is.

$ $

\noindent\textbf{Main Regularity Theorem.} {\sl
Let $\u$ be a solution to the Signorini problem in $B_1^+$ then 
$$
\|\u\|_{C^{1,1/2}(B_{1/2}^+)}\le C\|\u\|_{L^2(B_1^+)}
$$}

$ $

The free boundary regularity result is the same as in \cite{ACS} but we will 
need to introduce some notation before we can state the Theorem. The precise 
formulation can be found in Corollary \ref{freebdryregcor} in 
section \ref{freeboundregsec}.

The regularity theorems are, for simplicity, only formulated for solutions 
in $B_1^+$ with the constraint $u^3\ge 0$. The more general 
problem to minimize
$$
\int_{D}|\nabla \u+\nabla^\bot \u|^2+\lambda\div(\u)^2
$$
in the set
$$
K=\big\{ \u=(u^1,u^2,u^3)\in W^{1,2}(D:\R^3);\; u^3(x)\ge \psi(u^1(x),u^2(x)),
$$
$$
\u \textrm{ satisfies appropriate boundary conditions}\big\},
$$
where $D$ is some $C^{1,\beta}$ domain and $\psi$ is a $C^{1,\beta}$ function
with $\beta>1/2$, can be handled by a perturbation argument as in
for instance \cite{ASW}.

$ $

\noindent\textbf{Notation:}

\noindent At times we will write $n$ for the dimension. But some proofs
in the paper has been written only for $n=3$. This is for simplicity,
since the curl operator is more explicit in $\R^3$. The pedantic (here used
with no negative connotation) reader can always think that $n=3$.

\noindent $\Pi=\{x;\; x_n=0\}$ is the boundary of $\R^n_+$ .

\noindent We will use bold face $\u$, $\v$, $\w$, $\p$ etc. to denote 
vector valued functions $\u=(u^1,u^2,u^3, ..., u^n)$, $\v=(v^1,v^2,..., v^n)$ etc.

\noindent For a continuous function $\u=(u^1,u^2,..., u^n)$ we set 
$\Lambda_\u=\Lambda=\{x;\; x_n=0,\; u^n(x)=0\}$.

\noindent For a continuous function $\u=(u^1,u^2,...,u^n)$ we set 
$\Omega_\u=\Omega=\Pi\setminus \Lambda_\u$.

\noindent For a continuous function $\u$ we define its free boundary 
$\Gamma_\u=\Gamma=\overline{\Lambda_\u}\cap \overline{\Omega_\u}$.

\noindent $\nabla \u$ is the matrix:
$$
\nabla \u=
\left[
\begin{array}{llll}
\frac{\partial u^1}{\partial x_1} & \frac{\partial u^1}{\partial x_2} &
\cdots &\frac{\partial u^1}{\partial x_n} \\
\frac{\partial u^2}{\partial x_1} & \frac{\partial u^2}{\partial x_2} &
\cdots & \frac{\partial u^2}{\partial x_n} \\
\vdots & \vdots & & \vdots \\
\frac{\partial u^n}{\partial x_1} & \frac{\partial u^n}{\partial x_2} &
\cdots & \frac{\partial u^n}{\partial x_n} 
\end{array}
\right].
$$

\noindent $\nabla^\bot \u$ is the transpose of the matrix $\nabla \u$:
$$
\nabla^\bot \u=
\left[
\begin{array}{llll}
\frac{\partial u^1}{\partial x_1} & \frac{\partial u^2}{\partial x_1} &
\cdots & \frac{\partial u^n}{\partial x_1} \\
\frac{\partial u^1}{\partial x_2} & \frac{\partial u^2}{\partial x_2} &
\cdots & \frac{\partial u^n}{\partial x_2} \\
\vdots & \vdots & & \vdots \\
\frac{\partial u^1}{\partial x_n} & \frac{\partial u^2}{\partial x_n} &
\cdots & \frac{\partial u^n}{\partial x_n} 
\end{array}
\right].
$$

\noindent We will often use a prime to indicate the projection of an 
$n$-dimensional vector into an $(n-1)$-dimensional vector: $x'=(x_1,x_2,\cdots, x_{n-1})$,
$\nabla'=\big(\frac{\partial}{\partial x_1},\frac{\partial}{\partial x_1}, \cdots, \frac{\partial }{\partial_{n-1}}  \big)$ etc. At times we will slightly abuse notation and write
$x'=(x_1,x_2,\cdots, x_{n-1},0)$ and $\nabla'=\big(\frac{\partial}{\partial x_1},\frac{\partial}{\partial x_2}, \cdots, \frac{\partial}{\partial x_{n-1}},0  \big)$. It will always be clear from context what we intend. 

\noindent We use the notation $\nabla''=\big(0,\frac{\partial}{\partial x_2},\frac{\partial}{\partial x_3},...,\frac{\partial}{\partial x_{n-1}} ,0\big)$.
More generally for any vector $\xi\in \Pi$ we use the notation
$\nabla_\xi''=\nabla -e_n\frac{\partial }{\partial x_n}-\xi_w/|\xi_w|^2
(\xi_w\cdot \nabla )$ which is just the gradient restricted to the subspace 
orthogonal to $e_n$ and $\xi$.

\noindent $\Pr(\u,r)$ is a projection operator onto affine functions defined in
Definition \ref{PrDef}.

\noindent By $W^{k,p}(D)$ we mean the normal Sobolev space. We will often be
quite informal when assigning vector valued functions to this space and
write $(u^1,u^2,u^3, \cdots, u^n)\in W^{k,p}$ instead of $(u^1,u^2,u^3, \cdots, u^n)\in (W^{k,p})^n$ etc.

\noindent By $\|\u\|_{\tilde{L}^2(\Omega)}$ we mean the norm:
$$
\|\u\|_{\tilde{L}^p(\Omega)}=\bigg( \frac{1}{|\Omega|}\int_{\Omega}|\u|^p\bigg)^{1/p},
$$
defined in Definition \ref{definitiontildeL2}.

\noindent The homogeneous solutions to the Lame system $\p_{1/2}$, $\p_{3/2}$,
$\p^\xi_{1/2}$, $\p_{3/2}^\xi$, $\tilde{\p}_{3/2}$,... are defined in
Definition \ref{2dsolutiondef}.

\section{Some Simplifying Conventions.}\label{secsimpconv}

Our main goal is to prove that the solutions are $C^{1,1/2}$. It
will however simplify our exposition significantly if we assume 
that we have $C^{1,\beta}$ regularity. The techniques developed
in the following pages is certainly strong enough to prove that 
solutions are $C^{1,\beta}$. But to formally prove the $C^{1,\beta}$ regularity
we would have to work through
the same string of Lemmas and Theorems twice and end up with an article 
twice as long. We will therefore assume that the solutions are $C^{1,\beta}$ 
without proof. But we have indicated
in an appendix how to prove the following lemma. The exposition in the appendix
is rather terse and we will freely refer to results proved in the main body
of the paper. Hopefully there is enough information in the appendix for a 
thorough reader to reconstruct the proof.

Another $C^{1,\beta}$ proof was published in \cite{schumann}. I have not
been able to verify that proof due to a lack knowledge of psedodifferential 
operators.

\begin{lem}\label{c1beta}
Let $\u$ be a solution to the Signorini problem
$$
\begin{array}{ll}
\Delta \u+\frac{2+\lambda}{2} \nabla \textrm{div}(\u)=0 & \textrm{ in } B_1^+ \\
u^3=0 & \textrm{ on }\Lambda \\
\frac{\partial u^3}{\partial x_3}+\frac{\lambda}{4} \textrm{div}(\u)=0 & \textrm{ on }
\Pi\setminus\Lambda \\
\frac{\partial u^i}{\partial x_3}+\frac{\partial u^3}{\partial x_i}=0 & \textrm{ on } \Pi \textrm{ for }i=1,2 \\
u^3\ge 0 & \textrm{ on } \Pi  \\
\frac{\partial u^3}{\partial x_3}+\frac{\lambda}{4}\div(\u)\le 0 & \textrm{ on } \Pi.
\end{array}
$$
Then 
$\u\in C^{1,\beta}(B_{1/2}^+)$ 
for some $\beta>0$ and $\u$ satisfies the following estimate
$$
\|\u\|_{C^{1,\beta}(B_{1/2}^+)}\le C\|u\|_{L^{2}(B_1^+)}.
$$
\end{lem}

One of the advantages to have $C^{1,\beta}$ regularity in what follows is that
it will significantly simplify our exposition. It is easy to verify that if
$\u$ is a solution to the Signorini problem then
$$
\u+
\left[\begin{array}{l}
b_{1} \\
b_{2} \\
0
\end{array}\right]
+
\left[\begin{array}{lll}
a_{11} & a_{12} & 0 \\
a_{21} & a_{22} & 0 \\
0 & 0 & a_{33}
\end{array}\right]
\left[ \begin{array}{l}
x_1 \\ x_2 \\ x_3
\end{array}\right],
$$
is also a solution for any constants $b_1,b_2$ and  $a_{ij}$ such that
$$
a_{33}+\frac{\lambda}{4}(a_{11}+a_{22}+a_{33})=0.
$$
If $\u\in C^{1,\beta}$ then we can make the following, informal, 
standing assumption.

\noindent{\bf Standing Assumption:} {\sl Let $\u$ be a solution to the 
Signorini problem and assume that $x^0$ is a free boundary point of interest
(such as a point that we make a blow-up at). Then we will assume that
$|\u(x^0)|=0$ and that $\nabla \u(x^0)=0$.}

A more formal way of handling this would be to only consider 
our solutions modulo affine functions and define 
$$
\sup_{B_r(x^0)}|\u|\le C_0
$$ 
if there is a function $\v$ in the same equivalence class as $\u$ 
such that the estimate holds etc. 

In appendix 1 where we indicate how to prove the $C^{1,\beta}$-lemma 
we will not rely on the standard assumption but instead define a
projection operator $\Pr$ (see Definition \ref{PrDef}) and consider
$\u-\Pr(\u,r)$. In the main body of the paper this extra $\Pr(\u,r)-$term
would only clutter down our already complicated expressions too much. Therefore 
we will rely on the standing assumption.

At times we will refer to the Liouville Theorem without explanation. 
Whenever that is done we refer to the following simple result.

$ $

\noindent{\bf Liouville's Theorem.} {\sl Let $u$ be any function defined 
in $\R^n$ that satisfies the following estimates for all $k\in \mathbb{N}$, 
$R>1$ and some constant $C_0$
\begin{equation}\label{qoodetsforliou}
\sup_{B_R}|D^k u|\le \frac{C_0}{R^{k+n/2}}\|u\|_{L^2(B_{2R})}.
\end{equation}
Assume furthermore that $u$ satisfies the growth condition
$$
\|u\|_{L^2(B_{R})}\le C_1R^{k_0+n/2+\alpha}
$$
for all $R>1$, some $\alpha<1$, some $k_0\in \mathbb{N}$ and some constant 
$C_1$. Then $u$ is a polynomial of order $k_0$.}

$ $

The proof is trivial. This applies in particular to harmonic functions
of polynomial growth for which the estimates (\ref{qoodetsforliou}) are 
standard.

$ $

\textbf{Something about the dimension $n$.} All the results in this paper
are valid in $\R^n$ for any $n\ge 2$. However the main technical difficulties
arise in $\R^3$. Also, some proofs will rely on the curl operator that is
much more explicit in $\R^3$. Therefore some of the proofs are written
only for $n=3$. This is an attempt to balance the clarity of the
exposition without avoiding any of the intrinsic difficulties of the problem 
which arise in $\R^3$. 

In $\R^n$ we can define the curl operator on a vector field $\u$ according to
$$
\textrm{curl}(\u)=\big[ * (d\u^b)\big]^\sharp,
$$
where $*$ is the Hodge star, $^b$ the flat and $^\sharp$ the sharp operator.
With this definition we would still have curl$(\nabla f)=0$ etc. and all
the proofs would still work. Hopefully, the assumption that $n=3$
will increase the clarity enough to motivate the loss of generality.

In the later sections of the paper, where we do not use the curl operator, we 
will write $n$ instead of $3$. This is to indicate that the technique
is not simplified by the assumption $n=3$. The reader should always remember that
we, in order to avoid working with $[ * (d\u^b)]^\sharp$, 
always assume that $n=3$.

\section{Weak Regularity.}

In order to prove $W^{2,2}$ estimates for solutions to the Signorini problem
we need the well known Korn's inequality found for instance in \cite{Ciarlet}.

\begin{lem}\textsc{[Korn's Inequality]}
Let $\u: B_1^+\to \R^3$ then
$$
\bigg(\int_{B_1^+}|\nabla \u|^2\bigg)^{1/2}\le C\bigg(\int_{B_1^+}|\nabla \u+\nabla^\bot \u|^2\bigg)^{1/2}
$$
whenever the right hand side is defined. In particular 
$$
\bigg(\int_{B_1^+}|\nabla \u+\nabla^\bot \u|^2\bigg)^{1/2}
$$
is a semi norm on $W^{1,2}(B_1^+)$.
\end{lem}

Next we state a simple Lemma. The proof is standard and therefore omitted.
\begin{lem}\label{w12est}
Let $\u$ be a solution to the Signorini problem, that is $\u$ minimizes (\ref{signorini}) in $K$, then 
$$
\|\u\|_{W^{1,2}(B_{1/2}^+)}\le C\|\u\|_{L^2(B_{3/4}^+)}.
$$
\end{lem}

\begin{lem}
Let $\u$ be a solution to the Signorini problem, that is $\u$ minimizes (\ref{signorini}) in $K$,
then $\u\in W^{2,2}(B_{1/2}^+)$ and
$$
\|\u\|_{W^{2,2}(B_{1/2}^+)}\le C\|\u\|_{L^2(B_{3/4}^+)}.
$$
\end{lem}
\textsl{Proof:} Let $\xi\in C^\infty_0(B_{3/4})$ be a standard cut off function
$\xi=1$ in $B_{1/2}$, $\xi\ge 0$, $|\nabla \xi|\le 8$. Then
$$
\v=\v_{h,t}(x)=\u(x)+t\xi^2\frac{\u(x+e_ih)-\u(x)}{h}=\u(x)+t\xi^2\u_h
$$
is a competitor for minimality in $K$ if $0\le t \le h$ and $i=1,2$.
Thus
$$
0\le \int_{B_1^+}\nabla \u\cdot \nabla \v+\nabla^\bot \u\cdot \nabla \v+\nabla \u\cdot \nabla^\bot \v+\nabla^\bot \u\cdot \nabla^\bot \v+\lambda \div(\u)\div(\v).
$$
Differentiating at $t=0$ we get
$$
0\le \int_{B_1^+}\xi^2\Big( \nabla \u \cdot \nabla \u_h+\nabla^\bot \u\cdot\nabla \u_h+ \nabla \u\cdot \nabla^\bot \u_h+
$$
\begin{equation}\label{firstindiffquot}
\nabla^\bot \u\cdot \nabla^\bot \u_h +\lambda\div(\u)\div(\u_h)\Big)+
\end{equation}
$$
2\xi\Big(\u_h\cdot \nabla \u \cdot \nabla\xi+\u_h\nabla^\bot\u\cdot\nabla \xi+\nabla \xi\cdot \nabla \u\cdot \u_h+\nabla \xi \cdot \nabla^\bot \u\cdot \u_h+\lambda \div(\u)\nabla \xi \cdot \u_h \Big).
$$
If we denote $\u^h(x)=\u(x+e_ih)$ then $\u^h$ is a minimizer in $B_{1-h}^+$
and 
$$
\v_{h,t}=\v=\v^h+t\xi^2\v^h_{-h}
$$
is an admissible competitor for minimization. 
Therefore by differentiation at $t=0$
and using that $\u^h$ is a minimizer we can conclude that
\begin{equation}\label{secondindiffquot}
0\le \int_{B_{1-h}^+}\xi^2\Big( \nabla \u^h \cdot \nabla \u^h_{-h}+\nabla^\bot \u^h\cdot\nabla \u^h_{-h}+ \nabla \u^h\cdot \nabla^\bot \u^h_{-h}+\nabla^\bot \u^h\cdot \nabla^\bot \u^h_{-h} +\lambda\div(\u^h)\div(\u^h_{-h})\Big)+
\end{equation}
$$
2\xi\Big(\u^h_{-h}\cdot \nabla \u^h \cdot \nabla\xi+\u^h_{-h}\nabla^\bot\u^h\cdot\nabla \xi+\nabla \xi\cdot \nabla \u^h\cdot \u^h_{-h}+\nabla \xi \cdot \nabla^\bot \u^h\cdot \u^h_{-h}+\lambda \div(\u^h)\nabla \xi \cdot \u^h_{-h} \Big).
$$
Next we notice that 
$$
\u^h_{-h}(x)=\frac{\u(x)-\u(x+he_i)}{h}=-\u_h(x).
$$
Adding (\ref{firstindiffquot}) and (\ref{secondindiffquot}), rearranging 
the terms and dividing by $h$, we may conclude that
$$
\int_{B_1^+}\xi^2\Big(\big|\nabla \u_h+\nabla^\bot \u_h\big|^2+\lambda\div(\u_h)^2 \Big)\le
$$
$$
-2\int_{B_1^+}2\xi\Big(\u_h\cdot \nabla \u_h \cdot \nabla\xi+\u_h\nabla^\bot\u_h\cdot\nabla \xi+\nabla \xi\cdot \nabla \u_h\cdot \u_h+\nabla \xi \cdot \nabla^\bot \u_h\cdot \u_h+\lambda \div(\u_h)\nabla \xi \cdot \u_h \Big)
\le
$$
$$
\frac{1}{2}\int_{B_1^+}\xi^2\Big(|\nabla \u_h+\nabla^\bot \u_h|^2+\lambda\div(\u_h)^2 \Big)+ C\int_{B_1^+}|\u_h|^2|\nabla \xi|^2.
$$
In particular it follows, by letting $h\to 0$, that
$$
\int_{B_{1/2}^+}\Big(|\nabla \partial_{i}\u+ \nabla^\bot\partial_i \u_h|^2+\lambda\div(\partial_i \u)^2 \Big)\le
C\int_{B_1^+}\big| \partial_i \u\big|^2.
$$
By Kohn's inequality this implies that 
$\partial_i \u\in W^{1,2}(B_{1/2}^+;\R^3)$ for $i=1,2$. It directly
follows that $\u$ solves the Lame system in $B_{1/2}^+$ with boundary data given
by the restriction of a $W^{2,2}$ function to the boundary and $W^{2,2}$ 
regularity of $\u$ follows. The estimate given is a consequence of 
Lemma \ref{w12est} and the above. \qed

\section{Global Solutions, part 1. \\ Reduction of the System.}\label{GlobalSolPart1}

In this section we make a very useful reduction of solutions, with controlled 
growth, of the Lame system in the upper half space 
into a system with two unknown functions $\xi$ and $\tau$. Later we will
even be able to express the solution in terms of one harmonic function $\tau$.

We will call a solution $\u$ in the upper half space 
$\R^3_+=\R^3\cap \{x_3>0\}$ a global solution. In this section we will
only consider solutions in $\R^3$ for simplicity. In particular we will
utilise the \textsl{curl} operator which is much easier to express in $\R^3$
than in higher dimensions. There is however nothing in this section
that requires the dimension to be three and the reader may verify
that all the proofs works also in $\R^n$.  

\begin{lem}\label{ReduktionOfthe System}
Let $\u$ be a global solution to the Signorini problem and
$$
\liminf_{r\to \infty}\frac{\ln\big(\|\u\|_{\tilde{L}^2(B_r^+)}\big)}{\ln(r)}<2.
$$
Then there exist functions $\xi$ and $\tau$ such that
\begin{equation}\label{simpleruform}
\u=\nabla \xi +e_3 \tau
\end{equation}
and
\begin{equation}\label{xiequation}
\Delta \xi +\frac{\lambda+2}{\lambda+4}\frac{\partial \tau}{\partial x_3}=0 \quad \textrm{ in }\R^3_+
\end{equation}
\begin{equation}\label{tauequation}
\Delta \tau =0\quad \textrm{ in }\R^3_+
\end{equation}
\begin{equation}\label{bdryxitauonLambda}
\frac{\partial \xi}{\partial x_3}=\tau=0 \quad\textrm{ on } \Lambda_{\u}
\end{equation}
\begin{equation}\label{bdry1forxitauonOmega}
\tau=-2\frac{\partial \xi}{\partial x_3} \quad\textrm{ on } \Omega_{\u}
\end{equation}
\begin{equation}\label{bdry2forxitauonOmega}
\frac{\partial \tau}{\partial x_3}=-\frac{2\lambda+8}{3\lambda+8}\frac{\partial^2 \xi}{\partial x_3^2} \quad \textrm{ on }  \Omega_{\u}.
\end{equation}
\end{lem}
\textsl{Proof:} Let $\u$ be as in the Lemma and define
$$
v=\div(\u)\in W^{1,2}_{loc}(\R^3_+),
$$
$$
\w=\textrm{curl}(\u)=\big[u^3_2-u^2_3,u^1_3-u^3_1,u^2_1-u^1_2 \big]\in W^{1,2}_{loc}(\R^3_+).
$$
It is easy to see that $\Delta w^i=\Delta v=0$ for $i=1,2,3$. Moreover 
on $\Lambda_\u$ we have
$$
\frac{\partial w^3}{\partial x_3}= \frac{\partial}{\partial x_1}\Big(  \frac{\partial u^2}{\partial x_3}\Big) - \frac{\partial}{\partial x_2}\Big( \frac{\partial u^1}{\partial x_3}\Big)=0.
$$
And on $\Omega_\u$ 
$$
\frac{\partial w^3}{\partial x_3}= \frac{\partial^2 u^2 }{\partial x_3\partial x_1}-
\frac{\partial^2 u^1 }{\partial x_3\partial x_2}=
\frac{\partial}{\partial x_1}\Big(  -\frac{\partial u^3}{\partial x_2}\Big) - \frac{\partial}{\partial x_2}\Big( -\frac{\partial u^3}{\partial x_1}\Big)=0,
$$
where we have used (\ref{Bdry3}). In particular $w^3$ satisfies
$$
\begin{array}{ll}
\Delta w^3 =0 & \textrm{ in } \R^3_+ \\
\frac{\partial w^3}{\partial x_3}=0 & \textrm{ on } \Pi \\
\sup_{B_R^+}|w^3|\le C R^\alpha & \textrm{for }R\ge 1 \textrm{ and an }\alpha<1. 
\end{array}
$$
By the Liouville Theorem it follows that $w^3=$constant. But from Lemma
\ref{c1beta} it follows that for $R<1$ and some $\beta>0$  
$$
\sup_{B_R^+}|w^3|\le C R^{\beta}
$$
which implies that $w^3=0$.

We may therefore conclude that there exist a $\xi$ such that
\begin{equation}\label{poordog}
\Big( \frac{\partial \xi}{\partial x_1}, \frac{\partial \xi}{\partial x_2}\Big)=\big( u^1,u^2\big).
\end{equation}

Next we consider the equation for $w^1$
$$
0=\Delta w^1= \Delta \Big(\frac{\partial u^3}{\partial x_2}-\frac{\partial^2 \xi}{\partial x_2 \partial x_3} \big)=\frac{\partial }{\partial x_2}\Big( \Delta u^3- \Delta \frac{\partial \xi}{\partial x_3} \Big)
$$
$$
=
\frac{\partial }{\partial x_2}\Big( - \frac{2+\lambda}{2}\frac{\partial}{\partial x_3}\div(u)- \Delta \frac{\partial \xi}{\partial x_3} \Big).
$$
A similar consideration for $\Delta w^2$ implies that
$$
\frac{\partial }{\partial x_1}\Big( - \frac{2+\lambda}{2}\frac{\partial}{\partial x_3}\div(u)- \Delta \frac{\partial \xi}{\partial x_3} \Big)=0.
$$
That is
$$
\nabla' \Big(  \Delta \frac{\partial \xi}{\partial x_3}+\frac{2+\lambda}{2}\frac{\partial}{\partial x_3}\div(u)\Big)=0.
$$
Therefore 
$$
\Delta \frac{\partial \xi}{\partial x_3}+\frac{2+\lambda}{2}\frac{\partial}{\partial x_3}\div(u)=f(x_3)
$$
for some function $f(x_3)$. But $\xi$ is not determined up to functions in 
the $x_3$ variable, that is since the only condition on $\xi$ so far is that 
$\nabla' \xi =(u^1,u^2)$, so we may choose $\xi$ so that 
$$
\Delta \frac{\partial \xi}{\partial x_3}+\frac{2+\lambda}{2}\frac{\partial}{\partial x_3}\div(u)=0.
$$
In particular $u^3$ and $\frac{\partial \xi}{\partial x_3}$ satisfies the 
same partial differential equation, $\Delta\cdot =-\frac{\lambda+2}{2}\frac{\partial }{\partial x_3}\div(\u)$, and differ thus by a harmonic function which we will denote $\tau$.
The equations (\ref{simpleruform}) and (\ref{tauequation}) follows.

To show (\ref{xiequation}) we just notice that 
$$
0=\Delta \u + \frac{\lambda+2}{2}\nabla \div(\u)=\nabla \Big(\Delta \xi+\frac{\lambda+2}{2}\Big( \Delta \xi +\frac{\partial \tau}{\partial x_3}\Big) \Big).
$$
It immediately follows that
\begin{equation}\label{whatafuck}
\Delta \xi +\frac{\lambda+2}{\lambda+4}\frac{\partial \tau}{\partial x_3}=c_0,
\end{equation}
for some constant $c_0$.
By making the substitutions 
$$
\xi \to \xi+\frac{c_0}{2}x_3^2+\frac{a}{2}x_3^2
$$ 
and
$$
\tau \to \tau-\frac{\lambda+4}{\lambda+2}a x_3
$$
we may assume that the constant $c_0$ in (\ref{whatafuck})
is zero. Equation (\ref{xiequation}) follows. We want to point out that
the constant $a$ is arbitrary and that $\tau$ is therefore not determined 
up to linear functions $ax_3$, a fact that we will use later.

On $\Lambda_\u$ we have
\begin{equation}\label{rfv1}
0=u^3=\frac{\partial \xi}{\partial x_3}+\tau
\end{equation}
and
$$
0=\frac{\partial u^i}{\partial x_3}=\frac{\partial^2 \xi}{\partial x_i \partial x_3}
$$
for $i=1,2$. It follows that
\begin{equation}\label{rfv2}
\frac{\partial \xi}{\partial x_3}=c_i,
\end{equation}
where the constant $c_i$ may differ from component to component of $\Lambda_\u$.

The boundary conditions (\ref{bdryxitauonLambda}) follows from (\ref{rfv1}) and (\ref{rfv2}).

On $\Omega_\u$ we have
$$
0=\frac{\partial u^i}{\partial x_3}=\frac{\partial u^3}{\partial x_i}=
2\frac{\partial^2 \xi}{\partial x_i\partial x_3}+\frac{\partial \tau}{\partial x_i}.
$$
It follows that 
\begin{equation}\label{shitface}
\tau=-2\frac{\partial \xi}{\partial x_3}+\tilde{c}_i,
\end{equation}
where the constant $\tilde{c}_i$ may differ in different components of $\Omega_\u$.

By making the substitution 
\begin{equation}\label{oneofmanyxisubst}
\xi\to \xi+c_1 x_3
\end{equation}
we see that we may choose $\xi$
so that that the constant $c_i$ is zero for one component of $\Lambda_\u$.
In particular the boundary conditions (\ref{bdryxitauonLambda}), in that 
component $\Lambda_1$ of $\Lambda_\u$, follows from (\ref{rfv1}) and 
(\ref{rfv2}).

Next we notice that if $\Omega_i$ is any component of $\Omega_\u$ adjacent
to $\Lambda_1$ then by $C^{1,\beta}$ continuity of $\u$ it follows that
$\tilde{c}_i=0$. That is (\ref{bdry1forxitauonOmega}) holds in $\Omega_i$.
In particular if (\ref{bdryxitauonLambda}) is true in one component $\Lambda_i$
of $\Lambda_\u$ then (\ref{bdry1forxitauonOmega}) holds for all components
$\Omega_j$ adjacent to $\Lambda_i$. That is if $c_i=0$ in (\ref{rfv2}) for 
some component $\Lambda_i$ then $\tilde{c}_j=0$ in (\ref{shitface}) for each 
$j$ such that $\Omega_j$ is a component adjacent to $\Lambda_i$.

Conversely if (\ref{bdry1forxitauonOmega})  holds in a component 
$\Omega_i$ (that is $\tilde{c}_i=0$ in 
(\ref{shitface})) then by $C^{1,\beta}$ regularity (\ref{bdryxitauonLambda}) is 
true for each adjacent component $\Lambda_i$. So if we make the the substitution
(\ref{oneofmanyxisubst}) then it follows that $c_i=0$ and $\tilde{c}_i=0$
for all components of $\Lambda_\u$ and $\Omega_\u$. The boundary conditions
(\ref{bdry1forxitauonOmega}) and (\ref{bdryxitauonLambda}) holds on $\Pi$.

Finally we have on $\Omega$ that
$$
0=\frac{\partial u^3}{\partial x_3}+\frac{\lambda}{4}\div{u}=\frac{\partial^2 \xi}{\partial x_3^2}+
\frac{\partial \tau}{\partial x_3}+ \frac{\lambda}{4}\Big(\Delta \xi + \frac{\partial \tau}{\partial x_3} \Big)=\frac{\partial^2 \xi}{\partial x_3^2}+\frac{3\lambda+8}{2\lambda+8}\frac{\partial \tau}{\partial x_3},
$$
where we have used (\ref{xiequation}). This implies (\ref{bdry2forxitauonOmega}). \qed

Next we want to show that $\tau$ is an $x_3$-derivative of a harmonic function.

\begin{cor}\label{simplerformcor}
Let $\u$ be a global solution to the Signorini problem and
$$
\liminf_{r\to \infty}\frac{\|\u\|_{\tilde{L}^2(B_r^+)}}{\ln(r)}<2.
$$
Then there exist functions $\xi$ and $\tau$ such that
\begin{equation}\label{simpleruformcor}
\u=\nabla \xi +e_3 \frac{\partial \tau}{\partial x_3}
\end{equation}
and
\begin{equation}\label{xiequationcor}
\Delta \xi +\frac{\lambda+2}{\lambda+4}\frac{\partial^2 \tau}{\partial x_3^2}=0 \quad \textrm{ in }\R^3_+
\end{equation}
\begin{equation}\label{tauequationcor}
\Delta \tau =0\quad \textrm{ in }\R^3_+
\end{equation}
\begin{equation}\label{bdryxitauonLambdacor}
\frac{\partial \xi}{\partial x_3}=\frac{\partial \tau}{\partial x_3}=0 \quad\textrm{ on } \Lambda_{\u}
\end{equation}
\begin{equation}\label{bdry1forxitauonOmegacor}
\frac{\partial \tau}{\partial x_3}=-2\frac{\partial \xi}{\partial x_3} \quad\textrm{ on } \Omega_{\u}
\end{equation}
\begin{equation}\label{bdry2forxitauonOmegacor}
\frac{\partial^2 \tau}{\partial x_3^2}=-\frac{2\lambda+8}{3\lambda+8}\frac{\partial^2 \xi}{\partial x_3^2} \quad \textrm{ on }  \Omega_{\u}.
\end{equation}
\end{cor}
\textsl{Proof:} Let $\xi$ and $\tau$ be as in Lemma \ref{ReduktionOfthe System}
and let $\chi$ be the solution of 
$$
\begin{array}{ll}
\Delta \chi=-\frac{\lambda+2}{\lambda+4}\frac{\partial \tau}{\partial x_3} &
\textrm{ in } \R^3_+ \\
\frac{\partial \chi}{\partial x_3}=0 & \textrm{ on } \Pi.
\end{array}
$$
Then
$$
\begin{array}{ll}
\Delta 2(\chi-\xi)=0 &
\textrm{ in } \R^3_+ \\
\frac{\partial 2(\chi-\xi)}{\partial x_3}=2\frac{\partial \xi}{\partial x_3}=-\tau & \textrm{ on } \Pi.
\end{array}
$$
In particular if we denote 
$$
\tilde{\tau}=\tau +2\frac{\partial (\chi-\xi)}{\partial x_3}
$$
then
$$
\begin{array}{ll}
\Delta \tilde{\tau}=0 & \textrm{ in } \R^3_+ \\
\tilde{\tau}=0 & \textrm{ on } \Pi.
\end{array}
$$
Moreover $\sup_{B_R^+}|\tilde{\tau}|\le CR^{1+\alpha}$ so by Liouville's Theorem
$\tilde{\tau}=a x_3$ for some constant $a$. That is
\begin{equation}\label{almosttauisaderivative}
\tau=ax_3-2\frac{\partial (\chi-\xi)}{\partial x_3}.
\end{equation}
But as we pointed out in the discussion right after (\ref{whatafuck})
$\tau$ is only determined up to linear functions $ax_3$. We may therefore
choose a $\tau$ so that $a=0$ in (\ref{almosttauisaderivative}).

We have thus shown that $\tau$ in Lemma \ref{ReduktionOfthe System}
is expressible as the $x_3$-derivative of a harmonic function. The
corollary follows. \qed

\begin{cor}
For each $R>0$ we have, with $\tau$ as in Corollary \ref{simplerformcor}, 
$$
\frac{\partial \tau}{\partial x_3}\in W^{2,2}(B^+_R).
$$
\end{cor}
\textsl{Proof:} By Corollary \ref{simplerformcor}
$$
\frac{\partial \xi}{\partial x_3}+\frac{\partial \tau}{\partial x_3}=u^3\in 
W^{2,2}(B_R^+).
$$
It follows that
$$
\frac{\partial^2 \tau}{\partial x_i\partial x_3}=\frac{\partial u^3}{\partial x_i}-\frac{\partial^2 \xi}{\partial x_i\partial x_3}=\frac{\partial u^3}{\partial x_i}-\frac{\partial u^i}{\partial x_3}\in W^{1,2}(B_R^+),
$$
for $i=1,2$, where we have used that
$$
\frac{\partial \xi}{\partial x_i}=u^i.
$$
By the trace theorem we therefore have
$$
\frac{\partial^2\tau}{\partial x_i\partial x_3}\in W^{1/2,2}(\Pi\cap B_R).
$$
That in turn results in
$$
\frac{\partial \tau}{\partial x_3}\in W^{3/2,2}(\Pi\cap B_R).
$$
But $\frac{\partial \tau}{\partial x_3}$ is harmonic so we may conclude that
$$
\frac{\partial \tau}{\partial x_3}\in W^{2,2}(B_R^+)
$$
\qed

\textbf{Remark:} It is not hard to show that $\tau\in W^{3,2}(B_R)$, but we 
have no need for that stronger statement in what follows.

\begin{cor}\label{taubisiszeroonOmega}
Let $\xi$ and $\tau$ be as in Corollary \ref{simplerformcor}. Then
$$
\frac{\partial^2 \tau}{\partial x_3^2}=0
$$
on $\Omega_\u$.
\end{cor}
\textsl{Proof:} Let 
$$
v=\xi-\frac{\lambda+2}{2(\lambda+4)}\frac{\partial \tau}{\partial x_3}x_3.
$$
Then
$$
\Delta v=\Delta \xi -\frac{\lambda+2}{\lambda+4}\frac{\partial^2\tau}{\partial x_3^2}=0
$$
in $\R^3_+$. Moreover, on $\Pi$,
$$
\frac{\partial v}{\partial x_3}=-\frac{\lambda+3}{\lambda+4}\frac{\partial \tau}{\partial x_3}.
$$
This implies that
$$
v=-\frac{\lambda+3}{\lambda+4}\tau
$$
or equivalently that
\begin{equation}\label{aonetimer}
\xi=\frac{\lambda+2}{2(\lambda+4)}\frac{\partial \tau}{\partial x_3}x_3-\frac{\lambda+3}{\lambda+4}\tau.
\end{equation}
Using (\ref{bdry2forxitauonOmegacor}) we conclude that on $\Omega$
$$
-\frac{3\lambda+8}{2\lambda+8}\frac{\partial^2 \tau}{\partial x_3^2}=
\frac{\partial^2 \xi}{\partial x_3^2}=\frac{\lambda+2}{\lambda+4}\frac{\partial^2\tau}{\partial x_3^2}-\frac{\lambda+3}{\lambda+4}\frac{\partial^2\tau}{\partial x_3^2}=-\frac{1}{\lambda+4}\frac{\partial^2\tau}{\partial x_3^2}.
$$
Rearranging terms we may conclude that
$$
\frac{3\lambda+6}{2\lambda+8}\frac{\partial^2\tau}{\partial x_3^2}=0
$$
on $\Omega$ and the Corollary follows.\qed

We will need to refer to (\ref{aonetimer}) on several occasions
later so it is convenient to formulate that equality as a Corollary.

\begin{cor}\label{NEWCOR4}
Let $\xi$ and $\tau$ be as in Corollary \ref{simplerformcor}. Then
$$
\xi=\frac{\lambda+2}{2(\lambda+4)}\frac{\partial \tau}{\partial x_3}x_3-\frac{\lambda+3}{\lambda+4}\tau.
$$
\end{cor}

\begin{lem}\label{forreflectionofv}
Let $\u$ be a global solution to the Signorini problem as in Corollary
\ref{simplerformcor} and denote $v=\div(\u)$. Then
$$
v=\frac{2}{\lambda+4}\frac{\partial^2 \tau}{\partial x_3^2}
$$
where $\tau$ is as in Corollary \ref{simplerformcor}.

In particular, by Corollary \ref{taubisiszeroonOmega},
$$
v=0 \quad \textrm{ on }\Omega.
$$
\end{lem}
\textsl{Proof:} With the notation of Corollary \ref{simplerformcor}
we have
$$
v=\div(\u)=\Delta \xi+\frac{\partial^2 \tau}{\partial x_3^2}=\frac{2}{\lambda+4}\frac{\partial^2 \tau}{\partial x_3^2},
$$
where we have used (\ref{xiequationcor}). \qed

\section{A result by Benedicks.}

In this section we will remind ourselves of a Theorem due to M. Benedicks 
\cite{Benedicks}. We will formulate the theorem slightly differently from 
Benedick's for convenience. We will however not be able to directly apply
the theorem. Therefore we give a slightly different version of it, whose 
proof is a simple consequence of \cite{FriedlandHayman}. We will also
prove that global solutions of the Signorini problem with a bound at 
infinity are uniquely determined by the set $\Lambda_\u$. Later we will
refine this result somewhat and it will not be used in the rest of the paper.

\begin{prop}\label{BenThm}
Let $\Lambda\subset \Pi$ be a given set in $\Pi$ and  
$$
\mathcal{P}_\Delta(\Lambda)=\{u\in W^{1,2}(\R^n\setminus \Lambda);\; \Delta u=0
\textrm{ in }\R^n\setminus \Lambda,\;
u\ge 0\textrm{ in }\R^n\setminus \Lambda, u=0\textrm{ on }\Lambda\}.
$$
Then the set $\mathcal{P}_\Delta(\Lambda)$ is a one or two dimensional set.

Moreover if $u$ is even in $x_n$ then $\mathcal{P}_\Delta(\Lambda)$
is one dimensional.
\end{prop}
For a proof see \cite{Benedicks}. 

\begin{lem}\label{hassign}
Let $\Lambda\subset \Pi$ be a given set in $\Pi$ and  
$$
\begin{array}{ll}
\Delta u=0 & \textrm{ in }\R^3_+ \\
u=0 & \textrm{ on } \Lambda \\
\frac{\partial u}{\partial x_3}=0 & \textrm{ on } \Pi\setminus \Lambda \\
\sup_{B_R^+}|u|\le CR^\alpha & \textrm{ for } R>0 \textrm{ and an }\alpha<1
\end{array}
$$
then $u$ has a sign. That is $u\ge 0$ or $u\le 0$. In particular Proposition 
\ref{BenThm} applies.
\end{lem}
\textsl{Proof:} We may extend $u$ by an even reflection to
$$
\bar{u}(x)=\left\{
\begin{array}{ll}
u(x) & \textrm{ if } x_3\ge 0 \\
u(x_1,x_2,-x_3) & \textrm{ if } x_3<0.
\end{array}
\right.
$$
Then $\bar{u}$ is harmonic in $\R^3\setminus \Lambda$. If $u$ does not have a 
sign then $\bar{u}^\pm \ne 0$. But $\sup_{B_R}|\bar{u}^\pm|\le CR^\alpha$.
However, since the supports of $\bar{u}^\pm$ are disjoint it follows by 
Friedland Hayman's Theorem \cite{FriedlandHayman} that at least one of 
$\bar{u}^\pm$ must have
at least linear growth at infinity. This is a contradiction. Therefore 
we can conclude that either $\bar{u}^+=0$ or $\bar{u}^-=0$. \qed

Next we prove a version of Benedicks' Theorem for the Signorini problem.
We will however need a refined version later and we will therefore not 
need the following Lemma in what follows.

\begin{lem}\label{Bennedicks}
Let $\mathcal{P}$ be the set of $W^{2,2}$ solutions to
\begin{equation}\label{1i}
\begin{array}{ll}
\Delta \u+\frac{2+\lambda}{2} \nabla \textrm{div}(\u)=0 & \textrm{ in } \R^{3}_+
\end{array}
\end{equation}
\begin{equation}\label{1ii}
\begin{array}{ll}
u^3=0 & \textrm{ on }\Lambda
\end{array}
\end{equation}
\begin{equation}\label{1iii}
\begin{array}{ll}
\frac{\partial u^3}{\partial x_n}+\lambda \textrm{div}(U)=0 & \textrm{ on }
\Pi\setminus\Lambda
\end{array}
\end{equation}
\begin{equation}\label{1iv}
\begin{array}{ll}
\frac{\partial u^i}{\partial x_3}+\frac{\partial u^3}{\partial x_i}=0 & \textrm{ on } \Pi \textrm{ for }i=1,2
\end{array}
\end{equation}
\begin{equation}\label{1v}
\begin{array}{ll}
\sup_{B_R^+}|\u|\le C(1+R^{1+\alpha}) & \textrm{for some }\alpha<1 \textrm{ and all }R\ge 1
\end{array}
\end{equation}
\begin{equation}\label{1vi}
\begin{array}{ll}
\sup_{B_r^+}|\u|\le C r^{1+\beta} & \textrm{for some }\beta>0 \textrm{ and all }r<1.
\end{array}
\end{equation}
Here $\Lambda$ is a fixed set in $\Pi$.

Then $\mathcal{P}$ is a one dimensional set. That is 
$\mathcal{P}=\{t\v;\; t\in \R\}$ for some $\v\in W^{2,2}$.
\end{lem}
\textsl{Proof:} Since $\u\in \mathcal{P}$ implies that $\u\in W^{2,2}$
we know that $v\equiv \textrm{div}(\u)\in W^{1,2}$. Taking the 
divergence of (\ref{1i}) it follows that $\Delta v=0$ in 
$\R^3\setminus \Lambda$. Using (\ref{1iv}) we see that 
$$
\frac{\partial^2 u^i}{\partial x_j\partial x_n}+\frac{\partial^2 u^3}{\partial x_j\partial x_i}=0
$$ 
on $\Pi$ for $i,j=1,2$ and thus 
\begin{equation}\label{normalderivativesofthediv}
\frac{\partial v}{\partial x_3}=-\frac{\partial^2 u^3}{\partial x_1^2}-\frac{\partial^2 u^3}{\partial x_1^2}+\frac{\partial^2 u^3}{\partial x_3^2}
\end{equation}
on $\Pi$. From (\ref{1ii}) we see that $\Delta' u^3=0$ on $\Lambda$
$$
\frac{\partial v}{\partial x_3}=\frac{\partial^2 u^3}{\partial x_3^2}
\textrm{ on }\Lambda
$$
Therefore (\ref{1i}) implies that 
$$
0=\Delta u^3+\frac{2+\lambda}{2}\frac{\partial v}{\partial x_3}=\frac{4+\lambda}{2}\frac{\partial^2 u^3}{\partial x_3^2}
$$ 
on $\Lambda$, where we also used $\Delta' u^3=0$ on $\Lambda$. Thus, using (\ref{normalderivativesofthediv})
\begin{equation}\label{normalderivativeofvonLambda}
\frac{\partial v}{\partial x_3}=0 \quad \textrm{ on }\Lambda.
\end{equation}
Thus $v$ satisfies the following boundary value problem
\begin{equation}\label{boundaryvalueprobfordiv}
\begin{array}{ll}
\Delta v=0 & \textrm{ in } \R^3_+ \\
\frac{\partial v}{\partial x_3}=0 & \textrm{ on } \Lambda.
\end{array}
\end{equation}
Using (\ref{1v}) and (\ref{1vi}) we deduce that
\begin{equation}\label{growthofthediv1}
\sup_{B_R^+}|v|\le C\big( 1+R^\alpha\big)
\end{equation}
where $\alpha<1$ and
\begin{equation}\label{growthofthediv2}
\sup_{B_r^+}|v|\le C r^\beta
\end{equation}
where $\beta>0$ and for every $r<1$.

We make the following claim.

\textbf{Claim 1:} {\sl $v=0$ on $\Omega$.}

\textsl{Proof of claim 1:} This was proved in Lemma \ref{forreflectionofv}.

$ $

For the next claim we notice that we may define a antisymmetric solution
in $\R^3\setminus \Lambda$ by reflecting $v$
$$
v(x)=\left\{
\begin{array}{ll}
v(x) & \textrm{ if } x_3\ge 0 \\
v(x_1,x_2,-x_3) & \textrm{ if } x_3< 0.
\end{array}
\right.
$$
Then $v$ is harmonic in $\R^3\setminus \Lambda$ and satisfies the 
Neumann condition in (\ref{boundaryvalueprobfordiv}) on both sides of 
$\Lambda$.

$ $

\textbf{Claim 2.}{\sl The set of antisymmetric solutions to 
(\ref{boundaryvalueprobfordiv}), (\ref{growthofthediv1}) and 
(\ref{growthofthediv2}) is one dimensional.}

\textsl{Proof of claim 2:} By antisymmetry it is enough to show that
solutions to
$$
\begin{array}{ll}
\Delta v=0 & \textrm{ in } \R^3_+ \\
v=0 & \textrm{ on } \Pi\setminus \Lambda \\
\frac{\partial v}{\partial x_3}=0 & \textrm{ on } \Lambda
\end{array}
$$
are one dimensional. This is a special case of a Proposition \ref{BenThm}. 
The claim follows.

$ $

We are now ready to prove the Lemma. It is clearly enough to show
that if $\u$ and $\v$ are two  solutions normalized
so that $\textrm{div}(\u)=\textrm{div}(\v)$ then $\u=\v$ since 
by claim 1 and 2 such a normalization always exists for non-vanishing
solutions.

Let $\u,\v\in \mathcal{P}$ and $\textrm{div}(\u)=\textrm{div}(\v)$. Then
$\w=\u-\v\in \mathcal{P}$ and $\textrm{div}(\w)=0$. From (\ref{1i})
we conclude that 
\begin{equation}\label{wisharmonic}
\Delta \w=0. 
\end{equation}
From (\ref{1iv}) it follows, as in the 
beginning of this proof, that
$\frac{\partial w^i}{\partial x_3}=0$ on $\Pi$ for $i=1,2$. Using
(\ref{1v}) and the Liouville Theorem  we may conclude that
$$
\begin{array}{l}
w^1(x)=a_1 x_1+b_1x_2 \\
w^2(x)=a_2 x_1+b_2x_2
\end{array}
$$
for some constants $a_1,a_2,b_1$ and $b_2$. Using (\ref{1vi}) we
see that $a_1=a_2=b_1=b_2=0$. From (\ref{1iv}) we may deduce that
$$
\frac{\partial w^3}{\partial x_i}=0
$$
for $i=1,2$ since $\frac{\partial w^i}{\partial x_i}=0$ for $i=1,2$.
Therefore $w^3=\textrm{constant}=a$ on $\Pi$. Also by (\ref{wisharmonic})
we have that $\Delta w^3=0$ so by the Liouville Theorem again we may conclude 
from (\ref{1v}) that $w^3(x)=a+bx_3$ for some $a,b\in \R$. But
$$
0=\div(\w)=b
$$
on $\Omega$ which implies that $w^3$ is constant. Using (\ref{1vi}) we may 
deduce that $w^3=0$. The Lemma follows. \qed

\begin{lem}\label{Onepoint}
Let $\u$ be a solution to the Signorini problem in $\mathbb{R}^2_+$ and
$$
\limsup_{r\to \infty}\frac{\ln\big( \|\u\|_{\tilde{L}^2(B_r)}\big)}{\ln(r)}<2
$$
then $\Gamma$ contains at most one point.
\end{lem}
\textsl{Proof:} Taking the divergence of $L\u=0$ we see that
$\textrm{div}(\u)$ is harmonic. From claim 1 in the proof of Lemma 
\ref{Bennedicks} we also know that $\textrm{div}(\u)$ is antisymmetric
so 
\begin{equation}\label{divzeroonlambda}
\begin{array}{ll}
\textrm{div}(\u)=0 & \textrm{ on } \Sigma.
\end{array}
\end{equation}
Also as in Lemma \ref{Bennedicks} $\textrm{div}(\u)$ has a sign.

Similarly let 
$$
\textrm{curl}(\u)=\frac{\partial u^1}{\partial x_2}-
\frac{\partial u^2}{\partial x_1}.
$$
Then $\Delta \textrm{curl}(\u)=0$ and 
\begin{equation}\label{curlzerolambda}
\begin{array}{ll}
\textrm{curl}(\u)=\frac{\partial u^1}{\partial x_2}-\frac{\partial u^2}{\partial x_1}=-2\frac{\partial u^2}{\partial x_1}=0 & \textrm{ on } \Lambda,
\end{array}
\end{equation}
where we have used that 
$$
\frac{\partial u^1}{\partial x_2}-\frac{\partial u^2}{\partial x_1}=0 \textrm{ and } \frac{\partial u^2}{\partial x_1}=0
$$
on $\Lambda$.
Moreover 
$$
\textrm{curl}\Big( (1+(2+\lambda)/2)\textrm{div}(\u), \textrm{curl}(\u)\Big)=
\Delta u^2+\frac{2+\lambda}{2}\frac{\partial \textrm{div}(\u)}{\partial x_2}=0.
$$
That implies that there exist a $v$ such that
$$
\nabla v= \big( (1+(2+\lambda)/2)\textrm{div}(\u), \textrm{curl}(\u)\big).
$$
It is easy to see that
$$
\Delta v= \Delta u^1+\frac{2+\lambda}{2}\frac{\partial \textrm{div}(\u)}{\partial x_1}=0.
$$
Moreover
\begin{equation}\label{boundaryforpotential}
\begin{array}{ll}
\frac{\partial v}{\partial x_2}=\textrm{curl}(\u)=0 & \textrm{ on } \Lambda \\
\frac{\partial v}{\partial x_1}=\textrm{div}(\u)=0 & \textrm{ on } \Sigma.
\end{array}
\end{equation}
Since both $\textrm{div}(\u)$ and $\textrm{curl}(\u)$ have a sign, by
Lemma \ref{hassign}, so does
$\frac{\partial v}{\partial x_1}$ and $\frac{\partial v}{\partial x_2}$.
We will assume that $\frac{\partial v}{\partial x_1}\ge 0$ and, 
the other case is handled similarly.

Next we assume that $a,b\in \Gamma$ and $a\ne b$. If two such points exist 
then we may chose them so that $a$ is the left boundary of an interval 
$(a,a_0)\subset \Sigma$ and $b$ is the right endpoint in an interval
$(b_0,b)\subset \Sigma$. Observe that from (\ref{boundaryforpotential})
it follows that $v$ is constant in each component of $\Sigma$ and thus
a solution to the thin obstacle problem in $B_r^+(a)$ and $B_r^+(b)$ if
$r$ is small enough.

The asymptotic expansion at free boundary points for the thin obstacle problem 
is known \cite{AC} and we may conclude that in polar coordinates
$$
v(x_1-a,x_2)=v(a,0)-\alpha r^{k+1/2}\sin\big( (k+1/2)\phi \big)+O(r^{k+1+1/2})
$$ 
and
$$
v(x_1-b,x_2)=v(b,0)-\beta r^{k+1/2}\cos\big( (l+1/2)\phi \big)+O(r^{l+1+1/2})
$$ 
for some $\alpha,\beta\in \mathbb{R}_+$ and integers $k,l\ge 1$. This implies 
that
$$
\frac{\partial v}{\partial x_2} \le 0
$$
close to $a$ and
$$
\frac{\partial v}{\partial x_2} \ge 0
$$
close to $b$, with strict inequality away from $\{x_2=0\}$. 
Since $\frac{\partial v}{\partial x_2}$ has a sign we
get a contradiction. \qed

\section{Global Solutions, part 2.}

In this section we explicitly calculate the global homogeneous solutions
In $\R^2_+$ to (\ref{Lame})-(\ref{Bdry3}). 

\begin{lem}\label{eigenfunctionsinR2}
Let $\u(x,y)=\big(u(x,y),v(x,y)\big)$ be a global homogeneous solution of order $1+\alpha>0$ to the Lame system with $\lambda\ne-2$ in $\mathbb{R}^2_+$:
$$
\begin{array}{ll}
\Delta u+\frac{2+\lambda}{2}\frac{\partial \div (\u)}{\partial x}=0 & \textrm{ in } \mathbb{R}^2_+ \\
\Delta v+\frac{2+\lambda}{2}\frac{\partial \div (\u)}{\partial y}=0 & \textrm{ in } \mathbb{R}^2_+ \\
v(x,0)=0 & \textrm{ on } \{x>0\} \\
\frac{\partial u}{\partial y}(x,0)=0 &  \textrm{ on } \{x>0\} \\
\frac{\partial u}{\partial y}(x,0)+\frac{\partial v}{\partial x}(x,0)=0 & \textrm{ on } \{x<0\}\\
\frac{\partial v}{\partial y}(x,0)+\frac{\lambda}{4}\div(\u)(x,0)=0 & \textrm{ on }\{x<0\} \\
u\in W^{2,2}_{loc}(\R^2_+) &
\end{array}
$$
then $\alpha\in \mathbb{N}$ or $\alpha\in \mathbb{N}+\frac{1}{2}=\{1/2,3/2,5/2,...\}$, and 
in polar coordinates, $x=r\cos(\phi)$ and $y=r\sin(\phi)$, we have for some 
$a\in \R$
\begin{enumerate}
\item 
$$
\begin{array}{l}
u(r,\phi)=r^{1+\alpha}\Big(\frac{10+3\lambda+\alpha(2+\lambda)}{8(1+\alpha)}\cos\big((1+\alpha)\phi \big)-\frac{2+\lambda}{8}\cos\big((1-\alpha)\phi \big)\Big)a \\
v(r,\phi)=r^{1+\alpha}\Big( \frac{2-\lambda-\alpha(2+\lambda)}{8(1+\alpha)}\sin\big((1+\alpha)\phi  \Big)
-\frac{2+\lambda}{8}\sin\big((1-\alpha)\phi \big)  \Big)a.
\end{array}
$$
if $\alpha\in \mathbb{N}$ and 
\item
$$
\begin{array}{l}
u(r,\phi)=r^{1+\alpha}\Big(\frac{6+\lambda+\alpha(2+\lambda)}{8(1+\alpha)}\cos\big((1+\alpha)\phi \big)-\frac{2+\lambda}{8}\cos\big((1-\alpha)\phi \big)\Big)a\\
v(r,\phi)=r^{1+\alpha}\Big( \frac{6+\lambda-\alpha(2+\lambda)}{8(1+\alpha)}\sin\big((1+\alpha)\phi \big) 
-\frac{2+\lambda}{8}\sin\big((1-\alpha)\phi \big)  \Big)a.
\end{array}
$$
if $\alpha\in \mathbb{N}+\frac{1}{2}$.
\end{enumerate}
\end{lem}
\textsl{Proof:} Denoting $w=\div(\u)\in W^{1,2}$ as we have done before we see
that
$$
\begin{array}{ll}
\Delta w =0 & \textrm{ in } \R^2_+ \\
\frac{\partial w}{\partial y}(x,0)=0 & \textrm{ on }\{x>0\}.
\end{array}
$$
Also $w$ will be homogeneous of order $\alpha$. That is, in polar coordinates,
\begin{equation}\label{wForm}
w(r,\phi)=r^\alpha \big(a \cos(\alpha \phi)+b\sin(\alpha \phi) \big).
\end{equation}
Using that 
$$
\frac{\partial w}{\partial \phi}(r,0)=0
$$
we can deduce that $b=0$.

Next we consider the ordinary differential equation
$$
\div(\u(r,\phi))=w(r,\phi).
$$
Using that $\u$ is homogeneous of order $1+\alpha$ it is easy to see that 
$u$ is of the form 
$$
u(r,\phi)=
$$
$$
= r^{1+\alpha}\Big(a_u \cos\big((1+\alpha)\phi \big)+b_u \sin\big((1+\alpha)\phi \big)+c_u \cos\big((1-\alpha)\phi \big)+d_u \sin\big((1-\alpha)\phi \big) \Big) 
$$
and that
$$
v(r,\phi)=
$$
$$
= r^{1+\alpha}\Big(a_v \cos\big((1+\alpha)\phi \big)+b_v \sin\big((1+\alpha)\phi \big)+c_v \cos\big((1-\alpha)\phi \big)+d_v \sin\big((1-\alpha)\phi \big) \Big).
$$
Also $(u,v)$ solves
\begin{equation}\label{uEQ}
\Delta u+\frac{2+\lambda}{2}\frac{\partial w}{\partial x}=0 \textrm{ in }R^2_+
\end{equation}
\begin{equation}\label{vEQ}
\Delta v+\frac{2+\lambda}{2}\frac{\partial w}{\partial y}=0 \textrm{ in }R^2_+
\end{equation}
\begin{equation}\label{vBDRYphi0}
v(r,0)=0 
\end{equation}
\begin{equation}\label{uBDRYphi0}
\frac{\partial u}{\partial \phi}(r,0)=0
\end{equation}
\begin{equation}\label{BDRY1phiPi}
\frac{1}{r}\frac{\partial v}{\partial \phi}(r,\pi)-\frac{\lambda}{4}w(r,\pi)=0
\end{equation}
\begin{equation}\label{BDRY2phiPi}
\frac{\partial v}{\partial r}(r,\pi)+\frac{1}{r}\frac{\partial u}{\partial \phi}(r,\pi)=0
\end{equation}
and
\begin{equation}\label{DIVEQ}
\div\big((u,v)\big)=w \textrm{ in } \R^2_+.
\end{equation}
From (\ref{vBDRYphi0}) we may deduce that 
$$
c_v=-a_v
$$ 
and from (\ref{uBDRYphi0}) it follows that
$$
b_u=-\frac{1-\alpha}{1+\alpha}d_u.
$$
Using (\ref{uEQ}), (\ref{wForm}) and $b=0$ we may deduce that
$$
d_u=0
$$
and
$$
c_u=-\frac{2+\lambda}{8}a.
$$
Similarly from (\ref{vEQ}) we deduce that
$$
a_v=0
$$
and
$$
d_v=-\frac{(2+\lambda)}{8}a.
$$
Equation (\ref{DIVEQ}) implies that
$$
b_v=\frac{6+\lambda}{4(1+\alpha)}a-a_u.
$$

Next we use equation (\ref{BDRY2phiPi}) which implies that
either
$$
e^{2\pi \alpha}=1
$$
that is $\alpha\in \mathbb{N}$ or
\begin{equation}\label{EQforau1}
a_u=\frac{6+\lambda+\alpha(2+\lambda)}{8(1+\alpha)}a.
\end{equation}

From equation (\ref{BDRY1phiPi}) we deduce that either
$$
e^{2\pi \alpha}=-1
$$
that is $\alpha\in \mathbb{N}+\frac{1}{2}$ or
\begin{equation}\label{EQforau2}
a_u=\frac{10+3\lambda+\alpha(2+\lambda)}{8(1+\alpha)}a.
\end{equation}

Both (\ref{EQforau1}) and (\ref{EQforau2}) holds only if $\lambda=-2$.
Therefore we must have either $\alpha\in \mathbb{N}$ and (\ref{EQforau2})
holds or $\alpha\in \mathbb{N}+\frac{1}{2}$ and (\ref{EQforau1}) holds.

In case $\alpha\in \mathbb{N}+\frac{1}{2}$ then (\ref{EQforau1}) implies that
$$
u(r,\phi)=r^{1+\alpha}\Big(\frac{6+\lambda+\alpha(2+\lambda)}{8(1+\alpha)}\cos\big((1+\alpha)\phi \big)-\frac{2+\lambda}{8}\cos\big((1-\alpha)\phi \big)\Big)a
$$
and
$$
v(r,\phi)=r^{1+\alpha}\Big( \frac{6+\lambda-\alpha(2+\lambda)}{8(1+\alpha)}\sin\big((1+\alpha)\phi \big) 
-\frac{2+\lambda}{8}\sin\big((1-\alpha)\phi \big)  \Big)a.
$$

And if $\alpha\in \mathbb{N}$ then
$$
u(r,\phi)=r^{1+\alpha}\Big(\frac{10+3\lambda+\alpha(2+\lambda)}{8(1+\alpha)}\cos\big((1+\alpha)\phi \big)-\frac{2+\lambda}{8}\cos\big((1-\alpha)\phi \big)\Big)a
$$
and
$$
v(r,\phi)=r^{1+\alpha}\Big( \frac{2-\lambda-\alpha(2+\lambda)}{8(1+\alpha)}\sin\big((1+\alpha)\phi \big) 
-\frac{2+\lambda}{8}\sin\big((1-\alpha)\phi \big)  \Big)a.
$$
\qed

We will also define a normalized $L^2$ norm that scales like the $L^{\infty}$
norm which will be convenient later.

\begin{defi}\label{definitiontildeL2}
We will use the notation $\tilde{L}^p(\Omega)$ for the average $L^p$ space with norm
$$
\|\u\|_{\tilde{L}^p(\Omega)}=\bigg( \frac{1}{|\Omega|}\int_{\Omega}|\u|^p\bigg)^{1/p}
$$
where $|\Omega|$ is the measure of $\Omega$.
\end{defi}

Inspired by Lemma \ref{eigenfunctionsinR2} we make the following definition 
of the global normalised homogeneous two dimensional solutions.

\begin{defi}\label{2dsolutiondef}
We will denote by $\p_{1/2}(x_1,x_n)$, $\p_{1}(x_1,x_n)$, $\p_{3/2}(x_1,x_n)$, 
etc. the homogeneous solution to the Signorini problem in $\mathbb{R}^2_+$
with $\|\p_{1+\alpha}\|_{\tilde{L}^2(B_1^+)}=1$ of order $1/2$, $1$, $3/2$, etc. 
as specified in Lemma \ref{eigenfunctionsinR2}.

We will also use the notation 
$\p_{k+1/2}^\xi(x)=|\xi|\p_{k+1/2}\big(\frac{\xi}{|\xi|}\cdot x',x_3\big)$ where $\xi\in \Pi$.
And 
$$
\tilde{\p}_{k+1/2}(x_1,x_n)=\frac{\partial \p_{k+3/2}(x_1,x_n)}{\partial x_n}.
$$
\end{defi}

\begin{lem}\label{eigenfinctforhalfspace}
Let $\w=(w^1,w^2,w^3)\in W^{1,2}(B_1^+)$ be a solution to the following linear problem
$$
\begin{array}{ll}
\Delta \w+\frac{2+\lambda}{2}\nabla \div(\w)=0 & \textrm{ in } B_1^+ \\
w^3(x_1,x_2,0)=0 & \textrm{ on } \{x_1>0\} \cap \Pi \\
\frac{\partial w^i}{\partial x_3}+\frac{\partial w^3}{\partial x_i}=0 & \textrm{ on }\Pi \textrm{ for }i=1,2 \\
\frac{\partial w^3}{\partial x_3}+\frac{\lambda}{4}\div(\w)=0 & \textrm{ on }
\{x_1<0\}\cap\Pi \\
\|\w\|_{L^\infty(B_1^+)}\le 1
\end{array}
$$
then $\w=\sum_{i=1}^\infty a_i \q_i$ where $\q_i$ is a homogeneous solution of order $i/2$ to the same problem.

Furthermore if $\w\in W^{2,2}(B_1^+)$ then
$$
\w=a\p_{3/2}+\sum_{2=0}^\infty a_i \q_i.
$$
\end{lem}
For a brief sketch of a proof see Appendix 2.

\section{Flatness of the Free Boundary}\label{sectionFlatness}

In this section we introduce the first fundamental idea in the paper and
show that at non-regular points the free boundary is flat.

From here on we will no longer, with the exception of Lemma 
\ref{fundamentalthinobstlemma}, need any explicit 
calculations using the curl operator and we will therefore write $\R^n$ instead 
of $\R^3$. Some of the ideas in this section to control the growth of blow up 
sequences comes from \cite{ASU}.

\begin{prop}\label{flatness}
Let $\u$ be a solution to the Signorini problem in $B_1^+$ and assume that
$0\in \Gamma$. Assume furthermore that
\begin{equation}\label{somekindofnonregularitycondition}
\liminf_{r\to 0}\frac{\ln\big(\|\u\|_{\tilde{L}^2(B_r^+)} \big)}{\ln (r)}<2
\end{equation}
then there exists a sub-sequence $r_j\to 0$ such that
$$
\frac{u(r_j x)}{\|\u\|_{\tilde{L}^2(B_{r_j}^+)}}\to \v,
$$
where $v$ is a global solution to the Signorini problem and furthermore, after
a rotation, $\Lambda_v=\{x\in \mathbb{R}^{n-1};\; x_1 \ge 0\}$.
\end{prop}
\textsl{Proof:} Assume that the limit in 
(\ref{somekindofnonregularitycondition}) is less than two and call the limit  
$\gamma<2$ and let $\alpha=\gamma/2$, then $0<\alpha<1$ and $1+\alpha>\gamma$, 
we also let $u^j$ be as in the lemma, $r_j \to 0$ be a sequence such that
$$
\bigg{\|} \frac{\u(r_j x)}{jr_j^{1+\alpha}}\bigg{\|}_{\tilde{L}^2(B_1^+)}=1.
$$
Such sequences exist since 
$$
\limsup_{r\to 0}\bigg{\|}\frac{\u(x)}{r^{1+\alpha}}\bigg{\|}_{\tilde{L}^2(B_r^+)}=\infty.
$$
We may also choose $r_j$ maximal in the sense that
\begin{equation}\label{controlofgrowthflatproof}
\|\u^j(x)\|_{\tilde{L}^2(B_r^+)}\le jr^{1+\alpha}
\end{equation}
for $r\ge r_j$. We make the blow-up 
$$
\v^j(x)=\frac{\u^j(r_j x)}{jr_j^{1+\alpha}}.
$$
Then, for a sub-sequence $\v^j\to \v^0$, locally and weakly in $W^{2,2}$
and locally strongly in $W^{1,2}$.

Also 
$$
\v^{i,j}\equiv \frac{\partial \v^j}{\partial x_i}
$$
will converge locally in $C^\beta\cap W^{1,2}$ to a solution $\v^{i,0}$
to the following mixed boundary value problem 
$$
\begin{array}{ll}
\Delta \v^{i,0} +\frac{\lambda+2}{2}\nabla \div (\v^{i,0})= 0 & \textrm{ in } \mathbb{R}^n_+ \\
e_n\cdot \v^{i,0}=0 & \textrm{ on } \Pi\cap \{e_n\cdot \v^{i,0}= 0\} \\
\frac{\partial e_n\cdot \v^{i,0}}{\partial x_n}+\frac{\lambda}{4} \textrm{div}(\v^{i,0})=0 & \textrm{ on }
\Pi \cap \{e_n\cdot \v^{i,0}> 0 \} \\
\frac{\partial e_k\cdot \v^{i,0}}{\partial x_n}+\frac{\partial e_n\cdot \v^{i,0}}{\partial x_k}=0 & \textrm{ on } \Pi \textrm{ for }k=1,2, ..., n-1.
\end{array}
$$
if $i=1,2, ..., n-1$ and 
$$
\begin{array}{ll}
\Delta \v^{n,0} +\frac{\lambda+2}{2}\nabla \div (\v^{n,0})= 0 & \textrm{ in } \mathbb{R}^n_+ \\
e_n\cdot \v^{n,0}=0 & \textrm{ on } \Pi\cap \{e_n\cdot \v^{i,0}> 0\} \\
\frac{\partial e_n\cdot \v^{i,0}}{\partial x_n}+\frac{\lambda}{4} \textrm{div}(\v^{i,0})=0 & \textrm{ on }
\Pi\cap \{e_n\cdot \v^{i,0}= 0\} \\
\frac{\partial e_k\cdot v^{i,0}}{\partial x_n}+\frac{\partial e_n\cdot \v^{i,0}}{\partial x_k}=0 & \textrm{ on } \Pi \textrm{ for }k=1,2,...,n-1.
\end{array}
$$

Using (\ref{controlofgrowthflatproof}) we may also conclude that
$$
\sup_{B_R^+}|\v^0|\le CR^{1+\alpha}.
$$
From the Corollary \ref{NEWCOR4} we can conclude that
$$
\v^0=\nabla \Big( 
\frac{\lambda+2}{2(\lambda+4)}\frac{\partial \tau}{\partial x_3}x_3-\frac{\lambda+3}{\lambda+4}\tau \Big)+e_3\frac{\partial \tau}{\partial x_3}
$$
where 
$$
\begin{array}{ll}
\Delta \tau =0 & \textrm{ in } \R^n_+ \\
\frac{\partial \tau}{\partial x_n}=0 & \textrm{ on } \Lambda_{\v^0} \\
\frac{\partial^2 \tau}{\partial x_n^2}=0 & \textrm{ on } \Omega_{\v^0} \\
\frac{\partial \tau}{\partial x_n} \in W^{2,2}(B_R^+) & \textrm{ for each }R>0.
\end{array}
$$
Naturally the function  
$$
\zeta(x)=\frac{\partial \tau}{\partial x_n}
$$
will solve
$$
\begin{array}{ll}
\Delta \zeta =0 & \textrm{ in } \R^3_+ \\
\zeta=0 & \textrm{ on } \Lambda_{\v^0} \\
\frac{\partial \zeta}{\partial x_n}=0 & \textrm{ on } \Omega_{\v^0} \\
\zeta \in W^{2,2}(\overline{B_R^+}) & \textrm{ for each }R>0.
\end{array}
$$
Moreover,
$$
\zeta_i(x)=\frac{\partial \zeta}{\partial x_i} \quad \textrm{ for }i=1,2,...,n-1 
$$
will solve
\begin{equation}\label{zetai}
\begin{array}{ll}
\Delta \zeta_i =0 & \textrm{ in } \R^3_+ \\
\zeta_i=0 & \textrm{ on } \Lambda_{\v^0} \\
\frac{\partial \zeta_i}{\partial x_n}=0 & \textrm{ on } \Omega_{\v^0} \\
\zeta \in W^{1,2}(B_R^+) & \textrm{ for each }R>0.
\end{array}
\end{equation}
and
\begin{equation}\label{zetaigrowth}
\sup_{B_R^+}|\zeta_i|\le CR^{\alpha}.
\end{equation}
By Benedicks' Theorem we know that the set of solutions to (\ref{zetai})
and (\ref{zetaigrowth}) is a one dimensional set. We may conclude that
there are constants $a_1,a_2,...,a_{n-1}$ such that
$$
a_1\zeta_1=a_2\zeta_2=...=a_{n-1}\zeta_{n-1}
$$
and therefore that
$$
\eta \cdot \nabla' \zeta =0
$$
for every 
$$
\eta\in \{\eta\in \mathbb{R}^n;\; \eta \cdot (a_1,a_2,...,a_{n-1},0)=0\}\cap \Pi.
$$

By rotating the coordinate system we may assume that
$$
\zeta(x)=\zeta(x_1,x_n).
$$
We can directly conclude that
$$
\tau(x)=\tilde{\tau}(x_1,x_n)+\bar{\tau}(x').
$$
Where $\Delta \tilde{\tau}=\Delta\bar{\tau}=0$. We thus have that,
as in the proof of Corollary \ref{taubisiszeroonOmega},
$$
\v^0=\nabla \Big( 
\frac{\lambda+2}{2}(\lambda+4)\frac{\partial \tau}{\partial x_n}x_n
-\frac{\lambda+3}{\lambda+4}\tau\Big)+e_n\frac{\partial \tau}{\partial x_n}
$$
$$
=\nabla \Big( 
\frac{\lambda+2}{2}(\lambda+4)\frac{\partial \tilde{\tau}}{\partial x_n}x_n
-\frac{\lambda+3}{\lambda+4}\tilde{\tau}\Big)+e_n\frac{\partial \tilde{\tau}}{\partial x_3}+\frac{\lambda+3}{\lambda+4}\nabla\bar{\tau}.
$$
But $\bar{\tau}$ is harmonic and
$$
\sup_{B_R}|\bar\tau|\le CR^{2+\alpha}.
$$
By the Liouville Theorem it follows that $\nabla \bar\tau$ is an affine 
function. But by our standing assumption the affine part of $\v^0$ is zero.
We may thus conclude that
\begin{equation}\label{dimred}
\tau(x)=\tau(x_1,x_n).
\end{equation}

\textbf{Claim:} {\sl Neither $\Lambda_{\v^0}$ nor $\Sigma_{\v^0}$ are empty.}

\textsl{Proof of the Claim:} Lets assume that $\Lambda_{\v^0}=\emptyset$
then
$$
\begin{array}{ll}
\Delta \v^{0} +\frac{\lambda+2}{2}\nabla \div (\v^{0})= 0 & \textrm{ in } \mathbb{R}^n_+ \\
e_n\cdot \v^{0}=0 & \textrm{ on } \Pi \\
\frac{\partial e_k\cdot v^{0}}{\partial x_n}=0 & \textrm{ on } \Pi \textrm{ for }k=1,2.
\end{array}
$$
Moreover
$$
\sup_{B_R^+}|\v^0|\le R^{1+\alpha}.
$$
In particular
$$
\w(x)=\left\{
\begin{array}{ll}
\v^0(x) & \textrm{ if } x_n>0 \\
\left[
\begin{array}{l}
e_1\cdot \v^0(x_1,x_2,..., -x_{n}) \\
e_2\cdot \v^0(x_1,x_2,..., -x_n) \\
\vdots \\
-e_n\cdot \v^0(x_1,x_2,..., -x_n)
\end{array}
\right] & \textrm{ if } x_n<0
\end{array}
\right.
$$
will solve
$$
\Delta \w+\frac{\lambda+2}{2}\nabla \div(\w)=0 \quad \textrm{ in } \R^3
$$
and
$$
\sup_{B_R}|\w|\le R^{1+\alpha}.
$$
It follows, from Liouville's Theorem, that $\w$ is a plane. 
This contradicts that $0\in \Gamma_\u$.

$ $

Using (\ref{dimred}) we may 
consider $\v^0$ as a solution in $\mathbb{R}^2_+$ and use Lemma 
\ref{Onepoint}. It follows that the free boundary consists of one point.
Extending $\v^0$ to $\mathbb{R}^n_+$ again we see that the free boundary is a
plane in $\Pi$. The proposition follows. \qed

\section{Almost Optimal Regularity}

We can now easily deduce that the solutions are $C^{1,\beta}$
for each $\beta<1/2$.

\begin{lem}\label{almostoptholderreg}
Let $\u$ be a solution to the Signorini problem in $B_1^+$ and assume that
$\|\u\|_{\tilde{L}^2(B_1^+)}=1$ and that $0\in \Gamma_\u$ then for each $\alpha< 1/2$ there 
exist a constant $C_\alpha$ such that
$$
\sup_{B_r^+}|\u|\le C_\alpha r^{1+\alpha}.
$$
\end{lem}
\textsl{Proof:} If not then we can find an $\alpha<1/2$ and a sequence of 
solutions 
$\u^j$ and $r_j\to 0$ such that
$$
\sup_{B_{r_j}^+}|\u^j|=jr_j^{1+\alpha},
$$
and 
$$
\sup_{B_R^+}|\u^j|\le jR^{1+\alpha}
$$
for each $R\ge r_j$. Make the blow-up
$$
\tilde{\u}^j=\frac{\u^j(r_jx)}{jr_j^{1+\alpha}}
$$
then, using Lemma \ref{c1beta}, for a sub-sequence $\tilde{\u}^j\to \u^0$ 
in $C^{1,\beta}_{loc}(\mathbb{R}^n_+)$.
From Proposition \ref{flatness} we can conclude that
$$
\u^0(x)=\u^0(x_1,x_n)
$$ 
and that after a rotation
and $\Sigma$ and $\Lambda$ are complementary half spaces, furthermore
we have $\sup_{B_R^+}|\u^0|\le R^{1+\alpha}$ for $R>1$. It follows, 
from Lemma \ref{eigenfinctforhalfspace}, that $\u^0=0$
which contradicts $\sup_{B_1^+}|\u^0|=\lim_{j\to \infty}\big(\sup_{B_1^+}|\tilde{\u}^j|\big)=1$.
\qed

\begin{cor}\label{addeddetailin6}
Let $\u$ be a solution to the Signorini  problem in $B_1^+$ then for 
each $\alpha<1/2$
there exists a constant $C_\alpha$ such that
$$
\|\u\|_{C^{1,\alpha}(B_{1/2}^+)}\le C_\alpha \| \u \|_{\tilde{L}^2(B_1^+)}.
$$
\end{cor}
\textsl{Proof:} We may assume that $\|\u\|_{L^2}=1$. 
Let $d(x^0)=\textrm{dist}(x^0,\Gamma)$ then
$$
\sup_{B_{d(x^0)}(x^0)}|\u|\le C_\alpha d(x^0)^{1+\alpha}.
$$
Thus
$$
\u_{d(x^0)}=\frac{\u(d(x^0)x+x^0)}{d(x^0)^{1+\alpha}}
$$
is a solution in $B_1$ and $\sup_{B_1}|\u_{d(x^0)}|\le C_\alpha$. We might need,
and in that case we do, either evenly or oddly reflect $\u$ across $\Pi$ in 
order for $\u_{d(x^0)}$ to be defined in the entire unit ball. It follows
that $|\nabla \u_{d(x^0)}(0)|\le CC_\alpha$. Scaling back we get
$|\nabla \u(x)|\le CC_\alpha d(x)^\alpha$ which implies that $\u\in C^{1,\alpha}$.
\qed

\section{Fundamental and Technical Results.}\label{fundsec}

With this section we start to get a little more technical and we will
start to lay the foundation for the flatness improvement results that leads to 
optimal regularity and free boundary regularity.

\begin{lem}\label{subtraction}
Let $\u^j$ be a sequence of solutions to the Signorini 
problem in $B_1^+$ and 
\begin{equation}\label{starinhell}
\inf_{\xi\in \Pi}\|\u-\p_{3/2}^\xi\|_{\tilde{L}^2(B_1^+)}=
\|\u-\p_{3/2}\|_{\tilde{L}^2(B_1^+)}=\delta_j\to 0
\end{equation}
assume furthermore that
$$
\v^j(x)=\frac{\u(x)-\p_{3/2}}{\delta_j}\to \v^0 \textrm{ weakly in }L^2
$$
then $\v^0$ solves
$$
\begin{array}{ll}
\Delta \v^0+ \frac{\lambda+2}{2}\nabla\div(\v^0) =0 & \textrm{in }B_1^+ \\
e_n\cdot \v^{0}=0 & \textrm{ on } \Pi\cap \{x_1 > 0\} \\
\frac{\partial e_n\cdot \v^{0}}{\partial x_n}+\frac{\lambda}{4} \textrm{div}(\v^{0})=0 & \textrm{ on }
\Pi\cap \{ x_1 >0 \} \\
\frac{\partial e_k\cdot \v^{0}}{\partial x_n}+\frac{\partial e_n\cdot \v^{0}}{\partial x_k}=0 & \textrm{ on } \Pi \textrm{ for }k=1,2,...,n-1,
\end{array}
$$
\begin{equation}\label{twostarinhell}
\inf_{a\in \mathbb{R}}\|\v^0-a\p_{3/2}\|_{\tilde{L}^2(B_1^+)}=
\|\v^0\|_{\tilde{L}^2(B_1^+)}.
\end{equation}
\end{lem}
\textsl{Proof:} That $\v^0$ converges weakly to a solution of the
system is simple so we will only prove (\ref{twostarinhell}).

By (\ref{starinhell}) we have
$$
\int_{B_1^+}\p_{3/2}\big( \u^j-\p_{3/2}\big)=0
$$
so by weak convergence we have
\begin{equation}\label{threestarinhell}
\int_{B_1^+}\p_{3/2}\v^0=0.
\end{equation}
Now $\|\v^0-a\p_{3/2}\|_{\tilde{L}^2(B_1^+)}$ is convex in $a$ so it has only one
minimum. Therefore (\ref{threestarinhell}) implies that
the minimum is at $a=0$. \qed

\begin{lem}\label{assinpropimplyassinprop}
Let $\u$ be a global solution to the Signorini problem and
$$
\inf_{\xi\in \Pi}\| \u-\p_{3/2}^\xi \|_{\tilde{L}^2(B_R^+)}\le \mu \|\u\|_{\tilde{L}^2(B_R^+)}
$$
for all $R\ge 1$. Then if $\mu$ is small enough then $\u=\p_{3/2}^\xi$ for some $\xi\in \Pi$.

Also, if $\u$ is a solution to the Signorini problem in $B_1^+$ then for each 
$\epsilon>0$ there exist a $\mu_\epsilon>0$ such that if
$$
\inf_{\xi\in \Pi}\| \u-\p_{3/2}^\xi \|_{\tilde{L}^2(B_r^+)}\le \mu_\epsilon \|\u\|_{\tilde{L}^2(B_r^+)}
$$
for all $r\le 1$. Then  
$$
\|\u\|_{\tilde{L}^2(B_r)}\ge r^{3/2+\epsilon}.
$$
\end{lem}
\textsl{Proof:} The proof of the first and second parts are very 
similar so we will only prove the first part.

Let $\gamma>0$ be the real non-negative solution to $\gamma^{(n+3)/2}=1/2$. Then
$$
\mu \|\u\|_{\tilde{L}^2(B_R)}\ge \gamma^{n/2}\|\u-\p_{3/2}^{\xi_R}\|_{\tilde{L}^2(B_{\gamma R}^+)}\ge
$$
\begin{equation}\label{mutimesnormestimate}
\gamma^{n/2}\|\p_{3/2}^{\xi_R}-\p_{3/2}^{\xi_{\gamma R}}\|_{\tilde{L}^2(B_{\gamma R}^+)}-
\gamma^{n/2}\|\u-\p_{3/2}^{\xi_{\gamma R}}\|_{\tilde{L}^2(B_{\gamma R}^+)}\ge
\end{equation}
$$
\gamma^{n/2}\|\p_{3/2}^{\xi_R}-\p_{3/2}^{\xi_{\gamma R}}\|_{\tilde{L}^2(B_{\gamma R}^+)}-\mu \gamma^{n/2}\|\u\|_{\tilde{L}^2(B_{\gamma R}^+)}.
$$
Next we notice that 
$$
\gamma^{n/2}\|\p_{3/2}^{\xi_R}-\p_{3/2}^{\xi_{\gamma R}}\|_{\tilde{L}^2(B_{\gamma R}^+)}\ge 
c\gamma^{3/2}R^{3/2} |\xi_R- \xi_{\gamma R}|.
$$
Inserting this into (\ref{mutimesnormestimate}) results in
\begin{equation}\label{mutermfrombelow}
\mu \Big( \|\u\|_{\tilde{L}^2(B_R)} +\gamma^{n/2}\|\u\|_{\tilde{L}^2(B_{\gamma R}^+)} \Big)\ge
\gamma^{(n+3)/2}\big| \xi_R- \xi_{\gamma R} \big|.
\end{equation}
Next we use the triangle inequality to estimate for any $T\ge 1$
$$
\|\u\|_{\tilde{L}^2(B_T)}- \|\p_{3/2}^{\xi_T}\|_{\tilde{L}^2(B_T^+)}\le 
\|\u-\p_{3/2}^{\xi_T}\|_{\tilde{L}^2(B_T^+)}\le \mu\|\u\|_{\tilde{L}^2(B_T)}.
$$
That is
\begin{equation}\label{uestimatedbyp32withmu}
\|\u\|_{\tilde{L}^2(B_T)}\le \frac{1}{1-\mu}\|\p_{3/2}^{\xi_T}\|_{\tilde{L}^2(B_T^+)}.
\end{equation}
If we use (\ref{uestimatedbyp32withmu}) in (\ref{mutermfrombelow}) and that $\gamma^{(n+3)/2}=1/2$
we get
$$
\frac{\mu}{1-\mu}\big( 2+\sigma  \big)\ge(1-\sigma)
$$
where we have used the notation $\sigma |\xi_R|=|\xi_{\gamma R}|$. That is
$$
\sigma\ge 1-3\mu.
$$
We have shown that
$$
|\xi_{R}|\le \frac{1}{1-3\mu}|\xi_{\gamma R}|.
$$
Iterating this relation we get
$$
|\xi_{\gamma^{-k}}|\le \big(1-3\mu\big)^{-k}|\xi_1|.
$$
We may normalize so $|\xi_1|=1$ and use (\ref{uestimatedbyp32withmu}) to deduce that
\begin{equation}\label{niostjarnor}
\|\u\|_{\tilde{L}^2(B_{\gamma^{-k}}^+)}\le \frac{1}{1-\mu}\big( 1-3\mu\big)^{-k}\gamma^{-\frac{3k}{2}}.
\end{equation}
In particular if $\mu$ is small enough to satisfy
\begin{equation}\label{tiostjarnor}
(n+3)\frac{\ln\big( \frac{1}{1-3\mu}\big)}{\ln(2)}<1
\end{equation}
Then
\begin{equation}\label{growthatinfinityinsomelemma}
\lim_{R\to \infty}\frac{\ln\big( \|\u\|_{\tilde{L}^2(B_R^+)}\big)}{\ln(R)}<2.
\end{equation}
We may conclude, as in the argument of Proposition \ref{flatness} that $\u(x)=\u(x_1,x_n)$
in some coordinate system. From Lemma \ref{eigenfinctforhalfspace} we conclude that
$$
\u= a \p_{3/2}(x)+\sum_{i=2}^\infty a_i \q_i.
$$
From (\ref{growthatinfinityinsomelemma}) it follows that $a_i=0$ for all $i$ 
and the first part of the Lemma follows. The second part is similar. \qed

\begin{lem}\label{technicallemmaforrotest}
Assume that $\|\u\|_{W^{1,2}(B_1^+)}\le C_1$ and
$$
\inf_{\xi\in\Pi}\|\u-\p_{3/2}^\xi\|_{\tilde{L}^2(B_1^+)}\le \delta,
$$
where the minimizing $\xi$ satisfies $|\xi|=1$ assume furthermore that
$\delta$ and $\delta_1$ are small enough and that
\begin{equation}\label{starfordelta1}
\|\nabla'' \u\|_{\tilde{L}^2(B_1^+)}\le \delta_1.
\end{equation}
Then $|\xi-e_1|\le C\delta_1$ and
$$
\|\nabla_\xi'' \u\|_{\tilde{L}^2(B_1)}\le \big( CC_1 + 1 \big)\delta_1,
$$
here $\nabla_\xi''=\nabla -e_n\frac{\partial }{\partial x_n}-\xi_w/|\xi_w|^2
(\xi_w\cdot \nabla )$ is the gradient restricted to the subspace orthogonal to 
$e_n$ and $\xi$.
\end{lem}
\textsl{Proof:} We may rotate the coordinates so that 
$\xi=(\xi_1,\xi_2,0,0,...)$ where $\xi_2$ is very small, we
may also assume that $\xi_2\ge 0$.
Notice that when $\xi_2$ and $\delta$ are small then
\begin{equation}\label{derivativeinxi2}
\frac{\p_{3/2}^\xi-\p_{3/2}}{\delta}=\frac{\xi_2}{\delta}x_2\p_{1/2}+\q
\end{equation}
where 
$$
\|\q\|_{\tilde{L}^2(B_1^+)}=o\big(\xi_2/\delta\big).
$$
Let 
$$
\v(x)=\frac{\u(x)-\p_{3/2}^\xi(x)}{\delta}
$$
and
$$
\w(x)=\frac{\u(x)-\p_{3/2}(x)}{\delta}.
$$
Then by the minimizing property of $\p_{3/2}^\xi$ it follows that
\begin{equation}\label{statingminimalityofpxiintermsofv}
\|\v\|_{\tilde{L}^2(B_1^+)}\le \|\w\|_{\tilde{L}^2(B_1^+)}.
\end{equation}
From the assumption (\ref{starfordelta1}) it follows that 
$\|\nabla'' \w\|_{\tilde{L}^2(B_1^+)}\le \delta_1/\delta$ which by the
Poincare inequality implies that
\begin{equation}\label{poincareforwnablabis}
\|\w-\bar{\w}\|_{\tilde{L}^2(B_1^+)}\le C\frac{\delta_1}{\delta}
\end{equation}
where $\bar{\w}(y)=\bar{\w}(y_1,y_n)$ is the average of $\w$ on the 
$(n-1)$-plane $B_1^+\cap\{x_1=y_1,x_n=y_n\}$.
Form  (\ref{derivativeinxi2}) it follows that
$$
\int_{B_1^+}|\v|^2= \int_{B_1^+}\Big( |\w|^2+ 2\frac{\xi_2}{\delta}x_2\p_{1/2}\w+
\frac{(\xi_2)^2}{\delta^2}x_2^2|\p_{1/2}|^2\Big)+ o\big(\xi_2/\delta \big)\ge
$$
$$
\int_{B_1^+} |\w|^2+c\frac{(\xi_2)^2}{\delta^2}-2\frac{\xi_2}{\delta}
\sqrt{\int_{B_1^+}x_2^2\p_{1/2}^2}\sqrt{\int_{B_1^+}|\w-\bar{\w}|^2}
$$
where we have used that $\int_{B_1^+} x_2 \bar{\w}=0$. Applying 
(\ref{poincareforwnablabis}) and (\ref{statingminimalityofpxiintermsofv}) we 
may deduce that
$$
C\delta_1\ge \xi_2
$$
which implies the first conclusion in the Lemma.

The second conclusion follows by writing $\nabla''_\xi=\nabla''+\xi_2e_2\frac{\partial }{\partial x_2}+(1-\sqrt{1+\xi_2^2})e_1\frac{\partial }{\partial x_1}$
and thus
$$
\|\nabla''_\xi \u\|_{\tilde{L}^2(B_1^+)}\le \|\nabla'' \u\|_{\tilde{L}^2(B_1^+)}+
C|\xi_2|\|\nabla \u\|_{\tilde{L}^2(B_1^+)}\le \delta_1+CC_1\delta_1.
$$
\qed

The next Lemma looks more complicated than it is in reality. It just 
states that a gradient can not be written as something that is
not a gradient plus a small perturbation.

\begin{lem}\label{fundamentalthinobstlemma}
Let $\u$ solve the Signorini problem in $B_1^+$ and assume that
$\u=\p_{3/2}+R$ where 
\begin{equation}\label{ambushbycurl}
\nabla R=
\left[ \begin{array}{lll}
a_1p_{1/2}^1 & a_2p_{1/2}^1 & a_3\tilde{p}_{1/2}^1 \\
a_1p_{1/2}^2 & a_2p_{1/2}^2 & a_3\tilde{p}_{1/2}^2  \\
a_1p_{1/2}^3 & a_2p_{1/2}^3 & a_3\tilde{p}_{1/2}^3 
\end{array}\right]+
\left[ \begin{array}{lll}
m_{11}(x) & m_{12}(x) & m_{13}(x) \\
m_{21}(x) & m_{22}(x) & m_{23}(x) \\
m_{31}(x) & m_{32}(x) & m_{33}(x) 
\end{array}\right]
=
P+\tilde{M}
\end{equation}
and
\begin{equation}\label{1stcondlem13}
\|\tilde{M}\|_{L^2(B_1^+\cap\{x_3>\epsilon\})}\le \epsilon c \sqrt{\sum_{i=1}^3(a^i)^2}
\end{equation}
for any $\epsilon\in [0,1/4)$ and some small enough, but universal, $c$.
Assume furthermore that
\begin{equation}\label{2ndcondlem13}
|a^1-a^3|> \frac{1}{4}\big( |a^1|+|a^3|\big).
\end{equation}

Then $a^1=a^2=a^3=0$.
\end{lem}
\textsl{Proof:} If we apply the curl operator on both sides of 
(\ref{ambushbycurl}) we can deduce that
$$
0=\textrm{curl}\big( \nabla R\big)=\textrm{curl}\Bigg(
\left[ \begin{array}{lll}
a_1p_{1/2}^1 & a_2p_{1/2}^1 & a_3\tilde{p}_{1/2}^1 \\
a_1p_{1/2}^2 & a_2p_{1/2}^2 & a_3\tilde{p}_{1/2}^2  \\
a_1p_{1/2}^3 & a_2p_{1/2}^3 & a_3\tilde{p}_{1/2}^3 
\end{array}\right] \Bigg)
$$
$$
+\textrm{curl}\Bigg( \left[ \begin{array}{lll}
m_{11}(x) & m_{12}(x) & m_{13}(x) \\
m_{21}(x) & m_{22}(x) & m_{23}(x) \\
m_{31}(x) & m_{32}(x) & m_{33}(x) 
\end{array}\right] \Bigg).
$$
Rearranging terms and taking the 
$\tilde{L}^2(B_{1-2\epsilon}\cap \{x_3>2\epsilon\})$
norm on both sides we may conclude that
\begin{equation}\label{firstforveccurl}
\left\|\textrm{curl}\Bigg( 
\left[ \begin{array}{lll}
a_1p_{1/2}^1 & a_2p_{1/2}^1 & a_3\tilde{p}_{1/2}^1 \\
a_1p_{1/2}^2 & a_2p_{1/2}^2 & a_3\tilde{p}_{1/2}^2  \\
a_1p_{1/2}^3 & a_2p_{1/2}^3 & a_3\tilde{p}_{1/2}^3 
\end{array}\right]\Bigg)
\right\|_{\tilde{L}^2(B_{1-2\epsilon}\cap \{x_3>2\epsilon\})}=
\|\textrm{curl}(\tilde{M})\|_{\tilde{L}^2(B_{1-2\epsilon}\cap \{x_3>2\epsilon\})}.
\end{equation}
But
$$
\left\|\textrm{curl}\Bigg( 
\left[ \begin{array}{lll}
a_1p_{1/2}^1 & a_2p_{1/2}^1 & a_3\tilde{p}_{1/2}^1 \\
a_1p_{1/2}^2 & a_2p_{1/2}^2 & a_3\tilde{p}_{1/2}^2  \\
a_1p_{1/2}^3 & a_2p_{1/2}^3 & a_3\tilde{p}_{1/2}^3 
\end{array}\right]\Bigg)
\right\|_{\tilde{L}^2(B_{1-2\epsilon}\cap \{x_3>2\epsilon\})}=
$$
$$
=\left\|
\left[ \begin{array}{lll}
-a_2\frac{\partial p_{1/2}^1}{\partial x_3} & a_1\frac{\partial p_{1/2}^1}{\partial x_3}-a_3\frac{\partial \tilde{p}_{1/2}^1}{\partial x_1} & a_2\frac{\partial p_{1/2}^1}{\partial x_1}  \\ 
-a_2\frac{\partial p_{1/2}^2}{\partial x_3} & a_1\frac{\partial p_{1/2}^2}{\partial x_3}-a_3\frac{\partial \tilde{p}_{1/2}^2}{\partial x_1} & a_2\frac{\partial p_{1/2}^2}{\partial x_1}  \\ 
-a_2\frac{\partial p_{1/2}^3}{\partial x_3} & a_1\frac{\partial p_{1/2}^3}{\partial x_3}-a_3\frac{\partial \tilde{p}_{1/2}^3}{\partial x_1} & a_2\frac{\partial p_{1/2}^3}{\partial x_1}   
\end{array}\right]
\right\|_{\tilde{L}^2(B_{1-2\epsilon}\cap \{x_3>2\epsilon\})},
$$
where we have used that $\p_{1/2}(x)=\p_{1/2}(x_1,x_3)$. Next we notice that
by definition
$$
\tilde{\p}_{1/2}=\frac{\partial \p_{3/2}}{\partial x_3}
$$
and that
$$
\p_{1/2}=\frac{\partial \p_{3/2}}{\partial x_1}.
$$
In particular this implies that
$$
\frac{\partial p_{1/2}^i}{\partial x_3}=\frac{\partial \tilde{p}_{1/2}^i}{\partial x_1}
$$
for $i=1,2,3$. Therefore
\begin{equation}\label{secondforveccurl}
\left\|\textrm{curl}\Bigg( 
\left[ \begin{array}{lll}
a_1p_{1/2}^1 & a_2p_{1/2}^1 & a_3\tilde{p}_{1/2}^1 \\
a_1p_{1/2}^2 & a_2p_{1/2}^2 & a_3\tilde{p}_{1/2}^2  \\
a_1p_{1/2}^3 & a_2p_{1/2}^3 & a_3\tilde{p}_{1/2}^3 
\end{array}\right]\Bigg)
\right\|_{\tilde{L}^2(B_{1-2\epsilon}\cap \{x_3>2\epsilon\})}
\end{equation}
$$
=\left\|
\left[ \begin{array}{lll}
-a_2\frac{\partial p_{1/2}^1}{\partial x_3} & (a_1-a_3)\frac{\partial^2 p_{3/2}^1}{\partial x_1\partial x_3} & a_2\frac{\partial p_{1/2}^1}{\partial x_1}  \\ 
-a_2\frac{\partial p_{1/2}^2}{\partial x_3} & (a_1-a_3)\frac{\partial^2 p_{3/2}^2}{\partial x_1\partial x_3} & a_2\frac{\partial p_{1/2}^2}{\partial x_1}  \\ 
-a_2\frac{\partial p_{1/2}^3}{\partial x_3} & (a_1-a_3)\frac{\partial^2 p_{3/2}^3}{\partial x_1\partial x_3} & a_2\frac{\partial p_{1/2}^3}{\partial x_1}   
\end{array}\right]
\right\|_{\tilde{L}^2(B_{1-2\epsilon}\cap \{x_3>2\epsilon\})}
$$
$$
\ge c_0\sqrt{a_2^2\int_{B_{1-\epsilon}\cap \{x_3>\epsilon\} }\big|\nabla \p_{1/2} \big|^2+(a_1-a_3)^2\int_{B_{1-\epsilon}\cap \{x_3>\epsilon\} } \Big( \frac{\partial^2 \p_{3/2}}{\partial x_1\partial x_3}\Big)^2}
$$
$$
\ge c_1\sqrt{a_2^2+(a_1-a_3)^2}
\ge c_2 \sqrt{a_1^2+a_2^2+a_3^2},
$$
where we have used (\ref{2ndcondlem13}) in the last inequality.
Putting (\ref{firstforveccurl}) and (\ref{secondforveccurl}) together
we have thus shown that
$$
c_2 \sqrt{a_1^2+a_2^2+a_3^2}\le \|\textrm{curl}(\tilde{M})\|_{\tilde{L}^2(B_{1-2\epsilon}\cap \{x_3>2\epsilon\})}.
$$
Next we use that each column of $\tilde{M}$ solves the Lame system and therefore
we have the estimate
$$
\|\textrm{curl}(\tilde{M})\|_{\tilde{L}^2(B_{1-2\epsilon}\cap \{x_3>2\epsilon\})}
\le
C \|\nabla\tilde{M}\|_{\tilde{L}^2(B_{1-2\epsilon}\cap \{x_3>2\epsilon\})}\le
\frac{C}{\epsilon} \|\tilde{M}\|_{\tilde{L}^2(B_{1-\epsilon}\cap\{x_3>\epsilon \})}
$$
so we have shown that
$$
\sqrt{a_1^2+a_2^2+a_3^2}\le \frac{C}{c_2\epsilon} \|\tilde{M}\|_{\tilde{L}^2(B_{1}^+)}\le \frac{Cc}{c_2}\sqrt{a_1^2+a_2^2+a_3^2},
$$
where we have used (\ref{1stcondlem13}). But if $c\le c_2/(2C)$ this implies 
that
$$
\sqrt{a_1^2+a_2^2+a_3^2}\le \frac{1}{2}\sqrt{a_1^2+a_2^2+a_3^2}
$$ 
which implies that
$$
\sqrt{a_1^2+a_2^2+a_3^2}=0.
$$
\qed

\section{Regularity Improvement for the Tangential Gradient.}

In the next proposition, and also afterwards, we will denote 
$$
\nabla''=\big(0,\frac{\partial}{\partial x_2},\frac{\partial}{\partial x_3},...,\frac{\partial}{\partial x_{n-1}} ,0\big).
$$

\begin{lem}\label{localgainforobst}
Let $\u$ solve the Signorini problem in $B_1^+$ and $0<\gamma<1/2$ 
then there exists a $\delta_\gamma>0$ such that if 
\begin{equation}\label{odensdog}
\inf_{\xi\in\Pi}\|\u-\p_{3/2}^\xi\|_{\tilde{L}^2(B_1^+)}=\|\u-\p_{3/2}\|_{\tilde{L}^2(B_1^+)}
\le \delta_\gamma,
\end{equation}
and
$$
\|\nabla'' \u\|_{\tilde{L}^2(B_1^+)}\le \delta_\gamma
$$
then there exist an $0<s_\gamma<1$ such that
$$
\|\nabla'' \u\|_{\tilde{L}^2(B_{s_\gamma}^+)}\le 
s_{\gamma}^{1/2+\gamma}\|\nabla'' \u\|_{\tilde{L}^2(B_1^+)}.
$$
\end{lem}
\textsl{Proof:} Let $\u^j$ be a sequence as in the Lemma corresponding
to $\delta_j\to 0$ and a fixed $\gamma\in (0,1)$. Consider
$$
\left[
\begin{array}{l}
\v^{1,j} \\
\v^{2,j} \\
\v^{3,j} \\
\vdots \\
\v^{n-1,j} \\
\v^{n,j}
\end{array}\right]=
\left[\begin{array}{l}
\frac{1}{\|\nabla (\u^j-\p_{3/2})\|_{\tilde{L}^2(B_1^+)}}\frac{\partial \u^j-\p_{3/2}}{\partial x_1} \\
\frac{1}{\|\nabla'' \u^j\|_{\tilde{L}^2(B_1^+)}}\frac{\partial \u^j}{\partial x_2} \\
\frac{1}{\|\nabla'' \u^j\|_{\tilde{L}^2(B_1^+)}}\frac{\partial \u^j}{\partial x_3} \\
\vdots \\
\frac{1}{\|\nabla'' \u^j\|_{\tilde{L}^2(B_1^+)}}\frac{\partial \u^j}{\partial x_{n-1}} \\
\frac{1}{\|\nabla (\u^j-\p_{3/2})\|_{\tilde{L}^2(B_1^+)}}\frac{\partial \u^j-\p_{3/2}}{\partial x_n} 
\end{array}\right]
$$
that is $\v^{i,j}=\frac{1}{\|\nabla'' \u^j\|_{\tilde{L}^2(B_1^+)}}\frac{\partial \u^j}{\partial x_i}$
for $i=2,...,n-1$ and as in the displayed formula when $i=1$ or $i=n$.

Then for $i=2,3,...,n-1$ the function $v^i$ solves
\begin{equation}\label{flatthinlikeequation}
\begin{array}{ll}
\Delta \v^{i,j} +\frac{\lambda+2}{2}\nabla \div(\v^{i,j})=0 &\textrm{ in } B_1^+ \\
e_3\cdot \v^{i,j}=0 & \textrm{ on } \Lambda_{u^j}\to \{x_1>0\} \\

\frac{\partial e_3\cdot \v^{0}}{\partial x_3}+\frac{\lambda}{4} \textrm{div}(\v^{0})=0 & \textrm{ on } \Lambda_{u^j}\to \{x_1>0\} \\
\frac{\partial e_k\cdot v^{0}}{\partial x_3}+\frac{\partial e_3\cdot \v^{0}}{\partial x_k}=0 & \textrm{ on } \Sigma_{\u^j}\to \{x_1<0\} \textrm{ for }k=1,2
\end{array}
\end{equation}
and $\| \v^{i,j}\|_{L^\infty(B_{7/8}^+)},\|\v^{i,j}\|_{W^{1,2}(B_{7/8})}\le C$. 
We may also assume that $\lim_{j\to \infty}\v^{i,j}=\v^i$ by taking a sub-sequence if necessary.
By Lemma \ref{eigenfinctforhalfspace}
it follows that
\begin{equation}\label{someboringrep}
\v^i=\sum_{k=0}^\infty a_k^i \q_i
\end{equation}
where $\q_k$ are homogeneous of order $k/2$.

Also, since $\|\v^{i,j}\|_{\tilde{L}^2}\le 1$, we get that $\v^{1,j}$ 
converge weakly in $L^2$ to a solution $\v^1$  of 
(\ref{flatthinlikeequation}). 

Since $\v^{i,j}$ convergence strongly in $L^2(B_{1-\epsilon}^+)$
and $\v^{1,j}$ and $\v^{n,j}$ converges strongly in $L^2(B_{1-\epsilon}^+\cap\{x_n>\epsilon\})$
for each $\epsilon>0$ we may, in the set $B_{1-\epsilon}^+\cap\{x_n>\epsilon\}$, write
$\u^j=\p_{3/2}+R^j$ where
$$
\nabla R^j=
\left[ \begin{array}{l}
\|\nabla(\u^j-\p_{3/2})\|_{\tilde{L}^2(B_1)}\sum_{k=0}^\infty a_k^1\q_k \\
\|\nabla'' \u^j\|_{\tilde{L}^2(B_1)}\sum_{k=0}^\infty a_k^2\q_k \\
\|\nabla'' \u^j\|_{\tilde{L}^2(B_1)}\sum_{k=0}^\infty a_k^3\q_k \\
\vdots \\
\|\nabla'' \u^j\|_{\tilde{L}^2(B_1)}\sum_{k=0}^\infty a_k^{n-1}\q_k \\
\|\nabla( \u^j-\p_{3/2})\|_{\tilde{L}^2(B_1)}\sum_{k=0}^\infty a_k^n\q_k \\
\end{array}\right]+
\left[ \begin{array}{l}
o\big(\|\nabla(\u^j-\p_{3/2})\|_{\tilde{L}^2(B_1\cap \{x_n>\epsilon\})}\big)\\
o\big(\|\nabla'' \u^j\|_{\tilde{L}^2(B_1\cap \{x_n>\epsilon\})} \big)\\
o\big(\|\nabla'' \u^j\|_{\tilde{L}^2(B_1\cap \{x_n>\epsilon\})} \big)\\
\vdots \\
o\big(\|\nabla'' \u^j\|_{\tilde{L}^2(B_1\cap \{x_n>\epsilon\})} \big)\\
o\big(\|\nabla( \u^j-\p_{3/2})\|_{\tilde{L}^2(B_1\cap \{x_n>\epsilon\})}\big) \\
\end{array}\right]
$$ 
where the little-$o$ terms are considered to be little-$o$ in $L^2$ norm,
that is $f^j(x)=o(\|g^j\|_{\tilde{L}^2\cap \{x_n>\epsilon\}})$ if 
$\|f^j\|_{\tilde{L}^2(B_1\cap \{x_n>\epsilon\})}/\|g^j\|_{\tilde{L}^2\cap \{x_n>\epsilon\}}\to 0$ as $j\to \infty$.

Now consider $\u_s^j=\u^j(sx)/\sqrt{s}$ then $\u_s^j=s\p_{3/2}+R^j_s$ where
$$
\nabla R^j_s=
\left[ \begin{array}{l}
\|\nabla(\u^j-\p_{3/2})\|_{\tilde{L}^2(B_1)}a_0^1 \p_{1/2} \\
\|\nabla'' \u^j\|_{\tilde{L}^2(B_1)}a_0^2  \p_{1/2} \\
\|\nabla'' \u^j\|_{\tilde{L}^2(B_1)}a_0^3  \p_{1/2} \\
\vdots \\
\|\nabla'' \u^j\|_{\tilde{L}^2(B_1)}a_0^{n-1}  \p_{1/2} \\
\|\nabla( \u^j-\p_{3/2})\|_{\tilde{L}^2(B_1)}a_0^n \tilde{\p}_{1/2}
\end{array}\right]+
$$
$$
s\left[ \begin{array}{l}
\|\nabla(\u^j-\p_{3/2})\|_{\tilde{L}^2(B_1)}\sum_{k=1}^\infty s^{k-1} a_k^1\q_k \\
\|\nabla'' \u^j\|_{\tilde{L}^2(B_1)}\sum_{k=1}^\infty s^{k-1} a_k^2\q_k \\
\|\nabla'' \u^j\|_{\tilde{L}^2(B_1)}\sum_{k=1}^\infty s^{k-1} a_k^3\q_k \\
\vdots \\
\|\nabla'' \u^j\|_{\tilde{L}^2(B_1)}\sum_{k=1}^\infty s^{k-1} a_k^{n-1}\q_k \\
\|\nabla( \u^j-\p_{3/2})\|_{\tilde{L}^2(B_1)}\sum_{k=1}^\infty s^{k-1} a_k^n\q_k \\
\end{array}\right]+
s\left[ \begin{array}{l}
o\big(\|\nabla(\u^j-\p_{3/2})\|_{\tilde{L}^2(B_1\cap \{x_n>\epsilon\})}\big)\\
o\big(\|\nabla'' \u^j\|_{\tilde{L}^2(B_1\cap \{x_n>\epsilon\})} \big)\\
o\big(\|\nabla'' \u^j\|_{\tilde{L}^2(B_1\cap \{x_n>\epsilon\})} \big)\\
\vdots \\
o\big(\|\nabla'' \u^j\|_{\tilde{L}^2(B_1\cap \{x_n>\epsilon\})} \big)\\
o\big(\|\nabla( \u^j-\p_{3/2})\|_{\tilde{L}^2(B_1\cap \{x_n>\epsilon\})}\big) \\
\end{array}\right].
$$
We make the following claim.

\textbf{Claim:} {\sl Either $\|\nabla'' \u^j\|_{L^2(B_1^+)}=0$ or
$a_0^i=0$ for $i=1,...,n$ or both.}

\textsl{Proof of the claim:} We may assume that $j$ is very large and
$s$ very small. Moreover we may disregard the $o$-terms since they 
vanish in the limit. Writing
$$
P^j=\left[
\begin{array}{l}
\|\nabla(\u^j-\p_{3/2})\|_{\tilde{L}^2(B_1^+)}a_0^1 \p_{1/2} \\
\|\nabla'' \u^j\|_{\tilde{L}^2(B_1^+)}a_0^2  \p_{1/2} \\
\|\nabla'' \u^j\|_{\tilde{L}^2(B_1^+)}a_0^3  \p_{1/2} \\
\vdots \\
\|\nabla'' \u^j\|_{\tilde{L}^2(B_1^+)}a_0^{n-1}  \p_{1/2} \\
\|\nabla( \u^j-\p_{3/2})\|_{\tilde{L}^2(B_1^+)}a_0^n \tilde{\p}_{1/2}
\end{array}
\right]
$$
and
$$
\tilde{R}^j_s=
s\left[ \begin{array}{l}
\|\nabla(\u^j-\p_{3/2})\|_{\tilde{L}^2(B_1^+)}\sum_{k=1}^\infty s^{k-1} a_k^1\q_k \\
\|\nabla'' \u^j\|_{\tilde{L}^2(B_1^+)}\sum_{k=1}^\infty s^{k-1} a_k^2\q_k \\
\|\nabla'' \u^j\|_{\tilde{L}^2(B_1^+)}\sum_{k=1}^\infty s^{k-1} a_k^3\q_k \\
\vdots \\
\|\nabla'' \u^j\|_{\tilde{L}^2(B_1^+)}\sum_{k=1}^\infty s^{k-1} a_k^{n-1}\q_k \\
\|\nabla( \u^j-\p_{3/2})\|_{\tilde{L}^2(B_1^+)}\sum_{k=1}^\infty s^{k-1} a_k^n\q_k \\
\end{array}\right].
$$
Then 
$$
\nabla R_s^j=p^j+\tilde{R}^j_s+l.o.t.
$$
where $l.o.t.$ indicates lower order terms:that is terms of order 
$s^{3/2}, s^2,...$.
We hope to be able to use Lemma \ref{fundamentalthinobstlemma}.

First we need to show that the assumption (\ref{2ndcondlem13}) holds. Since 
$P^j$ is a gradient modulo lower order terms we may conclude that
$$
\frac{\partial \big( e_1\cdot P^j\big)}{\partial x_n}=\frac{\partial \big( e_n\cdot P^j\big)}{\partial x_1}
$$
which implies that $a_0^1=a_0^n$. Therefore (\ref{2ndcondlem13}) holds.

Next we need to verify that (\ref{1stcondlem13}) holds for some 
$\epsilon\in (0,1/4)$, that is
\begin{equation}\label{forLemma81}
\|\tilde{R}^j_s\|_{L^2(B_1^+)}\le c\bigg( \|\nabla(\u^j-\p_{3/2})\|_{\tilde{L}^2(B_1^+)}^2(a_0^1)^2
\end{equation}
$$
+\sum_{i=2}^{n-1}\|\nabla'' \u^j\|_{\tilde{L}^2(B_1)}^2(a_0^i)^2+ 
\|\nabla( \u^j-\p_{3/2})\|_{\tilde{L}^2(B_1)}^2(a_0^n)^2 \bigg)^{1/2}.
$$
If (\ref{forLemma81}) holds for some $s>0$ then $a^1_0=a^2_0=...=a^n_0=0$ by Lemma 
\ref{fundamentalthinobstlemma}.

If (\ref{forLemma81}) is not true then
\begin{equation}\label{somethingbeautiful}
\|\nabla'' \u^j\|_{\tilde{L}^2(B_1^+)}\sqrt{\sum_{i=2}^{n-1}(a_0^i)^2}\le cs \| \nabla(\u^j-\p_{3/2})\|_{\tilde{L}^2(B_1)}
\end{equation}
for some constant $c$. Since $s$ is arbitrary so the only two possibilities 
in (\ref{somethingbeautiful}) are that either 
$\|\nabla'' \u^j\|_{\tilde{L}^2(B_1^+)}=0$
and the Lemma is trivially true or $\sqrt{\sum_{i=2}^{n-1}(a_0^i)^2}=0$ which 
is what we claim. We may therefore assume that (\ref{forLemma81}) is true.

We may thus apply Lemma \ref{fundamentalthinobstlemma} and conclude that 
$a^1_0=a^2_0=...=a^n_0=0$. This finishes the proof of the claim.

$ $

We have therefore shown that 
$$
\v^{i,j}\to \v^i=\sum_{k=1}^\infty \q_k
$$
strongly in $L^2$. Since we have no $\q_1$ term in the sum and therefore the 
lowest homogeneity is $1$. It is easy to see that there exist an $s_\gamma$ 
such that
$$
\|\v^i\|_{\tilde{L}^2(B_{s_\gamma}^+)}\le \frac{1}{2}s^{1/2+\gamma}\|\v^i\|_{\tilde{L}^2(B_1^+)}
$$
for each $\gamma\in (0,1/2)$. The Lemma follows by strong convergence. \qed

\section{Decay of the Solution.}\label{decaysec}

In the next proposition we prove that the difference between $\u$ and $\p_{3/2}$
is small in $L^2$ norm then the difference decays geometrically. This is implies
regularity of the solutions in a standard way as will be shown later. 

\begin{prop}\label{needstobechangedagain}
Let $\u$ solve the Signorini problem in $B_1^+$ and $0<\gamma<1/2$ then 
there exists a
$\delta_\gamma>0$ such that if 
$$
\inf_{\xi\in \Pi}\|\u-\p_{3/2}^\xi \|_{\tilde{L}^2(B_1^+)}=
\|\u-\p_{3/2} \|_{\tilde{L}^2(B_1^+)} \le \delta_\gamma
$$
and
$$
\|\nabla'' \u\|_{\tilde{L}^2(B_1^+)}\le \|\u-\p_{3/2} \|_{\tilde{L}^2(B_1^+)}
$$
then there exist an $s_\gamma>0$ such that
\begin{equation}\label{Ivar}
\inf_{\xi\in \Pi}\|\u-\p_{3/2}^\xi \|_{\tilde{L}^2(B_{s_\gamma}^+)}\le 
s_{\gamma}^\gamma \|\u-\p_{3/2} \|_{\tilde{L}^2(B_1^+)}
\|\u\|_{\tilde{L}^2(B_{s_\gamma}^+)}
\end{equation}
and
\begin{equation}\label{wasprovedinthelemma}
\|\nabla'' \u\|_{\tilde{L}^2(B_{s_\gamma}^+)}\le 
s_{\gamma}^{1/2+\gamma}\|\nabla'' \u\|_{\tilde{L}^2(B_1^+)}.
\end{equation}
\end{prop}
\textsl{Proof:} Let us first point out that (\ref{wasprovedinthelemma}) was proved in 
Lemma \ref{localgainforobst}.

The proof is unfortunately quite long so we will divide it into several Lemmas.

\begin{lem}\label{Claim1}
Let $\u^j$ be a sequence of solutions as in Proposition 
\ref{needstobechangedagain} corresponding to $\delta_{\gamma,j}=\delta_j\to 0$ 
then there exist a modulus of continuity $\sigma$ such that
\begin{equation}\label{potpuri}
\inf_{\xi\in\Pi}\|\u^j-\p_{3/2}^\xi\|_{\tilde{L}^2(B_t^+)}\le \sigma(\delta_j)
\|\u^j\|_{\tilde{L}^2(B_t^+)}
\end{equation} 
for each $t<1$.
In particular, from Lemma \ref{assinpropimplyassinprop},  
for each $\u^j$ we have
$$
\limsup_{s\to 0}\frac{\inf_{\xi\in \Pi}\|\u^j-\p_{3/2}^\xi\|_{\tilde{L}^2(B_{s}^+)}}{\delta_j \|\u^j\|_{\tilde{L}^2(B_{s}^+)}}=0.
$$
\end{lem}
\textsl{Proof:} If (\ref{potpuri}) is not true then for some small 
$\mu>0$ there is a sequence $t_j$ such that
$$
\inf_{\xi\in\Pi}\frac{\|\u^j-\p_{3/2}^\xi\|_{\tilde{L}^2(B_{t_j}^+)}}{\|\u^j\|_{\tilde{L}^2(B_{t_j}^+)}}=\mu.
$$
We may assume that $t_j$ is the largest such $t$ corresponding to $\u^j$. Since $\delta_j\to 0$ it is easy
to see that $t_j\to 0$. We make the blow-up
$$
\w^j(x)=\frac{\u^j(t_j x)}{\|\u^j\|_{\tilde{L}^2(B_{t_j}^+)}}\to \w^0.
$$
From Lemma \ref{assinpropimplyassinprop} we conclude that
$$
\w^0=\p_{3/2}^\xi,
$$
for some $\xi\in \Pi$. This is a contradiction since
\begin{equation}\label{anotherstrongconvergence}
0= \inf_{\xi\in \Pi}
\frac{\|\w^0-\p_{3/2}^\xi\|_{\tilde{L}^2(B_{1}^+)}}{\|\w^0\|_{\tilde{L}^2(B_{1}^+)}}=
\lim_{j\to \infty}
\inf_{\xi\in\Pi}\frac{\|\u^j(t_j x)-\p_{3/2}^\xi\|_{\tilde{L}^2(B_{t_j}^+)}}{\|\u^j\|_{\tilde{L}^2(B_{t_j}^+)}}
=\mu.
\end{equation}
The second equality in (\ref{anotherstrongconvergence}) follows by
strong convergence of $\u^j(t_j x)/\|\u^j\|_{\tilde{L}^2(B_{t_j}^+)}$
(since $W^{1,2}$ compactly embeds into $L^2$) and of $\p_{3/2}^\xi$
(since these functions are contained in a finite dimensional subspace of $L^2$).
\qed

Before we state our next lemma let us describe the general idea of 
the proof of Proposition \ref{needstobechangedagain}.

The general idea to prove Proposition \ref{needstobechangedagain}  
is to argue by contradiction and assume that there exist $\u^j$, $\delta_j\to 0$
and $s_j\to 0$ such that $\|\u-\p_{3/2} \|_{\tilde{L}^2(B_1^+)}=\delta_j$, 
$\|\nabla'' \u\|_{\tilde{L}^2(B_1^+)}\le \delta_j$ and
\begin{equation}\label{generalprincipleforS}
\frac{\inf_{\xi\in \Pi}\|\u^j-\p_{3/2}^\xi\|_{\tilde{L}^2(B_{s_j}^+)}}{\delta_j \|\u^j\|_{\tilde{L}^2(B_{s_j}^+)}}=C_js_j^{\gamma_j}.
\end{equation}
Where $C_j\to \infty$ and $0<\gamma_j <1/2$. We will assume that the sequence 
$\gamma_j\to \gamma_0$, for some $\gamma_0$ that may be zero. However, as the 
proof will show, $\gamma_j\to 1/2$ or else we get a contradiction.

We may assume, if not the proposition is clearly true, that
$$
\frac{\inf_{\xi\in\Pi}\|\u^j-\p_{3/2}^\xi\|_{\tilde{L}^2(B_{\tilde{s}_j^+)}}}{\tilde{s}_j^{1/2}
\delta_j\|\u^j\|_{\tilde{L}^2(B_{\tilde{s}_j}^+)}}\to \infty
$$
for some sequence $\tilde{s}_j\to 0$. We also know from Lemma \ref{Claim1} that
$$
\lim_{s\to 0}\frac{\inf_{\xi\in\Pi}\|\u^j-\p_{3/2}^\xi\|_{\tilde{L}^2(B_s^+)}}{\delta_j\|\u^j\|_{\tilde{L}^2(B_s^+)}}=0
$$
for each $j$. Therefore we can choose $s_j\to 0$, $C_j\to \infty$ such that 
\begin{equation}\label{conditionsonsigmasjgammaj}
\begin{array}{ll}
\frac{\inf_{\xi\in\Pi}\|\u^j-\p_{3/2}^\xi\|_{\tilde{L}^2(B_{s_j}^+)}}{\delta_j\|\u^j\|_{\tilde{L}^2(B_{s_j}^+)}}\le C_j & \textrm{ if } s< s_j \\
\frac{\inf_{\xi\in\Pi}\|\u^j-\p_{3/2}^\xi\|_{\tilde{L}^2(B_{s_j}^+)}}{\delta_j\|\u^j\|_{\tilde{L}^2(B_{s_j}^+)}}\le C_js^{\gamma_j}  & \textrm{ if } s>s_j.
\end{array}
\end{equation}
We make the blow-up
\begin{equation}\label{letmeintroducevj}
\v^j=\frac{\u^j(s_j x)-s_j^{3/2}\p_{3/2}^{\xi_j}}{C_js_j^{\gamma_j}\delta_j\|\u^j\|_{\tilde{L}^2(B_{s_j}^+)}}.
\end{equation}
Then $\v^j\to \v^0$ weakly in $L^2(B_R^+)$ for each $R>1$ and also
$$
\|\v^j\|_{\tilde{L}^2(B_r^+)}=\frac{1}{C_j s_j^{\gamma_j}\delta_j \|\u^j\|_{\tilde{L}^2(B_{s_j}^+)}}
\|\u^j(s_j x)-s_j^{3/2}\p^{\xi_j}\|_{\tilde{L}^2(B_r^+)}\le
$$
$$
\frac{1}{C_j s_j^{\gamma_j}\delta_j \|\u^j\|_{\tilde{L}^2(B_{s_j}^+)}} \bigg(
\|\u^j(s_j x)-\p^{\xi_{s_j r}}\|_{\tilde{L}^2(B_{s_j r}^+)}+
r^{3/2}s_j^{3/2}\|\p_{3/2}^{\xi_j}-\p^{\xi_{s_j r}}\|_{\tilde{L}^2(B_{1}^+)}
\bigg)
$$
\begin{equation}\label{Icalledthis2inthehndwritt}
\le\left\{
\begin{array}{ll}
2\frac{\|\u^j\|_{\tilde{L}^2(B_{s_j r}^+)}}{\|\u^j\|_{\tilde{L}^2(B_{s_j}^+)}}+ \frac{r^{3/2}s_j^{3/2}}{C_j s_j^{\gamma_j}\delta_j}\frac{\|\p_{3/2}^{\xi_j}-\p_{3/2}^{\xi_{s_j r}}\|_{\tilde{L}^2(B_{1}^+)}}{\|\u^j\|_{\tilde{L}^2(B_{s_j}^+)}} &
\textrm{ if } r\le 1 \\
r^{\gamma_j}\frac{\|\u^j\|_{\tilde{L}^2(B_{s_j r}^+)}}{\|\u^j\|_{\tilde{L}^2(B_{s_j}^+)}}+\frac{r^{3/2}s_j^{3/2}}{C_j s_j^{\gamma_j}\delta_j}\frac{\|\p_{3/2}^{\xi_j}-\p_{3/2}^{\xi_{s_j r}}\|_{\tilde{L}^2(B_{1}^+)}}{\|\u^j\|_{\tilde{L}^2(B_{s_j}^+)}} & \textrm{ if }r>1,
\end{array}\right.
\end{equation}
The proof of Proposition \ref{needstobechangedagain} will be finished in
three steps.  First we will show, in Lemma \ref{SubClaim1} that (\ref{Icalledthis2inthehndwritt})
implies that $\|\v^0\|_{\tilde{L}^2(B_r^+)}\le C\big(r^{3/2+\gamma_0}+r^{3/2}\big)$. 
This will imply, by using Lemma \ref{eigenfinctforhalfspace} that either 
$\gamma_0\ge 1/2$ or $\v^0=0$. Then in Lemma \ref{SubClaim2} we will show
that $\v^j\to \v^0$ strongly which excludes that $\v^0=0$ that will
imply that $\gamma_0\ge 1/2$ so in particular, for each $\gamma < 1/2$ there 
has 
to be a $C_\gamma$ such that if $\u$ satisfies the conditions of the
proposition with $\delta$ small enough then
$$
\frac{\inf_{\xi\in \Pi}\|\u-\p_{3/2}^\xi\|_{\tilde{L}^2(B_{s_j}^+)}}{\delta \|\u\|_{\tilde{L}^2(B_{s}^+)}}\le C_\gamma s^{\gamma}
$$
for all $s<1$. 

\begin{lem}\label{SubClaim1}
Let $\v^j$ be as in (\ref{letmeintroducevj}), $s_j$, $\xi_j$, $c_j$
and $\delta_j$ chosen as in the discussion leading up to 
(\ref{letmeintroducevj}), in particular we assume that 
(\ref{conditionsonsigmasjgammaj}) holds. Also let $C>1$ then for each $r$ there
exist a $j_r$ such that
$$
\|\v^j\|_{\tilde{L}^2(B_r^+)}\le  
\left\{
\begin{array}{ll}
C r^{3/2} & \textrm{ if } r\le 1 \\
C r^{3/2 +\gamma_j} &\textrm{ if } r>1,
\end{array}\right.
$$
if $j>j_r$.
\end{lem}
\textsl{Proof:} We need to estimate the two terms in 
(\ref{Icalledthis2inthehndwritt}). Notice that by Lemma \ref{Claim1}
and Lemma \ref{assinpropimplyassinprop}, in particular the equations
(\ref{niostjarnor}) and (\ref{tiostjarnor}) in the proof with 
$\mu=\sigma(\delta_j)$ will imply that
\begin{equation}\label{3over2growthinthelimmit}
\frac{\|\u^j\|_{\tilde{L}^2(B_{s_j r}^+)}}{\|\u^j\|_{\tilde{L}^2(B_{s_j}^+)}}\le Cr^{3/2}
\end{equation}
when $\sigma(\delta_j)$ is small enough.

In order to estimate 
$$
\frac{r^{3/2}s_j^{3/2}}{C_j s_j^{\gamma_j}\delta_j}\frac{\|\p_{3/2}^{\xi_j}-\p_{3/2}^{\xi_{s_j r}}\|_{\tilde{L}^2(B_{1}^+)}}{\|\u^j\|_{\tilde{L}^2(B_{s_j}^+)}}
$$
we let
$$
\bar{\u}^j(x)=\frac{\u^j(s_j x)-s_j^{3/2}\p_{3/2}^{\xi_{s_j 2^k}}(x)}{\|\u^j-\p_{3/2}^{\xi_{s_j 2^k}}\|_{\tilde{L}^2(B_{2^ks_j})}}.
$$
Then 
$$
\| \bar{\u}^j\|_{\tilde{L}^2(B_{2^k})}=1.
$$
So $\bar{\u}^j\to \bar{\u}^0$ weakly in $L^2(B_{2^k})$. After a rotation we may assume that
$\xi_{s_j 2^k}/|\xi_{s_j 2^k}|=e_1$ and conclude that $\bar{\u}^j\to \bar{\u}^0$
strongly in $W^{1,2}(B_{2^k}\setminus \{\sqrt{|x_1|^2+|x_n|^2}\le t\})$
for any $t>0$. This is true since $\bar{\u}^j$ is a solution, with bounded 
$L^2-$norm, of the following Lame system
$$
\begin{array}{ll}
\Delta \bar{\u}^j+ \frac{\lambda+2}{2}\nabla \div(\bar{\u}^j)=0  &  \textrm{ in } B_{1/s_j}^+ \\
\frac{\partial \bar{\u}^j\cdot e_n}{\partial x_n}+ \frac{\lambda}{4}\div(\bar{\u}^j)=0 & \textrm{ on }\Pi \cap \{x_1< -t\} \\
\bar{\u}^j\cdot e_n=0 & \textrm{ on } \Pi \cap \{x_1 >t\} \\
\frac{\partial \bar{\u}^j\cdot e_i}{\partial x_n}+\frac{\partial \bar{\u}^j\cdot e_n}{\partial x_i}=0 & \textrm{ on }\Pi\setminus \{|x_1|<t\} \textrm{ for }i=1,2,...,n-1
\end{array}
$$
and for any $t>0$ if $j$ is large enough.

In particular, by Lemma \ref{eigenfinctforhalfspace} we may
conclude that
$$
\bar{\u}^0=\sum_{i=0}^{\infty}a_i \q_i.
$$
But $a_1=0$ since we subtracted $\p_{3/2}^{\xi_{s_j 2^k}}$ in the
definition of $\bar{\u}^j$  (see Lemma \ref{subtraction}). 

We therefore have, with the notation
$$
\inf_{\xi\in \Pi}\|\bar{\u}^j-\p_{3/2}^\xi\|_{\tilde{L}^2(B_1^+)}=
\|\bar{\u}^j-\p_{3/2}^{\bar{\xi}_j}\|_{\tilde{L}^2(B_1^+)},
$$
that
$$
\|\p_{3/2}^{\bar{\xi}_j}\|_{\tilde{L}^2(B_1^+)}=o(1)
$$
That is 
\begin{equation}\label{ompalompa}
\|\p_{3/2}^{\xi_j}-\p_{3/2}^{\xi_{s_j r}}\|_{\tilde{L}^2(B_1^+)}=s_j^{-3/2}\|\u_j-\p_{3/2}^{\xi_{s_jr}}\|_{\tilde{L}^2(B_{rs_j}^+)}\times o(1).
\end{equation}
Using equation (\ref{conditionsonsigmasjgammaj}) then (\ref{ompalompa}) 
can be estimated,
when $r>1$, by
$$
s_j^{-3/2}\|\u_j-\p_{3/2}^{\xi_{s_jr}}\|_{\tilde{L}^2(B_{rs_j}^+)}\le
$$
$$
s_j^{-3/2}C_js_j^{\gamma_j}\delta_j\|\u^j\|_{\tilde{L}^2(rs_j)}\le 
CC_js_j^{-3/2}\delta_j  s_j^{\gamma_j} r^{3/2}\|\u^j\|_{\tilde{L}^2(B_{s_j}^+)}
$$
where we used (\ref{3over2growthinthelimmit}) in the last inequality.
And when $r\le 1$ we may estimate (\ref{ompalompa}) according to
$$
s_j^{-3/2}\|\u_j-\p_{3/2}^{\xi_{s_jr}}\|_{\tilde{L}^2(B_{rs_j}^+)}\le
$$
$$
s_j^{-3/2}C_j\delta_j\|\u^j\|_{\tilde{L}^2(rs_j)}\le 
CC_js_j^{-3/2}\delta_j  r^{3/2}\|\u^j\|_{\tilde{L}^2(B_{s_j}^+)}.
$$
It follows that
$$
\frac{r^{3/2}s_j^{3/2}}{C_j s_j^{\gamma_j}\delta_j}\frac{\|\p_{3/2}^{\xi_j}-\p^{\xi_{s_j r}}\|_{\tilde{L}^2(B_{1}^+)}}{\|\u^j\|_{\tilde{L}^2(B_{s_j}^+)}}=o(1)
$$
when $r>1$. And 
$$
\frac{r^{3/2}s_j^{3/2}}{C_j\delta_j}\frac{\|\p_{3/2}^{\xi_j}-\p^{\xi_{s_j r}}\|_{\tilde{L}^2(B_{1}^+)}}{\|\u^j\|_{\tilde{L}^2(B_{s_j}^+)}}=o(1)
$$
when $r\le 1$.

Using this and (\ref{3over2growthinthelimmit}) in 
(\ref{Icalledthis2inthehndwritt}) we get
$$
\|\v^j\|_{\tilde{L}^2(B_r^+)}\le
\left\{
\begin{array}{ll}
Cr^{3/2} & \textrm{ if }r\le 1 \\
Cr^{3/2+\gamma_j} & \textrm{ if } r> 1,
\end{array}\right.
$$
when $j$ is large enough. \qed

\begin{lem}\label{SubClaim2}
Let $\v^j$ be as in (\ref{letmeintroducevj}) then
$\v^j\to \v^0$ strongly for some sub-sequence, in particular 
$\|\v^0\|_{\tilde{L}^2}=1$.

Moreover in some coordinate system $\v^0(x)=\v^0(x_1,x_n)$.
\end{lem}
\textsl{Proof:} As before we may rotate the coordinates 
so that $\xi_j/|\xi_j|=e_1$, we may need a different rotation for each $j$
but that will not affect the proof. Then, as before,
$$
\begin{array}{ll}
\Delta \v^j+\frac{2+\lambda}{2}\nabla \div(\v^j)=0 & \textrm{ in } B_{1/s_j}^+ \\
e_n\cdot \v^j = 0 & \textrm{ on } \Pi\cap \{x_1>t \}\cap B_{1/t} \\
\frac{\partial e_n\cdot \v^j}{\partial x_n}+\frac{\lambda}{4}\div(\v^j)=0  & \textrm{ on } \Pi\cap \{x_n< -t \}\cap B_{1/t} \\
\frac{\partial e_n\cdot \v^j}{\partial x_i}-\frac{\partial e_i\cdot \v^j}{\partial x_n}=0 & \textrm{ on }\Pi \textrm{ for } i=1,2,...,n-1
\end{array}
$$
for each $t>0$ provided that $j$ is large enough. So $\v^j\to \v^0$
in $W^{1,2}(B_{1/t}^+\setminus \{\sqrt{|x_1|^2+|x_n|^2}<t \})$ for
each $t>0$. The claim will follow if we can show that
$\|\v^j\|_{L^2(\{\sqrt{|x_1|^2+|x_n|^2}<t \}\cap B_R)}$ can be chosen arbitrarily small
for each $R$ by considering $t$ small enough and $j$ large enough.

Notice that by Lemma \ref{localgainforobst} and Lemma \ref{Claim1} we have 
$$
\|\nabla ''\u^j\|_{\tilde{L}^2(B_{s_\gamma}^+)}\le s_\gamma^{1/2+\gamma} 
\|\nabla ''\u^j\|_{\tilde{L}^2(B_{1}^+)}.
$$
By Lemma \ref{assinpropimplyassinprop} we have, when $\delta_j$ is small enough,
$$
\|\u^j\|_{\tilde{L}^2(B_t^+)}>t^{3/2+\epsilon}
$$
for $\epsilon>0$ being very small.
That implies that, after a small rotation of the coordinate system,
$$
\w=\frac{\u^j(s_\gamma x)}{\|\u^j\|_{\tilde{L}^2(B_{s_\gamma}^+)}}
$$
satisfies the conditions in Lemma \ref{localgainforobst}.
In particular when $\delta_j$ is small enough  we have for $\gamma>\gamma_j+1/2$
\begin{equation}\label{somecrapforxiandw}
\inf_{\xi\in \Pi}\|\w-\p_{3/2}^\xi\|_{\tilde{L}^2(B_1)}\le \delta_\gamma.
\end{equation}
Since $\|\nabla'' \w \|_{\tilde{L}^2(B_1^+)}\le s_\gamma^{1/2+\gamma}\delta_j$ and 
therefore by Lemma \ref{technicallemmaforrotest}
$$
\|\nabla''_{\xi_w}\w\|_{\tilde{L}^2(B_1^+)}\le Cs_\gamma^{1/2+\gamma}\delta_j
$$
where $\xi_w$ is the minimizer in (\ref{somecrapforxiandw}) and we used the notation
$\nabla''_{\xi_w}=\nabla -e_n\frac{\partial }{\partial x_n}-\xi_w/|\xi_w|^2(\xi_w\cdot \nabla )$.
By possibly decreasing $s_\gamma$ it follows that
$$
\|\nabla''_{\xi_w}\w\|_{\tilde{L}^2(B_1^+)}\le s_\gamma^{(1+\gamma)/2}\delta_j.
$$
Therefore the conditions of Lemma \ref{localgainforobst} are satisfied by $\w$.

Iterating this we see that 
$$
\|\nabla'' \v^j\|_{\tilde{L}^2(B_1)}\le C\frac{s_\gamma^{k(3/2+\gamma/2)}\|\nabla'' \u^j\|_{\tilde{L}^2(B_1^+)}}{C_js_\gamma^{k\gamma_j}\delta_j \|\u^j\|_{\tilde{L}^2(B_{s_\gamma^k)}}}\le
C\frac{s_\gamma^{k(1/2-\gamma_j-\epsilon)}}{C_jj}\to 0,
$$
where $k$ is chosen such that $s_\gamma^{k+1}< s_j\le s_\gamma^k$, 
$\epsilon$ is small (depending on $\delta$) and is the same constant as in 
Lemma \ref{assinpropimplyassinprop}. In particular,
since $\frac{\partial \v^j}{\partial x_i}$ solves an elliptic problem such as 
in equation (\ref{flatthinlikeequation}), this implies that 
$\nabla'' \v^j\to 0$ uniformly in $B_R^+$. This proves the last statement
in the Lemma. 

Next we we consider any $\kappa\in (0,R)$ and 
$e''\in \Pi\cap\{x_1=0\}$. The estimates on $\nabla'' \v^j$ implies that
$$
\|\v^j\|_{L^2(B_t(\kappa e''))}^2\le C\|\v^j+\sup|\nabla'' \v^j|\kappa\|_{L^2(B_t^+)}^2\le
$$
$$
Ct^{(n+3)}+C\big( \sup|\nabla'' \v^j|\big)^2\kappa^2 t^{n}.
$$
Therefore, by covering $B_R\cap \{\sqrt{|x_1|^2+|x_n|^2}\le t\}$
by $CR^{n-2}/t^{n-2}$ balls, we may conclude that
$$
\|\v^j\|_{L^2(\{\sqrt{|x_1|^2+|x_n|^2}<t \}\cap B_R)}^2\le C\big( t^5+t^2\big)R^3
$$
when $j$ is large enough. Choose $t\le \epsilon(CR^{n/2})^{-1}$ and
it follows that
$$
\|\v^j\|_{L^2(\{\sqrt{|x_1|^2+|x_n|^2}<t \}\cap B_R)}^2\le \epsilon.
$$
We may conclude that $\v^j\to \v^0$ strongly in $L^2(B_R^+)$ for every $R$.
\qed

By Lemma \ref{SubClaim2} it follows that
$\v^j\to \v^0(x_1,x_n)$ strongly in $L^2$. Since $\|\v^j\|_{\tilde{L}^2(B_1^+)}=1$
we can conclude that $\|\v^0\|_{\tilde{L}^2(B_1^+)}=1$. Moreover
by Lemma \ref{SubClaim1} we have
$$
\|\v^0\|_{\tilde{L}^2(B_r^+)}\le C 
\left\{
\begin{array}{ll}
r^{3/2} & \textrm{ if } r\le 1 \\
r^{3/2 +\gamma_j} &\textrm{ if } r>1,
\end{array}\right.
$$
which implies that $\v^0=\p_{3/2}$, however that is a contradiction to 
Lemma \ref{subtraction} since
we subtracted out the $\p_{3/2}$ part of $\u^j$ in the definition of $\v^j$.

We may therefore conclude that for each $\gamma<1/2$ there is a 
$C_\gamma$ such that
$$
\frac{\inf_{\xi\in \Pi}\|\u-\p_{3/2}^\xi\|_{\tilde{L}^2(B_{s_j}^+)}}{\delta \|\u\|_{\tilde{L}^2(B_{s}^+)}}\le C_\gamma s^{\gamma}
$$
for all $s<1$. The Proposition follows with slightly smaller $\gamma$ by 
choosing $s$ small enough.

\begin{cor}\label{refinedProp3}
Let $\u$ solve the Signorini problem in $B_1^+$ and $0<\gamma<1/4$ then there exists a
$\delta_\gamma>0$ such that if 
$$
\inf_{\xi\in\Pi}\|\u-\p_{3/2}^\xi\|_{\tilde{L}^2(B_1^+)}=\|\u-\p_{3/2}\|_{\tilde{L}^2(B_1^+)}
\le \delta_\gamma,
$$
and
$$
\|\nabla'' \u\|_{\tilde{L}^2(B_1^+)}\le \delta_\gamma
$$
then there exist an $s_\gamma$ such that
\begin{equation}\label{thirteenstar}
\max\big( s_\gamma^{-(1/2+\gamma)}\|\nabla'' \u\|_{\tilde{L}^2(B_{s_\gamma}^+)},s_\gamma^{-(3/2+\gamma)}\|\u-\p_{3/2}^\xi\|_{\tilde{L}^2(B_{s_\gamma}^+)}\big)\le 
\end{equation}
$$
\max\big( \|\nabla'' \u\|_{\tilde{L}^2(B_1^+)}, \|\u-\p_{3/2}\|_{\tilde{L}^2(B_1^+)}\big)
$$
where $\xi$ is the vector that minimizes 
$\inf_{\xi\in \Pi}\|\u-\p_{3/2}^\xi\|_{\tilde{L}^2(B_1^+)}$.
\end{cor}
\textsl{Proof:} In principle the proof consists of applying Lemma \ref{localgainforobst} and Proposition \ref{needstobechangedagain}. Unfortunately
this is not as straightforward as one might hope. We will have to
split the proof into four cases.

Properly speaking we only prove that there exist an
$\tilde{s}_\gamma\in(s_\gamma^{2n},s_\gamma)$ such that the Corollary holds where
$s_\gamma$ is the constant in Proposition \ref{needstobechangedagain}.
It is easy to see that this is implies the Corollary.

$ $

\textbf{Case 1:} {\sl If
\begin{equation}\label{GoodCaseinCor2}
\|\nabla_{\xi}'' \u\|_{\tilde{L}^2(B_1^+)}\le \|\u-\p_{3/2}\|_{\tilde{L}^2(B_1^+)}.
\end{equation}}

\textsl{Proof of the Corollary in Case 1:}
Proposition \ref{needstobechangedagain} directly implies that
$$
\max\big( s_\gamma^{-(1/2+\gamma)}\|\nabla'' \u\|_{\tilde{L}^2(B_{s_\gamma}^+)},
s_\gamma^{-(3/2+\gamma)}\inf_{\xi\in \Pi}\|\u-\p_{3/2}^\xi\|_{\tilde{L}^2(B_{s_\gamma}^+)}\big)\le 
$$
$$
\max\big( \|\nabla'' \u\|_{\tilde{L}^2(B_1^+)}, \|\u-\p_{3/2}\|_{\tilde{L}^2(B_1^+)}
\|\u\|_{\tilde{L}^2(B_{s_\gamma}^+)}\big)
$$
for $\gamma<1/2$.

By Lemma \ref{Claim1} the assumptions in Lemma \ref{assinpropimplyassinprop}
holds and we may thus deduce that
$$
\|\u\|_{\tilde{L}^2(B_{s_\gamma}^+)}\ge s_\gamma^{3/2+\epsilon}
$$
so the Corollary follows, with $\gamma-\epsilon$ in place of $\gamma$, if (\ref{GoodCaseinCor2}) is true. Since $\epsilon$ is arbitrarily small the Corollary follows
in the situation of case 1.

$ $

\textbf{Case 2:} {\sl If
\begin{equation}\label{BadCaseinCor2}
 \|\u-\p_{3/2}\|_{\tilde{L}^2(B_1^+)}\le \|\nabla_{\xi}'' \u\|_{\tilde{L}^2(B_1^+)}
\end{equation}
and
\begin{equation}\label{ThorsDog}
s_\gamma^{-1}\inf_{\xi\in \Pi}\|\u-\p_{3/2}^\xi\|_{\tilde{L}^2(B_{s_\gamma}^+)}\le
\|\nabla'' \u\|_{\tilde{L}^2(B_{s_\gamma}^+)}.
\end{equation}}

\textsl{Proof of the Corollary in Case 2:} From Lemma \ref{localgainforobst} 
we still have for $\gamma<1/4$
\begin{equation}\label{badcaseincor2}
s_\gamma^{-(1/2+2\gamma)}\|\nabla''\u\|_{\tilde{L}^2(B_{s_\gamma}^+)}\le 
\|\nabla'' \u\|_{\tilde{L}^2(B_1^+)}.
\end{equation}
Then (\ref{ThorsDog}) implies that
\begin{equation}\label{badcasecor2second}
s_\gamma^{-(1+\gamma)}\inf_{\xi\in \Pi}\|\u-\p_{3/2}^\xi\|_{\tilde{L}^2(B_{s_\gamma}^+)}
\le \|\nabla'' \u\|_{\tilde{L}^2(B_{s_\gamma}^+)} \le
s_\gamma^{1/2+\gamma}\|\nabla_{\xi}'' \u\|_{\tilde{L}^2(B_1^+)}.
\end{equation}
(\ref{badcaseincor2}) and (\ref{badcasecor2second}) implies the Corollary.

$ $

In order to state the third and the fourth case we need some notation.
First we notice that if we are not in case 1 or case 2 then
\begin{equation}\label{doaneasyproblem}
\|\u-\p_{3/2}\|_{\tilde{L}^2(B_1^+)} \le \|\nabla_{\xi}'' \u\|_{\tilde{L}^2(B_1^+)}
\end{equation}
and
\begin{equation}\label{somethingSNcoulddo}
\|\nabla'' \u\|_{\tilde{L}^2(B_{s_\gamma}^+)}
\le s_\gamma^{-1}\inf_{\xi\in \Pi}\|\u-\p_{3/2}^\xi\|_{\tilde{L}^2(B_{s_\gamma}^+)}.
\end{equation}
If (\ref{doaneasyproblem}) and (\ref{somethingSNcoulddo}) holds 
then
$$
\textrm{max}\big(
s_\gamma^{-(3/2+\gamma)}\inf_{\xi\in \Pi}\|\u-\p_{3/2}^\xi\|_{\tilde{L}^2(B_{s_\gamma}^+)},
s_\gamma^{-(1/2+\gamma)}\|\nabla'' \u\|_{\tilde{L}^2(B_{s_\gamma}^+)}
\big)
$$
$$
=s_\gamma^{-(3/2+\gamma)}\inf_{\xi\in \Pi}\|\u-\p_{3/2}^\xi\|_{\tilde{L}^2(B_{s_\gamma}^+)}.
$$
From Lemma \ref{technicallemmaforrotest} we can deduce that
$$
s_\gamma^{-1/2-2\gamma}\|\nabla_\xi'' \u\|_{\tilde{L}^2(B_{s_\gamma}^+)}\le
\big( CC_1+1\big)\|\nabla'' \u\|_{\tilde{L}^2(B_1^+)}.
$$
Or if $s_\gamma$ is small enough that
$$
s_\gamma^{-1/2-\gamma}\|\nabla_\xi'' \u\|_{\tilde{L}^2(B_{s_\gamma}^+)}\le
\|\nabla'' \u\|_{\tilde{L}^2(B_1^+)}.
$$
We thus have that 
$$
\u^1(x)=\frac{\u(s_\gamma x)}{\sup_{B_1}|\u(s_\gamma x)|}
$$
satisfies the conditions in Case 1, with 
$\delta s_\gamma^{n/2}$ in place of $\delta$.

Case 3 and Case 4, defined below
will follow from an iteration of this. In order to iterate we denote
$$
\u^j(x)=\frac{\u(s_\gamma^j x)}{\sup_{B_1}|\u(s_\gamma^j x)|}.
$$

\textbf{Case 3:} {\sl  Assume (\ref{doaneasyproblem}), 
(\ref{somethingSNcoulddo}) and that there exist a $j_0\le 2n$
such that
\begin{equation}\label{dogistofat1}
\|\nabla_{\xi_{j}}'' \u^j\|_{\tilde{L}^2(B_1^+)}\le \inf_{\xi\in \Pi} \|\u^j-\p_{3/2}^\xi\|_{\tilde{L}^2(B_1^+)}=\|\u^j-\p_{3/2}^{\xi_j}\|_{\tilde{L}^2(B_1^+)}
\end{equation}
for $j\le j_0$ and
\begin{equation}\label{darkside}
\inf_{\xi\in \Pi}\|\u^{j_0}-\p_{3/2}^\xi\|_{\tilde{L}^2(B_1^+)}=
\|\u^{j_0}-\p_{3/2}^{\xi_{j_0}}\|_{\tilde{L}^2(B_1^+)} \ge 
\|\nabla_{\xi}'' \u^{j_0}\|_{\tilde{L}^2(B_1^+)}.
\end{equation}}

\textsl{Proof of the Corollary in Case 3:} Observe that (\ref{dogistofat1})
implies that we may use Lemma \ref{localgainforobst} and 
Lemma \ref{technicallemmaforrotest} on $\u^j$ for $j<j_0$ and deduce 
that
\begin{equation}\label{darkside2}
\|\nabla''_{\xi^{j_0-1}} \u\|_{\tilde{L}^2(B_{s^{j_0}}^+)}\le
s_\gamma^{1/2+2\gamma}\|\nabla''_{\xi^{j_0-1}} \u\|_{\tilde{L}^2(B_{s^{j_0-1}}^+)}
\end{equation}
$$
\le s_\gamma^{1/2+\gamma}\|\nabla''_{\xi^{j_0-2}} \u\|_{\tilde{L}^2(B_{s^{j_0-1}}^+)}
\le s_\gamma^{2(1/2+\gamma)+\gamma}\|\nabla''_{\xi^{j_0-2}} \u\|_{\tilde{L}^2(B_{s^{j_0-2}}^+)}
$$
$$
\le
s_\gamma^{2(1/2+\gamma)}\|\nabla''_{\xi^{j_0-3}} \u\|_{\tilde{L}^2(B_{s^{j_0-2}}^+)}
\le \dots 
\le
s_\gamma^{j_0(1/2+\gamma)}\|\nabla'' \u\|_{\tilde{L}^2(B_1^+)}
$$
where $\xi^{k}$ is the minimizing vector in
$$
\inf_{\xi}\|\u^k-\p_{3/2}^k\|_{\tilde{L}^2(B_1^+)}.
$$

Next we use (\ref{darkside}) to conclude that
\begin{equation}\label{darkside3}
\max\big( \tilde{s}^{-(1/2+\gamma)}\|\nabla''_\xi \u\|_{\tilde{L}^2(B_{\tilde{s}}^+)},\tilde{s}_\gamma^{-(3/2+\gamma)}\|\u-\p_{3/2}^\xi\|_{\tilde{L}^2(B_{\tilde{s}^+})}\big)
\end{equation}
$$
\le \tilde{s}^{-(1/2+\gamma)}\|\nabla_{\xi}'' \u\|_{\tilde{L}^2(B_{\tilde{s}}^+)}
\le \|\nabla'' \u\|_{\tilde{L}^2(B_1^+)}
$$
where $\tilde{s}=s_\gamma^{j_0}$
The Corollary follows from (\ref{darkside3}) with $\tilde{s}$ in place of 
$s_\gamma$.

$ $

\textbf{Case 4:} {\sl  Assume (\ref{doaneasyproblem}), 
(\ref{somethingSNcoulddo}) and that for each $k\le 2n$
the following holds
\begin{equation}\label{dogistofat9}
\|\nabla_{\xi_k}'' \u^k\|_{\tilde{L}^2(B_1^+)}\le \|\u^k-\p_{3/2}^{\xi_k}\|_{\tilde{L}^2(B_1^+)}=\inf_{\xi\in\Pi}\|\u^k-\p_{3/2}^{\xi}\|_{\tilde{L}^2(B_1^+)}.
\end{equation}
}

\textsl{Proof of the Corollary in Case 4:} As in case 3 we have that
(\ref{dogistofat9}) implies that we may apply Lemma \ref{localgainforobst}
to $\u^k$. However $\u^k$ will also satisfy the conditions in case 1.
Applying case 1 on $\u^k$ for $k\le 2n$ we can conclude that
\begin{equation}\label{dogistoofat10}
\max\big( \tilde{s}^{-(1/2+\gamma)}\|\nabla''_\xi \u\|_{\tilde{L}^2(B_{\tilde{s}}^+)},\tilde{s}^{-(3/2+\gamma)}\|\u-\p_{3/2}^\xi\|_{\tilde{L}^2(B_{\tilde{s}}^+)} \big)\le 
\end{equation}
$$
\max\big( \|\nabla'' \u\|_{\tilde{L}^2(B_{s_\gamma}^+)}, \|\u-\p_{3/2}\|_{\tilde{L}^2(B_{s_\gamma}^+)}\big)
$$
where $\tilde{s}=s_\gamma^{2n}$. But
\begin{equation}\label{dogistoofat11}
\max\big( \|\nabla'' \u\|_{\tilde{L}^2(B_{s_\gamma}^+)}, \|\u-\p_{3/2}\|_{\tilde{L}^2(B_{s_\gamma}^+)}\big)\le
\end{equation}
$$
s_\gamma^{-n/2}\max\big( \|\nabla'' \u\|_{\tilde{L}^2(B_{1}^+)}, \|\u-\p_{3/2}\|_{\tilde{L}^2(B_{1}^+)}\big).
$$
The inequalities (\ref{dogistoofat10}) and (\ref{dogistoofat11}) implies that
$$
\max\big( \|\nabla'' \u\|_{\tilde{L}^2(B_{1}^+)}, \|\u-\p_{3/2}\|_{\tilde{L}^2(B_{1}^+)}\big)
$$
$$
\ge
s_\gamma^{n/2}\max\big( \tilde{s}^{-(1/2+\gamma)}\|\nabla''_\xi \u\|_{\tilde{L}^2(B_{\tilde{s}}^+)},\tilde{s}^{-(3/2+\gamma)}\|\u-\p_{3/2}^\xi\|_{\tilde{L}^2(B_{\tilde{s}}^+)} \big)
$$
$$
=\max\big( \tilde{s}^{-(1/2+\gamma-1/4)}\|\nabla''_\xi \u\|_{\tilde{L}^2(B_{\tilde{s}}^+)},\tilde{s}^{-(3/2+\gamma-1/4)}\|\u-\p_{3/2}^\xi\|_{\tilde{L}^2(B_{\tilde{s}}^+)} \big).
$$
But this is true for each $\gamma<1/2$ so
$$
\max\big( \tilde{s}^{-(1/2+\gamma)}\|\nabla''_\xi \u\|_{\tilde{L}^2(B_{\tilde{s}}^+)},\tilde{s}^{-(3/2+\gamma)}\|\u-\p_{3/2}^\xi\|_{\tilde{L}^2(B_{\tilde{s}}^+)} \big)
$$
$$
\le 
\max\big( \|\nabla'' \u\|_{\tilde{L}^2(B_{1}^+)}, \|\u-\p_{3/2}\|_{\tilde{L}^2(B_{1}^+)}\big)
$$
for each $\gamma< 1/4$. This finishes the proof for case 4 and the Corollary.
\qed

\section{Regularity of Solutions.}

We are now finally ready to prove that solutions are in fact $C^{1,1/2}$
which is the first main result of the paper and the main result of this section.

Before we prove the main theorem we will need one more small Lemma.

\begin{lem}\label{Claim1Theorem1}
Let $\u$ solve the Signorini problem in $B_1^+$ and 
$\|\u\|_{L^\infty(B_1^+)}=1$.
Then for every $\delta>0$ there exist a $C_\delta$ such that
if for some $r$ we have
$$
\frac{\|\u\|_{\tilde{L}^2(B_s^+)}}{C_\delta s^{3/2}}\le 1
$$
for $s\ge r$ and
$$
\frac{\|\u\|_{\tilde{L}^2(B_r^+)}}{C_\delta r^{3/2}}= 1.
$$
Then we have 
\begin{equation}\label{17star}
\inf_{\xi\in \Pi}\bigg{\|}\frac{\u(rx)}{C_\delta r^{3/2}}-\p^\xi_{3/2} \bigg{\|}_{\tilde{L}^2(B_1)}<\delta
\end{equation}
and assuming that the minimizing $\xi=|\xi|e_1$ we also have
\begin{equation}\label{17twostar}
\bigg{\|} \nabla'' \big(\frac{\u(rx)}{C_\delta r^{3/2}}\Big)\bigg{\|}_{\tilde{L}^2(B_1^+)}<\delta.
\end{equation}
\end{lem}
\textsl{Proof:} If the Lemma is not true then there exist $u^j$ and $r_j$
such that
\begin{equation}\label{controlofvjfrombaobexyz}
\frac{\|\u^j\|_{\tilde{L}^2(B_s^+)}}{j s^{3/2}}\le 1
\end{equation}
for $s\ge r_j$ and
$$
\frac{\|\u^j\|_{\tilde{L}^2(B_{r_j}^+)}}{j r_j^{3/2}}= 1.
$$
But 
\begin{equation}\label{17contradictionwithdeltastrictlypos1}
\inf_{\xi\in \Pi}\bigg{\|}\frac{\u^j(r_jx)}{j r_j^{3/2}}-\p^\xi_{3/2} \bigg{\|}_{\tilde{L}^2(B_1)}
>\delta
\end{equation}
or
\begin{equation}\label{17contradictionwithdeltastrictlypos2}
\bigg{\|} \nabla'' \frac{\u^j(r_jx)}{j r_j^{3/2}}\bigg{\|}_{\tilde{L}^2(B_1^+)}>\delta
\end{equation}
for some fixed $\delta>0$. Since $\|\u^j\|_{L^\infty(B_1^+)}=1$ we may deduce that
$r_j\to 0$. We make the blow-up
$$
\v^j(x)=\frac{\u^j(r_j x)}{jr_j^{3/2}}.
$$
Then $\v^j\to \v^0$ strongly in $W^{1,2}$ for some sub-sequence. Next we notice  
that
$\v^0=p_{3/2}^{\xi_0}$ for some $\xi_0\in \Pi$. In particular, for $R>1$, 
$\|\v^0\|_{\tilde{L}^2(B_R^+)}\le R^{3/2}$
by (\ref{controlofvjfrombaobexyz}) and $\v^0$ is a global solution to the 
Signorini problem. Arguing as in Proposition \ref{flatness} one 
readily deduces that $\v^0(x)$
is two dimensional and the assertion that $\v^0=\p_{3/2}^{\xi_0}$ follows. 
Rotating the coordinate system we may assume that $\xi_0=|\xi_0|e_1$ and obviously
$$
\inf_{\xi\in \Pi}\| \v^0-\p_{3/2}^\xi\|_{\tilde{L}^2(B_1^+)}=
\| \nabla'' \v^0\|_{\tilde{L}^2(B_1^+)}=0
$$
this together with strong convergence clearly contradicts (\ref{17contradictionwithdeltastrictlypos1}) and (\ref{17contradictionwithdeltastrictlypos2}) when $j$ is large enough. \qed

\begin{thm}\label{Coneonehalfforthin}
Let $\u$ solve the Signorini problem in $B_1^+$ then 
$$
\|\u\|_{C^{1,1/2}(B_{1/2}^+)}\le C\|\u\|_{L^2(B_1^+)}
$$
\end{thm}

\textsl{Proof:} It is enough to show the Theorem for $\|\u\|_{L^2(B_1^+)}=1$,
since we may always apply the proof to $\u/\|\u\|_{L^2(B_1^+)}$.
We will therefore assume that $\|\u\|_{L^\infty(B_1^+)}=1$ for the rest of the proof.
It is also enough to prove that
\begin{equation}\label{iwanttogohome}
\|\u\|_{\tilde{L}^2(B_r^+(x^0))}\le C  r^{3/2}
\end{equation}
for all $r\in (0,1/2)$ and $x^0\in \Gamma\cap B_{1/2}$. Once  
(\ref{iwanttogohome}) is proved we may argue as in Corollary \ref{addeddetailin6}
to show that $\u\in C^{1,1/2}$.
We will therefore assume that $0\in \Gamma$ and show that
$\|\u\|_{\tilde{L}^2(B_r^+(0))}\le C r^{3/2}$.

We choose $\delta<\delta_\gamma$ where $\gamma=1/8$ and $\delta_\gamma$ is as in
Corollary \ref{refinedProp3}. Then, by Lemma \ref{Claim1Theorem1}, there exist 
a $C_\delta$
with the properties of that Lemma. If $\u$ is as in the Theorem then either
$\|\u\|_{\tilde{L}^2(B_r^+)}\le C_\delta r^{3/2}$ for each $r\in (0,1)$ and we are done. Or there exist a largest $r$, lets call it $r_0$, such that
$$
\|\u\|_{\tilde{L}^2(B_{r_0}^+)}=C_\delta r_0^{3/2}.
$$
Consider $\v=\frac{\u(r_0 x)}{\|\u(r_0 x)\|_{\tilde{L}^2(B_{1}^+)}}$, which by 
Lemma \ref{Claim1Theorem1}
satisfies the assumptions of Corollary \ref{refinedProp3}. Using
(\ref{thirteenstar}) we see that
$$
\inf_{\xi\in \Pi}\| \v-\p_{3/2}^\xi\|_{\tilde{L}^2(B_{s_\gamma}^+)}\le 
s_\gamma^{3/2+\gamma} \delta.
$$
Rescale $\v_{s_\gamma}=\frac{\v(s_\gamma x)}{s_\gamma^{3/2}}$ and we get
$$
\inf_{\xi\in \Pi}\|\v_{s_\gamma}-\p_{3/2}^\xi\|_{\tilde{L}^2(B_1^+)}\le s_\gamma^\gamma
\delta \|\v\|_{\tilde{L}^2(B_{s_\gamma}^+)}.
$$
Also, by Corollary \ref{refinedProp3} and Lemma \ref{technicallemmaforrotest},
$$
\|\nabla''_\xi \v_{s_\gamma}\|_{\tilde{L}^2(B_1^+)}\le Cs_\gamma^\gamma \delta\le
s_\gamma^{\gamma/2}\delta,
$$
where we have decreased $s_\gamma$ so that the last inequality holds.
In particular this implies that $\v_{s_\gamma}$ satisfies the conditions
in Corollary \ref{refinedProp3} with $s_\gamma^{\gamma/2}\delta$ instead of 
$\delta$.

Next if we, as usual, let $\xi_r$ be the minimizer of
$$
\inf_{\xi\in \Pi}\|\u-\p_{3/2}^\xi\|_{\tilde{L}^2(B_r^+)}
$$
then
$$
\|\p^{\xi_{s_\gamma}}_{3/2}-\p_{3/2}^{\xi_1}\|_{\tilde{L}^2(B_1^+)}\le
\frac{1}{s_\gamma^{3/2}} \|\p^{\xi_{s_\gamma}}_{3/2}-\p_{3/2}^{\xi_1}\|_{\tilde{L}^2(B_{s_\gamma}^+)}\le
$$
$$
\frac{1}{s_\gamma^{3/2}} \|\v-\p^{\xi_{s_\gamma}}_{3/2}\|_{\tilde{L}^2(B_{s_\gamma}^+)}+
\frac{1}{s_\gamma^{3/2}} \|\v-\p^{\xi_{1}}_{3/2}\|_{\tilde{L}^2(B_{s_\gamma}^+)}\le
$$
$$
s_\gamma^{\gamma/2}+\frac{1}{s_\gamma^{(n+3)/2}} \|\v-\p^{\xi_{1}}_{3/2}\|_{\tilde{L}^2(B_{1}^+)}\le
\delta\Big(s_{\gamma}^{\gamma/2}+s_\gamma^{\gamma-n/2} \Big).
$$
Since $\v_{s_\gamma}$ also satisfies the conditions in Corollary \ref{refinedProp3}
with $s_\gamma^{\gamma/2}\delta$ for $\delta$ we may iterate this. If
$$
\v_{s_\gamma^k}=\frac{\v(s_\gamma^k x)}{s_\gamma^{3k/2}}
$$
then
$$
\big{\|} \v_{s_\gamma^k}-\p_{3/2}^{\xi_{s_\gamma^{k\gamma}}}\big{\|}_{\tilde{L}^2(B_1^+)}\le \delta s_\gamma^{k\gamma/2},
$$
$$
\|\nabla''_{\xi_{s_\gamma^k}}\v_{s_\gamma^k} \|_{\tilde{L}^2(B_1^+)}\le \delta s_\gamma^{k\gamma/2}
$$
and
$$
\| \p_{3/2}^{\xi_{s_\gamma^{k\gamma}}}-\p_{3/2}^{\xi_1}\|_{\tilde{L}^2(B_1^+)}\le
\sum_{j=1}^{k-1}
\| \p_{3/2}^{\xi_{s_\gamma^{(j+1)\gamma}}}-\p_{3/2}^{\xi_{s_\gamma^{\gamma j}}}\|_{\tilde{L}^2(B_1^+)}\le
$$
$$
\delta\Big( s_{\gamma}^{\gamma/2}+s_\gamma^{\gamma-n/2} \Big)\sum_{j=0}^{k-1}s_{\gamma}^{\gamma j/2}=
\delta\Big( s_{\gamma}^{\gamma/2}+s_\gamma^{\gamma-n/2} \Big)
\frac{1-s_\gamma^{k\gamma/2}}{1-s_\gamma^{\gamma/2}}\le \tilde{C}(\gamma)\delta.
$$
By the triangle inequality we therefore have
$$
\|\v_{s_\gamma^{k}}\|_{\tilde{L}^2(B_1^+)}\le 
\big{\|} \v_{s_\gamma^k}-\p_{3/2}^{\xi_{s_\gamma^{k\gamma}}}\big{\|}_{\tilde{L}^2(B_1^+)}+
\big{\|} \p_{3/2}^{\xi_1}-\p_{3/2}^{\xi_{s_\gamma^{k\gamma}}}\big{\|}_{\tilde{L}^2(B_1^+)}+
$$
$$
\big{\|} \p_{3/2}^{\xi_1}\big{\|}_{\tilde{L}^2(B_1^+)}\le
\delta s_\gamma^{k\gamma/2}+\tilde{C}(\gamma)\delta+
\big{\|} \p_{3/2}^{\xi_1}\big{\|}_{\tilde{L}^2(B_1^+)}.
$$
Noticing that
$$
\big{\|} \p_{3/2}^{\xi_1}\big{\|}_{\tilde{L}^2(B_1^+)}\le 2
$$
since $\|\v\|_{\tilde{L}^2(B_1^+)}\le 1$ we may deduce that
$$
\|\u\|_{\tilde{L}^2(B_{s_\gamma^k r_0}^+)}\le CC_\delta s_\gamma^{3k/2}.
$$
The theorem follows. \qed

\section{Free Boundary Regularity}\label{freeboundregsec}

In the previous section we proved that $\u\in C^{1,1/2}$. The proof was based on
the fact that if the asymptotic profile of $\u$ at a free boundary
point is $\p_{3/2}$ then the blow-up is unique. We can use exactly the
same reasoning to show that the free boundary is $C^{1,\alpha}$ close
to a point where the asymptotic profile is $\p_{3/2}$. This is done
in this section.

\begin{thm}\label{uniqueblowupofthin}
Let $\u$ solve the Signorini problem in $B_1^+$ and $0\in \Gamma$ then 
there exists a $\delta_0>0$ such that if 
$$
\inf_{\xi\in\Pi}\|\u-\p_{3/2}^\xi\|_{\tilde{L}^2(B_1^+)}=\|\u-\p_{3/2}\|_{\tilde{L}^2(B_1^+)}
\le \delta,
$$
and
$$
\|\nabla'' \u\|_{\tilde{L}^2(B_1^+)}\le \delta,
$$
and $\delta\le \delta_0$.
Then the limit
$$
\lim_{r\to 0}\frac{\u(rx)}{r^{3/2}}=\u_0
$$
exists is unique and furthermore
$$
\| \u_0-\p_{3/2}\|_{\tilde{L}^2(B_{1}^+)}\le C\delta.
$$
\end{thm}
\textsl{Proof:} Let $\u$ be as in the Theorem with $\delta_0$ small enough
then by Corollary 
\ref{refinedProp3}
$$
\inf_{\xi\in \Pi}\frac{\|\u-\p_{3/2}^\xi\|_{\tilde{L}^2(B_{s_\gamma})}}{\|\u\|_{\tilde{L}^2(B_{s_\gamma})}}\le \eta \frac{\delta}{\|\u\|_{\tilde{L}^2(B_{1})}}
$$
for small $\eta< <1$ depending only on $\delta_0$ and $n$, in particular
by choosing $\delta_0$ small enough we may make $\eta$ as small as we need.

Next we notice that 
$$
\|\p_{3/2}^{\xi_{s_\gamma}}-\p_{3/2}\|_{\tilde{L}^2(B_{s_\gamma}^+)}\le 
\|\u-\p_{3/2}^{\xi_{s_\gamma}}\|_{\tilde{L}^2(B_{s_\gamma}^+)}+
\|\u-\p_{3/2}\|_{\tilde{L}^2(B_{s_\gamma}^+)}\le
$$
$$
\eta\delta \|\u\|_{\tilde{L}^2(B_{s_\gamma})}+
\delta s_{\gamma}^{-n/2}\le 
C\eta \delta \big( s_\gamma^{3/2}+s_\gamma^{-n/2}\big).
$$
Therefore
$$
\|\p_{3/2}^{\xi_{s_\gamma}}-\p_{3/2}\|_{\tilde{L}^2(B_{1}^+)}\le 
C\eta \delta \big( 1+s_\gamma^{-(n+3)/2}\big).
$$
In particular, this together with (\ref{wasprovedinthelemma}) and 
Lemma \ref{technicallemmaforrotest} implies that 
$$
\frac{\u(s_\gamma x)}{\|\u\|_{\tilde{L}^2(B_{s_\gamma}^+)}}
$$
satisfies the conditions in the Theorem with $\delta_1\le C\eta \delta$.
If $\delta_0$ is small enough then $\eta<\frac{1}{2C}$ and we may conclude that $\delta_1\le \delta/2$.
We may thus iterate the above and deduce that
$$
\|\p_{3/2}^{\xi_{s_\gamma^k}}-\p_{3/2}\|_{\tilde{L}^2(B_1^+)}\le
\sum_{j=1}^k \|\p_{3/2}^{\xi_{s_\gamma^j}}-\p_{3/2}^{\xi_{s_\gamma^{j-1}}}\|_{\tilde{L}^2(B_1^+)}\le 
$$
$$
C\delta \big( 1+s_\gamma^{-(n+3)/2}\big)\sum_{j=1}^k2^{-j}\le
4C\big( 1+s_\gamma^{-(n+3)/2}\big)\delta.
$$
Which implies the Theorem. \qed

\begin{cor}\label{freebdryregcor}
Let us solve the Signorini problem in $B_1^+$ and assume that $0\in \Gamma$
assume furthermore that 
$$
\lim_{r\to 0}\frac{\u(rx)}{r^{3/2}}=\p_{3/2}^\xi
$$
for some vector $\xi\in \Pi$. Then there exist an $r_0>0$ such that $\Gamma \cap B_{r_0}$
is an $(n-2)-$dimensional $C^{1,\alpha}$ manifold.
\end{cor}
\textsl{Proof:} We may, by normalizing and rotating the coordinates, assume that
$$
\lim_{r\to 0}\frac{\u(rx)}{r^{3/2}}=\p_{3/2}.
$$
It therefore exist an $s$ such that 
$$
\|\u-\p_{3/2}\|_{\tilde{L}^2(B_s^+)}\le s_\gamma^{3/2}\delta_\gamma
$$
and 
$$
\|\nabla'' \u\|_{\tilde{L}^2(B_s^+)}\le s_\gamma^{1/2}\delta_\gamma
$$
where $\delta_\gamma$ is as in Corollary \ref{refinedProp3}.

Using (\ref{Ivar}) and $\|\u\|_{\tilde{L}^2(B_r^+)}\le Cr^{3/2}$, which follows from Theorem
\ref{Coneonehalfforthin}, we may induce as in Theorem \ref{Coneonehalfforthin} that
$$
\inf_{\xi\in \Pi}\|\u-\p_{3/2}^\xi\|_{\tilde{L}^2(B_r^{+})}\le Cr^{3/2+\gamma}\delta_\gamma
$$
for $r\le s$. Also, by Lemma \ref{assinpropimplyassinprop} and 
Proposition \ref{flatness}, if $x^0\in \Gamma \cap B_{r/2}$ for any $r<s$ then
$$
\lim_{t\to 0}\frac{\u(tx+x^0)}{\|\u\|_{\tilde{L}^2(B_t^+)}}=\p_{3/2}^\xi
$$
for some $\xi$ and also using Proposition \ref{needstobechangedagain} we can deduce that
$$
\inf_{\xi\in \Pi}\|\u-\p_{3/2}^\xi\|_{\tilde{L}^2(B_r^{+}(x^0))}\le
C\inf_{\xi\in \Pi}\|\u-\p_{3/2}^\xi\|_{\tilde{L}^2(B_{2r}^{+}(0))}\le Cr^{3/2+\gamma}\delta_\gamma.
$$
From this it follows that 
$$
\|\p_{3/2}^{\xi^{x^0}_r}-\p_{3/2}^{\xi_r}\|_{\tilde{L}^2(B_1^+)}\le C r^\gamma \delta_\gamma
$$
and from Theorem \ref{uniqueblowupofthin} we have
$$
\lim_{r\to 0}\frac{\u(rx+x^0)}{r^{3/2}}=\p_{3/2}^{\xi_0^{x^0}}
$$
where
$$
\|\p_{3/2}-\p_{3/2}^{\xi_0^{x^0}}\|_{\tilde{L}^2(B_1^+)}\le C |x^0|^\gamma \delta_\gamma.
$$
We have thus shown that the normal of $\Gamma$ changes in a H\"older continuous fashion
in $B_{s/2}$ which implies the Corollary. \qed

\section{Appendix 1: Proof of Lemma \ref{c1beta}.}

In this appendix we will indicate how to prove $C^{1,\beta}$
regularity of the solutions to the Signorini problem. The proof
follows the lines of the proof in the main body of the paper, but
it is significantly simpler. We will therefore only briefly
indicate some main points.

Since we do not know that the solutions are $C^{1,\beta}$ yet 
we will no longer make the ``standing assumption'' we did in 
section \ref{secsimpconv}.

We will also use the curl operator explicitly so we will only, for the sake 
of simplicity, formulate the proof in $\R^3$.

\begin{lem}\label{growthimplieslinear}
Let $\u\ne 0$ be a global solution to the Signorini problem and assume that
$$
\liminf_{r\to \infty}\frac{\ln\big( \|\u\|_{\tilde{L}^2(B_r^+)}\big)}{\ln(r)}< 3/2
$$
then $\u$ is a linear function.
\end{lem}
\textsl{Proof:} The proof is almost line for line the same 
as the proof of Lemma \ref{ReduktionOfthe System}. Following
that proof we consider $\textrm{curl}(\u)=\w$ and deduce that
$\w^3=$constant. Noticing that
$$
\sup_{B_R^+}|\w|\le C(1+R)^{\alpha-1}
$$
we may conclude that the constant is zero if $\alpha<1$. If $\alpha\ge 1$
we may without loss of generality subtract a linear function from
$\u$ such that $\w^3=0$. Equation (\ref{poordog}) follows. As in Lemma 
\ref{ReduktionOfthe System} we may conclude that (\ref{rfv2}) and 
(\ref{shitface}) holds even without the $C^{1,\beta}$ assumption.

In particular we may deduce that
\begin{equation}\label{nottobademail}
2\frac{\partial \xi}{\partial x_3}+\tau=\left\{
\begin{array}{l}
\frac{\partial \xi}{\partial x_3}=c_i\textrm{ in each component of }
\Lambda_\u \\
\tilde{c}_i \textrm{ in each component of }\Omega_\u.
\end{array}\right.
\end{equation}
It is also easy to see that 
$$
2\frac{\partial \xi}{\partial x_3}+\tau\in W^{2,2}_{loc}(\R^3_+)
$$
which by the trace Theorem implies that 
\begin{equation}\label{fromtracetheorem}
2\frac{\partial \xi}{\partial x_3}+\tau\in W^{3/2,2}_{loc}(\Pi).
\end{equation}
But (\ref{nottobademail}) implies that $\nabla \w=0$ almost everywhere
on $\Pi$ and we may therefore conclude from (\ref{fromtracetheorem})
that $2\frac{\partial \xi}{\partial x_3}+\tau$ is constant. In other 
words $c_i=c_j=\tilde{c}_k=\tilde{c}_l$ for all $i,j,k,l$.

We may conclude, as in the main body of the paper, that $\u(x)=\u(x_1,x_3)$
and that $\Gamma_\u$ contains at most one point. Linearity follows
from Lemma \ref{eigenfunctionsinR2}. \qed

The following Corollary is an easy consequence of Lemma 
\ref{growthimplieslinear}.

\begin{cor}\label{growthisone}
Let $\u\ne 0$ be a global solution to the Signorini problem and assume that
$$
\liminf_{r\to \infty}\frac{\ln\big( \|\u\|_{\tilde{L}^2(B_r^+)}\big)}{\ln(r)}\le 1
$$
then
$$
\liminf_{r\to \infty}\frac{\ln\big( \|\u\|_{\tilde{L}^2(B_r^+)}\big)}{\ln(r)}= 1.
$$
\end{cor}

\begin{defi}\label{PrDef}
We will denote the $L^2$-projection of $\u\in L^2(B_{r}^+(x^0);\R^3)$ onto the 
space $\mathcal{P}$ by $\Pr(\u,r,x^0)$. The space $\mathcal{P}$ we is the space
of affine functions $l$ satisfying
\begin{enumerate}
\item $l(x)$ is affine of the following form;
$$
l(x)=
\left[\begin{array}{l}
b_{1} \\
b_{2} \\
b_{3}
\end{array}\right]
+
\left[\begin{array}{lll}
a_{11} & a_{12} & 0 \\
a_{21} & a_{22} & 0 \\
0 & 0 & a_{33}
\end{array}\right]
\left[ \begin{array}{l}
x_1 \\ x_2 \\ x_3
\end{array}\right],
$$ 
\item\label{cond2inmathcalP} $a_{33}+\frac{\lambda}{4}(a_{11}+a_{22}+a_{33})$=0. 
\end{enumerate}
That is $\Pr(\u,r,x^0)$ is the element in $\mathcal{P}$ that satisfies
$$
\|\u-\Pr(\u,r,x^0)\|_{L^2(B_r^+(x^0))}=\inf_{\p \in \mathcal{P}}\|\u-\p\|_{L^2(B_r^+(x^0))}.
$$
When $x^0=0$ we will just write $\Pr(\u,r)$ for $\Pr(\u,r,0)$.
\end{defi}
\textbf{Remark:} Notice that the condition $a_{13}=a_{23}=a_{31}=a_{32}=0$ 
and (ii) just implies that
$l(x)$ satisfies the same boundary data as $\u$ in $\Omega_\u$.

{\sc Proof of Lemma \ref{c1beta}:} We start by proving that
if $\u$ is a solution in $B_1^+$ such that
\begin{equation}\label{normailzedsolution}
\|\u\|_{\tilde{L}^2(B_1^+)}=1
\end{equation}
then there exist an $r_0$ such that
$$
\|\u(r_j x)-\Pr(u,r_j)\|_{\tilde{L}^2(B_{r_j}^+)}\le C r^{1+\beta}
$$
for $r\le r_0$.

In particular we will show that if $\u^j$ is a sequence of solutions
satisfying (\ref{normailzedsolution}) and $0\in \Gamma$ then
\begin{equation}\label{lininfge1plus2beta}
\liminf_{j\to \infty,r_j\to 0}\frac{|\ln\big( \|\u^j(r_j x)-\Pr(u^j,r_j)\big)|\|_{\tilde{L}^2(B_{r}^+)}\big)}{|\ln(r_j)|}\ge 1+2\beta
\end{equation}
for some $\beta>0$. Once we have shown (\ref{lininfge1plus2beta})
the Lemma follows as Corollary \ref{addeddetailin6}.

We will assume the contrary that we have sequences $\u^j$ satisfying 
(\ref{normailzedsolution}) and $r_j\to 0$ such that 
$$
\lim_{j\to \infty}\frac{\ln\big( \|\u^j(r_j x)-\Pr(u^j,r_j)\|_{\tilde{L}^2(B_1^+)}\big)}{|\ln(r_j)|}=\alpha\le 1.
$$ 
The proof will progress in several steps.

\textbf{Step 1:} {\sl Arguing as we did in Lemma \ref{SubClaim1} we may find 
a sub-sequence of $\u^j$ and a sequence $r_j\to 0$ such that
$$
\v^j(x)=\frac{u^j(r_jx)-\Pr(u^j,r_j)}{\|u^j(r_jx)-\Pr(u^j,r_j)\|_{\tilde{L}^2(B_1^+)}}\to \u^0
$$
where
$$
\sup_{B_R^+}\|\u^0\|_{\tilde{L}^2(B_R^+)}\le C(1+R)^{1+\epsilon}
$$
for some small $\epsilon$.
In particular from Corollary \ref{growthisone} it follows that $\u^0$ is
linear.}

\textbf{Step 2:} {\sl The limit $\u^0$ from step 1 satisfies
$$
\u^0=\frac{1}{\big\|\frac{\lambda+4}{2\lambda}x_1+\frac{\lambda+4}{2\lambda}x_2-x_3\big\|_{\tilde{L}^2(B_1^+)}}\left[
\begin{array}{l}
\frac{\lambda+4}{2\lambda}x_1 \\ \frac{\lambda+4}{2\lambda}x_2 \\ -x_3
\end{array}
\right]\equiv \f(x).
$$}
\textsl{Proof of Step 2:} This is a simple consequence of the fact that
$\u^0$ is a linear solution satisfying $\Pr(\u^0,1)=0$.

\textbf{Step 3:} {\sl Let $\u^j$ be as in step 2 and $\mu$ some small 
constant. Then there exist a
sequence $x^j\in B_{1/2}\cap \Pi$ and a sequence of real numbers $s_j$
such that
\begin{equation}\label{kappaappend}
\inf_{\gamma\in \R}\frac{\big\| \u^j(s_jr_jx+x^j)-\Pr(\u^j,r_js_j)-\gamma \f(x)\big\|_{\tilde{L}^2(B_1^+)}}{\big\| \u^j(s_jr_jx+x^j)-\Pr(\u^j,r_js_j)\big\|_{\tilde{L}^2(B_1^+)}}=\mu.
\end{equation}}

\textsl{Proof of Step 3:} Notice that since $0\in \Gamma$ we have 
$\Lambda\cap B_{1/2}\ne 0$ so we may find a small ball $B_{\delta}(x^j)$
such that $e_3\cdot \u^j(r_jx)>0$ in $B_{\delta}\cap \Pi$.

It is not hard to show that for some sequence $t_k\to 0$ 
$$
\lim_{t_k\to 0}\frac{\u^j(t_kr_jx+x^j)-\Pr(\u^j, t_kr_j, x^j)}{\|\u^j(t_kr_jx+x^j)-\Pr(\u^j, t_kr_j, x^j)\|_{\tilde{L}^2(B_1^+)}}=\tilde{\u}^j
$$
where
$$
\begin{array}{ll}
\Delta \tilde{\u}^j+\frac{\lambda+2}{2} \nabla \div(\tilde{\u}^j)= 0 & \textrm{ in } B_1^+\\
\frac{\partial e_3\cdot \tilde{\u}^j}{\partial x_i}+\frac{\partial e_i\cdot \tilde{\u}^j}{\partial x_3}=0 & \textrm{ on } \Pi \textrm{ for }i=1,2 \\
\frac{\partial e_3\cdot \tilde{\u}^j}{\partial x_3}+\frac{\lambda}{4}\div(\tilde{\u}^j)=0 & \textrm{ on } \Pi \\
\Pr(\tilde{\u}^j,1,0)=0. &
\end{array}
$$
Standard regularity theory implies that $\tilde{\u}^j$ can be written as
a sum of polynomials and $\Pr(\tilde{\u}^j,1,0)=0$ implies that
the zeroth and first order polynomial are identically zero.

It follows that
$$
\inf_{\gamma\in \R}\frac{\big\| \u^j(t_kr_jx+x^j)-\Pr(\u^j,t_kr_j)-\gamma \f(x)\big\|_{\tilde{L}^2(B_1^+)}}{\big\| \u^j(t_kr_jx+x^j)-\Pr(\u^j,t_kr_j)\big\|_{\tilde{L}^2(B_1^+)}}\to 1
$$
as $k\to \infty$. But by step 2 we also have that
$$
\inf_{\gamma\in \R}\frac{\big\| \u^j(r_jx+x^j)-\Pr(\u^j,r_j)-\gamma \f(x)\big\|_{\tilde{L}^2(B_1^+)}}{\big\| \u^j(r_jx+x^j)-\Pr(\u^j,r_j)\big\|_{\tilde{L}^2(B_1^+)}}\to 0
$$
as $j\to \infty$. An argument of continuity shows that we may chose the
$s_j$ as claimed in the step.

\textbf{Step 4:} {\sl Let 
$$
\w^j=\frac{\u^j(r_js_j x+x^j)-\Pr(\u^j,r_js_j,x^j)}{\|\u^j(r_js_j x+x^j)-\Pr(\u^j,r_js_j,x^j)\|_{\tilde{L}^2(B_1^+)}},
$$
with $x^j$ and $s_j$ as in step 3.

Then there exists a sub-sequence such that $\w^j\to \w^0$ strongly in $W^{1,2}$
and weakly in $W^{2,2}$. Moreover if we chose the sequences appropriately
$\w^0$ will be a linear function.}

\textsl{Proof of Step 4:} It it follows 
from strong convergence and Step 3 that
\begin{equation}\label{notequaltof}
\inf_{\gamma\in \R}\frac{\big\| \w^0-\gamma e_3\f(x)\big\|_{\tilde{L}^2(B_1^+)}}{\big\| \w^0\big\|_{\tilde{L}^2(B_1^+)}}=\mu
\end{equation}
and, if we chose the $s_j$ as the largest $s\in (0,1)$ such that (\ref{kappaappend}) holds,
\begin{equation}\label{newyork}
\inf_{\gamma\in \R}\frac{\big\| \w^0-\Pr(\w^0,R)-\gamma e_3\f(x)\big\|_{\tilde{L}^2(B_R^+)}}{\big\| \w^0\big\|_{\tilde{L}^2(B_R^+)}}\le \mu,
\end{equation}
for each $R\ge 1$.

Also from the convergence it follows that
$$
\Pr(\w^0,1)=0.
$$  

Arguing as in Lemma \ref{assinpropimplyassinprop} it follows from
(\ref{newyork}) that 
\begin{equation}\label{twostarinappendix}
\lim_{R\to \infty}\frac{\ln\big(\|\w^0\|_{\tilde{L}^2(B_R^+)}\big)}{\ln\big( R\big)}
\le 1+\epsilon.
\end{equation}

Using (\ref{twostarinappendix}) and Lemma \ref{growthimplieslinear} we may 
conclude that $\w^0$ is a linear solution to the Signorini problem, as 
step 4 claims.

We have thus constructed a solution $\w^0$ such that $\Pr(\w^0,1)=0$ and
such that (\ref{notequaltof}) holds. It follows from
$\Pr(\w^0,1)=0$ and linearity that $\w^0=\gamma \f$ for some 
$\gamma$, but that contradicts (\ref{notequaltof}). Our argument of
contradiction is therefore complete and we have shown 
(\ref{lininfge1plus2beta}).

From (\ref{lininfge1plus2beta}) and similar argument as in 
Corollary \ref{addeddetailin6} we may conclude that
$$
\frac{1}{d(x^0)^{n+2+2\beta}}\int_{B_{d(x^0)}(x^0)\cap \R^3_+}\big| \u- \u(x^0)-(\nabla \u)_{r,x^0}\cdot(x-x^0)\big|^2\le C 
$$
for each ball $B_{d(x^0)}(x^0)$ with $|x^0|<1/2$, $r\le 1/4$ and
$d(x^0)=\textrm{dist}(x^0,\Gamma)$. It is now fairly
standard, using the interior regularity for the Lame system,  
to show that this implies that $\u^0\in C^{1,\beta}$.\qed

\section{Appendix 2: Sketch of the proof of Lemma \ref{eigenfinctforhalfspace}.}\label{AppendixEigenf}

In this appendix we will briefly indicate how to prove the 
eigenfunction expansion in Lemma  \ref{eigenfinctforhalfspace}. 

We will argue as in section \ref{GlobalSolPart1}.
Following the proof of Lemma \ref{ReduktionOfthe System} 
we denote $\w=\textrm{curl}(\u)$. Then $\Delta \w=0$. Noticing that a difference
quotient argument assures that $\frac{\partial^m \u}{\partial x_2^m}\in W^{1,2}$
We can conclude that 
$$
\frac{\partial u^2}{\partial x_1}=-\frac{\partial u^3}{\partial x_2}\in H^{1/2}(\Pi)
$$
and $\frac{\partial u^1}{\partial x_2}\in H^{1/2}(\Pi)$. Therefore
$$
\frac{\partial w^3}{\partial x_3}
$$
is defined on $\Pi$. In particular
$$
\begin{array}{ll}
\Delta w^3 =0 & \textrm{ in } B_1^+ \\
\frac{\partial w^3}{\partial x_3}=0 & \textrm{ on } \Pi\cap B_1. 
\end{array}
$$
By extending $w^3$ to the lower half ball by an even reflection
we see that $w^3$ may be expressed by a power series
$$
w^3(x)=\sum_{k=0}^\infty z_k(x)
$$
where $z_k$ is a homogeneous polynomial of order $k$.

We can thus find a $\xi$ and $(\chi^1,\chi^2)$ such that
$$
(u^1,u^2)=\nabla' \xi +(\chi^1,\chi^2)
$$
where
$$
\frac{\partial \chi^2}{\partial x_1}-\frac{\partial \chi^1}{\partial x_2}=
\sum_{k=0}^\infty z_k(x).
$$
We may thus express $\chi^1$ and $\chi^2$ by power series expressions
\begin{equation}\label{eigApp1}
\chi^i=\sum_{k=0}^\infty q^i_k(x) 
\end{equation}
for $i=1,2$.

If we consider the equations for $w^1$ and $w^2$ as in Lemma 
\ref{ReduktionOfthe System} we may conclude that
$$
\frac{\partial }{\partial x_1}\bigg( \Delta \frac{\partial \xi}{\partial x_3}+\frac{\lambda+2}{2}\div(\u)\bigg)=-\frac{\partial }{\partial x_3}\Delta \chi^1
$$
and
$$
\frac{\partial }{\partial x_2}\bigg( \Delta \frac{\partial \xi}{\partial x_3}+\frac{\lambda+2}{2}\div(\u)\bigg)=-\frac{\partial }{\partial x_3}\Delta \chi^2.
$$
It follows that
$$
u^3=\frac{\partial \xi}{\partial x_3}+\chi^3+\tau
$$
for some harmonic function $\tau$ and 
$$
\Delta \chi^3=\sum_{k=0}^\infty \tilde{z}_k(x)
$$
where $\tilde{z}_k$ are homogeneous polynomials of order $k$. 
We may therefore express
\begin{equation}\label{eigApp2}
\chi^3=\sum_{k=0}^\infty q_k^3(x).
\end{equation}
That is
\begin{equation}\label{eigApp3}
\u=\nabla \xi+(\chi^1, \chi^2,\chi^3)+\tau e_3
\end{equation}
where $\chi^i$ are analytic functions. Here
$$
\begin{array}{ll}
\Delta \frac{\partial \chi}{\partial x_i}+\frac{\lambda+2}{2}\frac{\partial }{\partial x_i}\div(\nabla \chi+\tau e_3)=0 & \textrm{ in } B_1^+ \textrm{ for }i=1,2 \\
\Delta \frac{\partial \chi}{\partial x_3}+\frac{\lambda+2}{2}\frac{\partial }{\partial x_3}\div(\nabla \chi+\tau e_3)=0 & \textrm{ in } B_1^+ \\
\Delta \tau =0 & \textrm{ in } B_1^+
\end{array}
$$ 
with the boundary values
\begin{equation}
\begin{array}{ll}
\frac{\partial \xi}{\partial x_3}+\tau=\chi^3=\sum_{k=0}^\infty q_k^3(x) & \textrm{ on }\{ x_1>0\}\cap \Pi \\
\frac{\partial^2 \xi}{\partial x_3^2}+\frac{\lambda}{4}\div(\nabla \xi+\tau e_3)=-\frac{\lambda}{4}\div\big( (\chi^1,\chi^2,\chi^3)\big)-\frac{\partial \chi^3}{\partial x_3} & \textrm{ on }\{x_1<0\}\cap \Pi \\
2\frac{\partial^2 \xi}{\partial x_i\partial x_3}=-\frac{\partial\chi^i}{\partial x_3}-\frac{\partial \chi^3}{\partial x_i}=-\frac{\partial}{\partial x_3}\sum_{k=0}^\infty q_k^i-\frac{\partial}{\partial x_i}\sum_{k=0}^\infty q_k^3  & \textrm{ on } \Pi \textrm{ for }i=1,2.
\end{array}
\end{equation}
We see that we can split the boundary values into their different homogeneity's
and write 
\begin{equation}\label{eigApp4}
\xi=\sum_{k=0}^\infty \xi^k+\tilde{\xi}
\end{equation}
\begin{equation}\label{eigApp5}
\tau=\sum_{k=0}^\infty \tau_k +\tilde{\tau},
\end{equation}
where $\xi_k$ and $\tau_k$ satisfy the boundary conditions
$$
\begin{array}{ll}
\frac{\partial \xi_{k+1}}{\partial x_3}+\tau_k=\chi^3=q_k^3(x) & \textrm{ on }\{ x_1>0\}\cap \Pi \\
\frac{\partial^2 \xi_{k+1}}{\partial x_3^2}+\frac{\lambda}{4}\div(\nabla \xi_{k+1}+\tau_k e_3)=-\frac{\lambda}{4}\div\big( (q^1_{k},q_k^2,q_k^3)\big)-\frac{\partial q_k^3}{\partial x_3} & \textrm{ on }\{x_1<0\}\cap \Pi \\
2\frac{\partial^2 \xi_{k+1}}{\partial x_i\partial x_3}=-\frac{\partial q^i_k}{\partial x_3}-\frac{\partial q^3_k}{\partial x_i} & \textrm{ on } \Pi \textrm{ for }i=1,2.
\end{array}
$$
It takes some calculation to see that we may chose $\xi_{k+1}$ and $\tau_k$
to be polynomials.

The functions $\tilde{\xi}$ and $\tilde{\tau}$ will satisfy 
homogeneous boundary conditions and we can therefore analyse them as
in Corollaries \ref{simplerformcor}, \ref{taubisiszeroonOmega} and \ref{NEWCOR4}
and conclude that there is a harmonic function $\bar{\tau}$
solving the boundary value problem
$$
\begin{array}{ll}
\Delta \bar{\tau}=0 & \textrm{ in }B_1^+ \\
\frac{\partial \bar{\tau}}{\partial x_3}=0 & \textrm{ on }\{ x_1>0\}\cap \Pi \\
\frac{\partial^2 \bar{\tau}}{\partial x_3^2}=0 & \textrm{ on }\{x_1<0\}\cap \Pi
\end{array}
$$
such that
\begin{equation}\label{eigApp6}
\tilde{\tau}=\frac{\partial \bar{\tau}}{\partial x_3}
\end{equation}
and
\begin{equation}\label{eigApp7}
\tilde{\xi}=\frac{\lambda+2}{2(\lambda+4)}\frac{\partial \bar{\tau}}{\partial x_3}x_3-\frac{\lambda+3}{\lambda+4}\bar{\tau}.
\end{equation}
Since $\bar{\tau}$ is a harmonic with homogeneous boundary data
we can make the expansion of $\bar{\tau}$ into a sum of homogeneous
eigenfunctions as
\begin{equation}\label{eigApp8}
\bar{\tau}=\sum_{k=0}^\infty E_{k,\bar{\tau}}(x),
\end{equation}
here $E_{k,\bar{\tau}}$ is a homogeneous
harmonic function, but not necessarily of order $k$. In general
such an eigenfunction expression of a harmonic function may contain
generalized eigenfunctions of growth $\ln(|x|)|x|^m$, but since our boundary
is so simple no such terms will appear in (\ref{eigApp8}).

Putting (\ref{eigApp1})-(\ref{eigApp8}) together we have shown that
$\u$ may be written as a sum of homogeneous functions.
$$
\u=\sum_{k=0}^\infty \textbf{H}_k(x)
$$
where $\textbf{H}_k$ is homogeneous, but not necessarily of order $k$.

The Lemma follows if we can show that each $\textbf{H}_{k}$ is homogeneous
of order $j/2$ for some $j\in \mathbb{N}$. To that end we notice that a 
difference quotient argument implies that 
\begin{equation}\label{higherregularityfor2ndderivative}
\frac{\partial^m \textbf{H}_k}{\partial x_2^m}\in W^{1,2}\textrm{ for each }m\in \mathbb{N}.
\end{equation}
Also $\frac{\partial^m \textbf{H}_k}{\partial x_2^m}$
will be homogeneous of $m$ orders less than $\textbf{H}_k$. But this together 
with \ref{higherregularityfor2ndderivative}) implies 
that if $m$ is large enough then 
$$
\frac{\partial^m \textbf{H}_k}{\partial x_2^m}=0.
$$
Therefore there exist an $m$ such that
$$
\frac{\partial^m \textbf{H}_k}{\partial x_2^m}=\textbf{L}(x_1,x_3).
$$
By the classification in Lemma \ref{eigenfunctionsinR2} it follows that
$\textbf{L}$ is homogeneous of order $j$ or $j+1/2$. The first part of the 
Lemma follows.

The second part of Lemma \ref{eigenfinctforhalfspace} is simple to prove. 
First we notice that in polar coordinates 
$\q_0=\p_{1/2}\notin W^{2,2}$ and therefore $a_0=0$ if $\w\in W^{2,2}$.
Next we notice that $q_1=ap_{3/2}+(\nu\cdot x)p_{1/2}$, 
$\nu\cdot e_1=\nu\cdot e_2=0$, 
which is only a $W^{2,2}$ function if $\nu=0$. The second part of the Lemma follows. \qed

\bibliographystyle{plain}
\bibliography{signorini.bib}
\end{document}